\documentclass[reqno,a4paper]{elsarticle}
\journal{Arxiv}

\usepackage{amssymb,amsmath,epsfig,fancyhdr,graphics,psfrag,graphicx,color,subfigure}
\definecolor{lightblue}{rgb}{0.22,0.45,0.70}
\definecolor{lightgreen}{rgb}{0.22,0.50,0.25}

\usepackage{hyperref}

\hypersetup{colorlinks=true}

\definecolor{darkred}{rgb}{0.82,0.15,0.20}
\definecolor{darkblue}{rgb}{0.82,0.15,0.12}

\numberwithin{equation}{section}
\numberwithin{figure}{section}
\newcommand\cero{\boldsymbol{0}}

\newcommand\bI{\mathbf{I}}
\newcommand\bP{\mathbf{P}}
\newcommand\bU{\mathbf{U}}
\newcommand\bV{\mathbf{V}}
\newcommand\bW{\mathbf{W}}

\newcommand\bR{\boldsymbol{R}}
	
\newcommand\beps{\boldsymbol{\varepsilon}}

\newcommand\bb{\boldsymbol{b}}
\newcommand\bk{\boldsymbol{k}}

\newcommand\nn{\boldsymbol{n}}

\newcommand\bsigma{\boldsymbol{\sigma}}
\newcommand\bPsi{\boldsymbol{\Psi}}

\newcommand{\bUpsilon}{\boldsymbol{\Upsilon}}

\newcommand\bu{\boldsymbol{u}}
\newcommand\bv{\boldsymbol{v}}
\newcommand\bw{\boldsymbol{w}}
\newcommand\bx{\boldsymbol{x}}

\newcommand\RR{\mathbb{R}}

\newcommand\bdiv{\mathop{\mathbf{div}}\nolimits}
\newcommand\vdiv{\mathop{\mathrm{div}}\nolimits}
\newcommand\bt{\boldsymbol{t}}

\newcommand\btheta{\boldsymbol{\theta}}

\newcommand{\ii}{\mathrm{i}}
\newcommand{\betaC}{\frac{\beta_3-\beta_2-(\beta_2+\beta_3)^3}{\beta_2+\beta_3}}
\newcommand{\betaDC}{\frac{D_2(\beta_3-\beta_2)-D_1(\beta_2+\beta_3)^3}{\beta_2+\beta_3}}

\newcommand{\AIVCI}{\rho\left( c_0( D_1 + D_2 ) + \frac{\kappa}{\eta} \right)}
\newcommand{\AIVCO}{\rho c_0 \beta_1\betaC}

\newcommand{\AIIICII}{\rho\left( \frac{\kappa}{\eta}(D_1 + D_2) + c_0 D_1 D_2 \right)}
\newcommand{\AIIICI}{c_0( 2\mu + \lambda) + \alpha^2 - \rho c_0 \beta_1\betaDC 
- \frac{\kappa\rho}{\eta}\beta_1\betaC}
\newcommand{\AIIICIbO}{c_0( 2\mu + \lambda) + \alpha^2 + \gamma c_0\left( (\beta_2 + \beta_3 )\theta_1^{(i)} + \frac{\beta_3}{( \beta_2 + \beta_3 )^2}\theta_2^{(i)} \right)}
\newcommand{\AIIICIbD}{c_0( 2\mu + \lambda) + \alpha^2 - \rho c_0 \beta_1(D_2 - D_1\beta_3^2) - \frac{\kappa\rho}{\eta}\beta_1(1-\beta_3^2)}
\newcommand{\AIIICIbT}{c_0( 2\mu + \lambda) + \alpha^2 + \rho c_0 \beta_1(D_2 + D_1 \beta_2^2)  
+ \frac{\kappa\rho}{\eta}\beta_1(1+\beta_2^2)}
\newcommand{\AIIICO}{\rho c_0 \beta_1^2(\beta_2+\beta_3)^2}

\newcommand{\AIICIII}{\frac{\kappa\rho}{\eta}D_1D_2}
\newcommand{\AIICII}{ \left( c_0( 2\mu + \lambda ) + \alpha^2 \right)( D_1 + D_2 ) 
+ \frac{\kappa}{\eta}( 2\mu + \lambda ) - \frac{\kappa\rho}{\eta}\beta_1\betaDC}
\newcommand{\AIICIIbO}{\left( c_0( 2\mu + \lambda ) + \alpha^2 \right)( D_1 + D_2 ) 
+ \frac{\kappa}{\eta}( 2\mu + \lambda )}
\newcommand{\AIICIIbD}{ \left( c_0( 2\mu + \lambda ) + \alpha^2 \right)( D_1 + D_2 ) 
+ \frac{\kappa}{\eta}( 2\mu + \lambda ) - \frac{\kappa\rho}{\eta}\beta_1(D_2 - D_1\beta_3^2)}
\newcommand{\AIICIIbT}{ \left( c_0( 2\mu + \lambda ) + \alpha^2 \right)( D_1 + D_2 ) 
+ \frac{\kappa}{\eta}( 2\mu + \lambda ) + \frac{\kappa\rho}{\eta}\beta_1(D_2 + D_1\beta_2^2)}
\newcommand{\AIICI}{\frac{\kappa\rho}{\eta}\beta_1^2(\beta_2+\beta_3)^2 
- \left( c_0( 2\mu + \lambda ) + \alpha^2 \right)\beta_1\betaC}
\newcommand{\AIICIbD}{\frac{\kappa\rho}{\eta}\beta_1^2\beta_3^2 
- \left( c_0( 2\mu + \lambda ) + \alpha^2 \right)\beta_1(1 - \beta_3^2)}
\newcommand{\AIICIbT}{\frac{\kappa\rho}{\eta}\beta_1^2\beta_2^2 
+ \left( c_0( 2\mu + \lambda ) + \alpha^2 \right)\beta_1(1 + \beta_2^2)}

\newcommand{\AICIII}{\frac{\kappa}{\eta}(2\mu + \lambda)(D_1+D_2) 
+ \left( c_0(2\mu+\lambda) + \alpha^2 \right)D_1D_2}
\newcommand{\AICIIIbis}{\frac{\kappa}{\eta}(2\mu + \lambda)(D_1+D_2) 
+ \left( c_0(2\mu+\lambda) + \alpha^2 \right)D_1D_2 + \frac{\gamma\kappa}{\eta}\left( (\beta_2 + \beta_3 )\theta_1^{(i)}D_2 +\frac{\beta_3}{( \beta_2 + \beta_3 )^2}\theta_2^{(i)}D_1 \right)}
\newcommand{\AICII}{\left( c_0(2\mu + \lambda) + \alpha^2 \right) \beta_1\betaDC 
+ \frac{\kappa}{\eta}( 2\mu + \lambda )\beta_1\betaC}
\newcommand{\AICIIbD}{\left( c_0(2\mu + \lambda) + \alpha^2 \right) \beta_1(D_2 - D_1\beta_3^2) 
+ \frac{\kappa}{\eta}( 2\mu + \lambda )\beta_1(1 - \beta_3^2)}
\newcommand{\AICIIbT}{\left( c_0(2\mu + \lambda) + \alpha^2 \right) \beta_1(D_2 + D_1\beta_2^2)
+ \frac{\kappa}{\eta}( 2\mu + \lambda )\beta_1(1 + \beta_2)^2}
\newcommand{\AICI}{\left( c_0 ( 2\mu + \lambda ) + \alpha^2 \right) \beta_1^2(\beta_2+\beta_3)^2}

\newcommand{\AOCIV}{D_1D_2}
\newcommand{\AOCIII}{\beta_1\betaDC}
\newcommand{\AOCII}{\beta_1^2(\beta_2+\beta_3)^2}

\newcommand{\cblue}{}

\allowdisplaybreaks

\begin{document}
	\hypersetup{
		linkcolor=lightblue,
		urlcolor=lightblue,
		citecolor=lightblue
	}
		
	\begin{frontmatter}
		
		\title{\cblue{Stability analysis for a new model of multispecies convection-diffusion-reaction in poroelastic tissue}
			\footnote{\textit{Funding:} This work has been partially supported by the London Mathematical Society through Scheme 5, Grant 51703.}}
				
		\author[geneva]{Luis Miguel De Oliveira Vilaca}
		\ead{luismiguel.deoliveiravilaca@unige.ch}
		
		\author[ci2ma,cr]{Bryan G\'omez-Vargas}
		\ead{bryan.gomezvargas@ucr.ac.cr}

		\author[iist]{Sarvesh Kumar}
		\ead{sarvesh@iist.ac.in}
		
		\author[ox,unach,monash]{Ricardo Ruiz-Baier}
		\ead{ricardo.ruizbaier@monash.edu}
		
		\author[iist]{Nitesh Verma}
		\ead{nitesh@iist.ac.in}
				
		\address[geneva]{Department of Genetics and Evolution, 
		University of Geneva, 4 Boulevard d'Yvoy, 1205 Geneva, Switzerland.}

\address[ci2ma]{CI$^{\,2}\!$MA and Departamento de Ingenier\'\i a Matem\' atica,
			Universidad de Concepci\' on, Casilla 160-C, Concepci\' on, Chile.}
			
		\address[cr]{Present address: Secci\'on de Matem\'atica, Sede de Occidente, 
		Universidad de Costa Rica, San Ram\'on, Alajuela, Costa Rica.}

		\address[iist]{Department of Mathematics, Indian Institute of Space
			Science and Technology, Trivandrum 695 547, India.}

\address[ox]{Mathematical Institute,
			University of Oxford, A. Wiles Building,
			Woodstock Road, Oxford OX2 6GG, UK.}

		\address[unach]{Universidad Adventista de Chile, Casilla 7-D, Chill\'an, Chile.}
		
		\address[monash]{Present address: School of Mathematics, Monash University, 
			9 Rainforest Walk, Clayton VIC 3800, Australia.}
				
		
\begin{abstract}
We perform the linear stability analysis for a new  \cblue{model} for poromechanical processes \cblue{with inertia} (formulated in mixed form using the solid deformation, fluid pressure, and total pressure) interacting with diffusing and reacting solutes \cblue{convected} in the medium. We find parameter regions that lead to \cblue{spatio-temporal instabilities of the coupled} system. \cblue{The} mutual dependences between deformation and diffusive patterns are of substantial relevance in the study of morphoelastic changes in biomaterials. We provide a set of computational examples in 2D and 3D (related to brain mechanobiology) that can be used to form a better understanding on how, and up to which extent, the deformations of the porous structure dictate the generation and suppression of spatial patterning dynamics, also related to the onset of mechano-chemical waves. 
\end{abstract}
		

\end{frontmatter}

	\section{\cblue{Introduction and problem statement}}
\subsection{Scope and related work}	
We propose a new model for the interaction between diffusing species and an underlying 
poroelastic structure. This work is composed by two main contributions. In \cite{verma01} we have recently explored the well-posedness of the coupled system and have addressed the stability of a mixed finite element discretisation. On the other hand, in the present  companion paper we focus more on \cblue{detailing the physical aspects of the model, on deriving a spectral linear 
stability analysis, and in providing} numerical examples dealing with growth and pattern formation, as well as with \cblue{an application} in traumatic brain injury. 

\cblue{Recent applications of poroelastic consolidation theory to the poromechanical characterisation of soft living tissues confined to the regime of infinitesimal strains include mainly 
 the formation and development of brain oedema \cite{reis19} and the importance of including pia matter \cite{stover16}. 
 Using poroelasticity to model soft tissues is of Un relevance since 
 the permeability of tissue constituents such as  collagenous membranes is typically in the orders of $10^{-14}$ to $10^{-12}$\,[m$^2$N$^{-1}$s$^{-1}$]. If one considers membranes having a thickness of a few hundred microns, then fluid exchange occurs in the range of seconds and therefore this flow can perfectly affect physiological tissue deformations due to cardiac cycle or breathing \cite{ehret17}. 
Thus we will work under the assumption that a compound of living cells forms a macroscopic linear poroelastic structure fully saturated with interstitial fluid.} 

As considered here, reaction-diffusion equations are coupled to the balances of mass and linear 
momentum of the fluid-solid mixture through convection \cblue{(by the velocity of the poroelastic solid), as well as through a modification in the reaction, which is} modulated 
by changes in volume. In turn, the solutes and the external forces drive the motion of the medium \cblue{by means of} contractile 
forces. 
Even if the present theoretical framework is motivated by examples in cell dynamics, applications sharing the same mathematical and mechano-chemical
	structure are numerous. These include the formation of inflammatory edema in the context of immune systems \cite{reis19}, oxygen diffusivity in cartilage \cite{mauk03}, contaminant
	transport \cite{arega08}, drug delivery 
in arteries \cite{calo08},  tumour
localisation and biomass growth \cite{sacco17}, 
	or chemically-controlled cell motion \cite{moee13}. In some of these phenomena, one can observe mechanically-induced transport
of the solutes. This effect occurs as the consolidation of the porous media increases
the flow of interstitial fluid which in turn contributes to the solute advective transport  \cite{royer10}.
On the other hand, the presence of chemical solutes in so-called active poroelastic materials
locally modifies morphoelastic properties \cite{rads14}, and these processes can be homogenised to obtain
macroscopic models of poroelasticity coupled with convection-reaction-diffusion equations \cite{collis17}.

Very often, these coupled models are of high-dimensions and \cblue{strongly nonlinear}, which impedes to obtain exact solutions in closed 
form. Even if a large variety of numerical methods exist for producing approximate solutions, appropriate methods 
(in the sense of being robust with respect to model parameters, being convergent, and replicating key properties of the 
underlying physico-chemical phenomena) have appeared only recently. 
As in \cite{lee19}, here we employ a mixed three-field formulation for poroelasticity, \cblue{and 
in \cite{verma01} we have carried out the theoretical analysis of convergence properties of the scheme, tailored for} the coupling with a primal formulation for the convection-diffusion system. 

\cblue{On the other hand, and apart from the question of actually solving the set of equations, carrying out a stability analysis can reveal the essential physical mechanisms of the proposed system with respect to the parameter values. This is particularly useful in the context of patterning systems to locate the parameter space where the model leads to stationary spatially unstable solutions. However, the increasing complexity of current models implies that the stability analysis is more and more analytically involved.  In the context of the present work, related studies have 
been performed on particular sub-systems such as  decoupled elasticity and diffusion \cite{moreo2010}. More recent works tend to integrate further complexity \cblue{by adding} multi-layered 
coupled systems \cite{catla2012}, incorporating domain or mechanical growth \cite{goriely2010}, the coupling between elasticity-diffusion \cite{neville06}, poroelasticity \cite{pourjafar2017,recho2019}, and also porelasticity-diffusion \cite{rads14}, which 
resembles more the idea we advocate in this work. 
The key contributions of this paper include a new three-dimensional model for the two-way coupling between poroelasticity and reaction-diffusion, the derivation and discussion of dispersion relations that indicate that the mechano-chemical feedback onsets Turing instabilities (with non-trivial wavenumber) for a range of coupling parameters, the formulation and numerical realisation of a locking-free finite element method, and a sample of numerical results including applications in brain injuries poromechanics. This work also represents an	extension with respect to recent three-field models of poroelasticity using total pressure. We demonstrate here the feasibility 
of the model and of the numerical method to reproduce a variety of coupling scenarios including pattern suppression, linear growth instability, and other morphological changes.}  
	
	The remainder of this paper is laid out as follows. The governing equations proposed in \cite{verma01} are recalled in what is left of this section. Then, in Section~\ref{sec:stability} we perform a linear stability analysis around a steady state with zero solid displacement, constant fluid pressure, and constant solute concentrations. To make the analysis as general as possible, we modify the momentum equilibrium that we presented in \cite{verma01}, now including also an acceleration term.  We use that to make some model comparisons. 
We proceed in Section~\ref{sec:FE} with recalling the scheme from the companion paper \cite{verma01}. Then we \cblue{close with} some illustrative numerical examples in 2D and 3D collected in Section~\ref{sec:results}.

\subsection{Poroelasticity of soft tissue}	
Let us consider flow of interstitial fluid through a porous medium that is subject to elastic deformations. We will consider that the process occurs in either two- or three-dimensional domains $\Omega\subset \RR^d$ with $d\in \{ 2,3\}$, and that the fluid does not enter nor leaves the body. As common in the study of flow in porous media, we \cblue{quantify the variation in 
fluid and solid states} in terms of locally averaged variables. Then, for a given time $t\in (0,t_{\mathrm{final}}]$, poromechanical quantities of interest are in this case the average displacement of the porous structure $\bu^s(t):\Omega\to \RR^d$ and the pressure head associated with the fluid flowing through the pores, $p^f(t):\Omega\to\RR$. 
We also suppose that gravitational forces have little effect in contributing to the momentum balances in comparison to other external body forces such as applied loads depending on space and time variables $\bb(t):\Omega\to \RR^d$. In the classical theory of consolidation the system allows to describe physical loading of porous layers and the change of hydraulic equilibrium in a
fluid-structure system. There, one assumes as well that the exerted stresses contain shear contributions by the solid phase whereas volumetric contributions appear from both solid and \cblue{fluid} phases (since the interstitial flow is considered governed by Darcy's law). This fact motivates the idea from \cite{oyarzua16,lee17} to introduce an auxiliary scalar unknown 
\begin{equation}
\psi = \alpha p^f - \lambda \vdiv\bu^s, \label{eq:defphi}
\end{equation}
representing the total pressure, or the volumetric part of the total Cauchy stress $\bsigma$ (specified in the constitutive equation \eqref{eq1:sigma}, below), where $\alpha$  is the so-called Biot-Willis
	consolidation (or pressure storage coupling) parameter. 
	
Denoting by $\ell(t):\Omega\to \RR$ a given volumetric fluid source \cblue{(or fluid sink, 
which is considered a datum in a system not necessarily in equilibrium)}, 
the conservation of total pore fluid content can be stated as an equation for the fluid pressure $p^f(t):\Omega\to\RR$
\begin{equation}\label{eq1:mass}
\biggl(c_0
	+\frac{\alpha^2}{\lambda}\biggr) \partial_t p^f -\frac{\alpha}{\lambda} \partial_t \psi
	- \frac{1}{\eta} \vdiv(\kappa  \nabla p^f )  = \ell,
\end{equation}
where $\kappa(\bx)$ is the permeability (or hydraulic conductivity) of the porous medium which \cblue{can be} anisotropic, $\eta$ is the constant viscosity of the
pore fluid, and $c_0$ is the constrained specific storage coefficient \cblue{(which encompasses both the porosity of the solid skeleton and the 
compressibility of the fluid or of the solid in the meso-scale)}.  

The equations of motion (balance of linear momentum and the constitutive equation relating stress and strains) consist in 
finding solid displacements $\bu^s(t):\Omega\to \RR^d$  such that 
	\begin{align}
	\bsigma & =2\mu \beps(\bu^s)- \psi \bI, \label{eq1:sigma}\\
	\rho\partial_{tt}\bu -\bdiv\bsigma & = \rho\bb, \label{eq1:momentum}
	\end{align}
	where the total pressure is defined in \eqref{eq:defphi}, 
$\beps(\bu)=\frac{1}{2}(\nabla\bu+\nabla\bu^\intercal)$ is the tensor of infinitesimal
strains, $\bI$ is the identity tensor, $\rho$ \cblue{denotes} the density of the saturated porous material, and $\mu,\lambda$ are the shear and dilation moduli associated with the
	constitutive law of the solid structure. These and all other model parameters are assumed constant, positive and bounded, except for the dilation modulus $\lambda$, which approaches infinity  for fully incompressible materials. 
\cblue{The analysis in} \cite{verma01} does not consider acceleration in the balance of linear momentum \eqref{eq1:momentum}, as one typically supposes that solid deformations are much slower than the fluid flow rate. Nevertheless we keep that term here, as we will also explore 
the influence of inertial effects in the context of linear stability analysis. 

\subsection{Macroscopic description of two-species motion}
Next we turn to the incorporation of two interacting species whose dynamics occurs \cblue{by diffusion and reaction, as well as convection by the velocity of the moving domain}. Alternatively, one could also suppose that the species are convected only by the fluid velocity (or by the filtration velocity). 
Simpler models are able to take advantage of one-dimensional geometries, or of a constant material density of the constituents (\textit{e.g.}  cells) in order to obtain closed-form expressions for the advecting velocity \cite{neville06}. Instead, here we use the transient form of the equations of motion \eqref{eq1:momentum},\eqref{eq1:sigma} to determine such velocity.

We therefore 
	consider the propagation of a generic species with concentration
	$w_1$, reacting with an additional species \cblue{having a concentration} $w_2$. The problem
	can be written as follows
	\begin{align}
	\label{eq:ADR1}
	\partial_t w_1 + \partial_t \bu^s \cdot \nabla w_1 -
	\vdiv\{ D_1(\bx)\, \nabla w_1 \}
	&=  f(w_1,w_2,\bu^s) ,\\
	\label{eq:ADR2}
	\partial_t w_2 + \partial_t \bu^s \cdot \nabla w_2 -
	\vdiv \{ D_2(\bx)\, \nabla w_2 \}
	&=  g(w_1,w_2,\bu^s),
	\end{align}
	where $D_1,D_2$ are positive definite matrices containing possibly anisotropy of self-diffusion. 
	The net reaction terms depend on parameters that account for the reproduction of species, the removal of species concentration due to reactive interactions, and the intrinsic changes from local modifications in volume (that is, how the pore microstructure evolves with deformation). 	For illustrative purposes, and as in \cite{verma01,neville06}, we can simply consider hypothetical kinetic specifications, which can also simplify the exposition of the linear stability analysis of Section~\ref{sec:stability}. We choose 
 a modification to the classical Schnackenberg model 
	\begin{align*}
	f(w_1,w_2,\bu^s) &=  \beta_1(\beta_2 - w_1 + w_1^2w_2) + \gamma\, w_1\, \partial_t\vdiv\bu^s,
	\\
	g(w_1,w_2,\bu^s) &=\beta_1(\beta_3 -w_1^2w_2) + \gamma\, w_2\, \partial_t\vdiv\bu^s,
	\end{align*}
	where $\beta_1,\beta_2,\beta_3,\gamma$ are positive rate constants.  As mentioned above,  the mechano-chemical feedback operates only by 
	 convection and \cblue{by} the last two terms defining $f,g$. 
	These terms are modulated by $\gamma>0$, and therefore they act as 
	a local source for a given species if the solid volume increases, otherwise the additional
	terms  contribute to removal of species concentration \cite{neville06}.
	
\subsection{Active stress}	
We assume that stresses are exerted by solid, by fluid, and by morphogens. Then the forces are condensed in a macroscopic balance equation for the mixture where we recall that  the solid phase is simply considered as an isotropic deformable porous medium and 
that the fluid phase only contributes volumetrically to the stress through the hydrostatic fluid pressure at the interstitium. Microscopic tension generation is here 
supposed to occur due to active stresses
	\eqref{eq1:sigma} and
	\begin{equation}\label{eq:total-stress}
	\bsigma_{\text{total}} = \bsigma + \bsigma_{\text{act}},
	\end{equation}
	where
	the active stress operates primarily on a given, constant direction
	$\bk$, and its intensity
	depends on a scalar field $r=r(w_1,w_2)$ and on a positive constant $\tau$, to be specified later on (see \textit{e.g.}  \cite{jones12})
	\begin{equation}\label{eq:active-stress}
	\bsigma_{\text{act}} = -\tau\, r \bk\otimes\bk.
	\end{equation}

\subsection{Initial and boundary conditions under different model configurations}
We employ appropriate initial data at rest
	\begin{equation*}
	w_1(0) = w_{1,0}, \quad w_2(0)=w_{2,0}, \quad \bu^s(0)= \cero, \quad \partial_t \bu^s(0)= \cero, \quad p^f(0) = 0, \quad \psi(0) = 0 \quad \text{in $\Omega\times\{0\}$.}
	\end{equation*}
Regarding boundary conditions,  the species concentrations will assume zero diffusive flux boundary conditions on the whole boundary
$$D_1(\bx)\nabla w_1\cdot \nn = 0 \quad\text{and}\quad D_2(\bx)\nabla w_2\cdot \nn=0  \qquad\qquad\text{on $\partial \Omega\times(0,t_{\text{final}}]$}.$$
For the poromechanics we adopt either Robin conditions for the deformations (mimicking the presence of supporting springs) 
and zero fluid flux everywhere 
on the boundary, 
\begin{equation}\label{eq:Robin}
[2\mu\beps(\bu^s) - \psi\,\bI + \bsigma_{\text{act}}]\nn + \zeta \bu^s = \cero \quad \text{and}  \quad \frac{\kappa}{\eta} \nabla p^f \cdot\nn = 0\qquad \qquad\text{on $\partial \Omega\times(0,t_{\text{final}}]$},\end{equation}
where $\zeta>0$ is the (possibly time-dependent) stiffness of the spring; 
or as in \cite{verma01} we can separate the boundary  $\partial\Omega = \Gamma\cup\Sigma$  into two parts $\Gamma$ and $\Sigma$
	where we prescribe clamped boundaries and zero fluid normal fluxes; and zero (total) traction together with constant fluid pressure, respectively
	\begin{align}\label{bc:Gamma}
	\bu^s = \cero\quad \text{and} \quad \frac{\kappa}{\eta} \nabla p^f \cdot\nn = 0\qquad \qquad&\text{on $\Gamma\times(0,t_{\text{final}}]$},\\
	\label{bc:Sigma}
	[2\mu\beps(\bu^s) - \psi\,\bI + \bsigma_{\text{act}}]\nn = \cero \quad\text{and}\quad p^f=0\qquad\qquad  &\text{on $\Sigma\times(0,t_{\text{final}}]$}.
	\end{align}
Each case will be specified in the tests of Sections~\ref{sec:stability} and \ref{sec:results}. 

	\section{Linear stability analysis and dispersion relation}\label{sec:stability}
Next we proceed to derive a linear stability analysis \cblue{following} \cite{neville06}. This analysis gives insight about interaction mechanisms between tissue deformation and diffusing solutes. The present development is \cblue{however more} involved, since we are including acceleration effects in the 
 momentum equilibrium equation. We skip as much as possible the lengthy details of the derivation, and concentrate only in the dispersion relation. Note that since the functions $f$ and $g$ are prescribed, we can 
specify a steady state given by $w_1 = w_{1,0} = \beta_2 + \beta_3$, $w_2 = w_{2,0} = \frac{\beta_3}{(\beta_2+\beta_3)^2}$, $p=p_0$, $\psi=\psi_0$ and $\bu=\cero$. 
\cblue{We will restrict} the analysis \cblue{to the case of an} infinite domain in $\mathbb{R}^d$, with $d=\{2,3\}$. We also maintain the dimensional form of the governing equations 
so that the analysis accommodates a large class of models.

\subsection{General form of the dispersion relation}
Following \textit{e.g.}  \cite{zak}, we can derive a dispersion relation 
that is eventually defined by the product of two distinct polynomials 
\begin{equation*}
P(\phi;k^2) = P_1(\phi;k^2)^{d-1}P_2(\phi;k^2), 
\end{equation*}
where $P_1(\phi;k^2)=\rho\phi^2+\mu k^2$, and where $d=\{2,3\}$ is the spatial 
dimension of the infinite domain $\Omega = \mathbb{R}^d$, 
where  the linear stability analysis of the coupled problem \eqref{eq1:mass}-\eqref{eq1:momentum} is performed. 

Since $P_1$ is a polynomial with pure imaginary roots, it does not have an influence on the stability of the steady state version 
of \eqref{eq1:mass}-\eqref{eq1:momentum}. Consequently, we can focus our attention on \cblue{the}  fifth-order polynomial 
\begin{equation}
P_2(\phi;k^2) = A_5(k^2)\phi^5 + A_4(k^2)\phi^4 + A_3(k^2)\phi^3 + A_2(k^2)\phi^2 + A_1(k^2)\phi + A_0(k^2), \label{eq:la-poly}
\end{equation}
defined by the terms
{\small 
\begin{align*}
A_5(k^2) &= \rho c_0, \\
A_4(k^2) &= \AIVCI  k^2 - \AIVCO, \\
A_3(k^2) &=  \AIIICII k^4 \\
&\quad + \left[ \AIIICI \right] k^2 \\
&\quad + \AIIICO - \ii\gamma\left[\sum_{j=1}^d \widehat{\Upsilon}_j k_j \right] c_0 (w_{1,0}\theta_1 + w_{2,0}\theta_2), \\
A_2(k^2) &= \AIICIII k^6 + \left[ \AIICII \right] k^4 \\
&\quad + \left[ \AIICI \right] k^2 \\
&\quad - \ii\gamma\left[\sum_{j=1}^d \widehat{\Upsilon}_j k_j \right]\Bigg[\left( \frac{\kappa}{\eta}(w_{1,0}\theta_1 + w_{2,0}\theta_2) + c_0(w_{1,0}\theta_1 D_2 + w_{2,0}\theta_2 D_1) \right)k^2 \\
&\qquad + c_0\bigg(-w_{1,0}\theta_2\frac{2\beta_1\beta_3}{\beta_2+\beta_3} + w_{2,0}\theta_1\beta_1(\beta_2+\beta_3)^2 + w_{1,0}\theta_1\beta_1(\beta_2+\beta_3)^2 - w_{2,0}\theta_2\frac{\beta_1(\beta_3-\beta_2)}{\beta_2+\beta_3}\bigg)\Bigg], \\
A_1(k^2) &= \left[ \AICIII \right] k^6 \\
&\quad - \left[ \AICII \right] k^4 \\
&\quad + \AICI k^2 - \ii\frac{\kappa\gamma}{\eta}\left[\sum_{j=1}^d \widehat{\Upsilon}_j k_j \right]\Bigg[ (w_{1,0}\theta_1 D_2 + w_{2,0}\theta_2 D_1)k^2 \\
&\qquad -w_{1,0}\theta_2\frac{2\beta_1\beta_3}{\beta_2+\beta_3} + w_{2,0}\theta_1\beta_1(\beta_2+\beta_3)^2 + w_{1,0}\theta_1\beta_1(\beta_2+\beta_3)^2 - w_{2,0}\theta_2\frac{\beta_1(\beta_3-\beta_2)}{\beta_2+\beta_3} \Bigg] k^2, \\
A_0(k^2) &= \frac{\kappa}{\eta}(2\mu+\lambda)k^4\left( \AOCIV k^4 - \AOCIII k^2 + \AOCII \right),
\end{align*}
}
where the coefficients $w_{j,0}$, for $j=\{1,2\}$, are the steady state concentrations of 
generic species $w_j$, 
and $\theta_j = (\partial_{\bw}\sigma_{act}(\bw_0))_j$ with $\sigma_{act}(\bw) = -\tau r(\bw)$; $\widehat{\Upsilon}_j = \Upsilon_{\cdot j} + \ii\widetilde{\Upsilon}_j$, with $\Upsilon_{\cdot j} = \sum_k \partial_{\bx_k}\bUpsilon_{kj}$, $\widetilde{\Upsilon}_j = \sum_k k_k \bUpsilon_{kj}$ \cblue{where} $\bUpsilon=\bk\otimes\bk$. 

For the rest of the linear analysis, we impose that $\bUpsilon=\bI$, with $\bI$ the identity matrix, $r^{(1)}(\bw) = w_1+w_2$, $r^{(2)}(\bw) = w_1^2$ and $\bb=\cero$. \cblue{Under such conditions}, only the coefficients $A_3,A_2,A_1$ are modified and they adopt the following forms
{\small\begin{align*}
A_3(k^2) &= \left[ \AIIICII \right] k^4 \\
&\quad + \Bigg[ \AIIICI \\
&\qquad + \gamma c_0\left( (\beta_2 + \beta_3 )\theta_1^{(i)} + \frac{\beta_3}{( \beta_2 + \beta_3 )^2}\theta_2^{(i)} \right) \Bigg] k^2 + \AIIICO,  \\
A_2(k^2) &= \AIICIII k^6 + \Bigg[ \AIICII \\
&\qquad + \gamma \left( \frac{\kappa}{\eta}\left( (\beta_2 + \beta_3 )\theta_1^{(i)} +\frac{\beta_3}{( \beta_2 + \beta_3 )^2}\theta_2^{(i)} \right) + c_0\left( (\beta_2 + \beta_3 )\theta_1^{(i)}D_2 +\frac{\beta_3}{( \beta_2 + \beta_3 )^2}\theta_2^{(i)}D_1 \right) \right) \Bigg] k^4 \\
&\quad + \Bigg[ \AIICI \\
&\qquad + \gamma\left( -2\beta_1\beta_3\theta_2^{(i)} + \beta_1\beta_3\theta_1^{(i)} + \beta_1( \beta_2 + \beta_3 )^3\theta_1^{(i)} - \frac{\beta_1\beta_3( \beta_3 - \beta_2 )}{( \beta_2 + \beta_3 )^3}\theta_2^{(i)} \right)\Bigg] k^2, \\
A_1(k^2) &= \Bigg[ \AICIII + \frac{\gamma\kappa}{\eta}\left( (\beta_2 + \beta_3 )\theta_1^{(i)}D_2 +\frac{\beta_3}{( \beta_2 + \beta_3 )^2}\theta_2^{(i)}D_1 \right) \Bigg] k^6 \\
&\quad - \Bigg[ \frac{\gamma\kappa}{\eta}\left( 2\beta_1\beta_3\theta_2^{(i)} - \beta_1\beta_3\theta_1^{(i)} - \beta_1( \beta_2 + \beta_3 )^3\theta_1^{(i)} + \frac{\beta_1\beta_3( \beta_3 - \beta_2 )}{( \beta_2 + \beta_3 )^3}\theta_2^{(i)} \right) \\
&\qquad + \AICII \Bigg] k^4 \\
&\quad + \AICI k^2.
\end{align*}
}
As the characteristic polynomial \eqref{eq:la-poly} is of high order, it  is challenging to determine analytically the main features of \cblue{the coupled set of equations. We will therefore 
solve the eigenvalue systems numerically.} 
\cblue{Before that, we note that for the case without inertia} ($\rho=0$), the polynomial $P(\phi;k^2)$ is only of order 3.  We concentrate on \cblue{separate cases including} or not the inertial term. Unless specified otherwise, throughout the analysis we will employ the following parameter values 
\begin{gather*}
D_1=0.05, \quad D_2=1.0, \quad \beta_1=170, \quad \beta_2=0.1305, \quad \beta_3=0.7695, \quad E=3\cdot10^4, \quad \nu=0.495, \\
 \rho=1, \quad c_0=1 \cdot 10^{-3}, \quad \kappa = 1 \cdot 10^{-4}, \quad \alpha=0.1, \quad \eta=1, \quad \gamma=1 \cdot 10^{-4}, \quad \cblue{\ell = 0},
\end{gather*}
which are relevant to the specifications in Tests 1-4 from Section~\ref{sec:results}. 
All the computations and graphs in the remainder of this section have been produced with 
an in-house MATLAB implementation. 

\subsection{Spatial homogeneous distributions}
For the case $k^2=0$, the characteristic polynomial $P_2(\phi;0)=0$ reduces to 
\begin{equation*}
P_2(\phi;0) = \phi^3\Bigg[ \rho c_0 \phi^2 - \AIVCO \phi + \AIIICO \Bigg].
\end{equation*}
Therefore its roots are either zero, or are defined by the second-order polynomial in square brackets. Owing to the
Routh-Hurwitz conditions, 
for any polynomial of order 2, a necessary and sufficient set of conditions 
can be stated so that the roots are in the space of complex non-positive real values \cblue{$\{z\in\mathbb{C}: \Re({z}) \leq 0 \}$}. For a general polynomial $P(\phi) = a_2\phi^2 + a_1\phi + a_0$, we need to satisfy that all $a_i>0$ (or all $a_i<0$). In \cblue{the present case, $a_2$ and  $a_0$} are positive by definition, and consequently the spatial homogeneous case is stable if and only if 
\begin{equation}
\beta_3 - \beta_2 < (\beta_2+\beta_3)^3. \label{eq:cond-hom}
\end{equation}
\cblue{Thus, the difference between the basal source rate $\beta_3$ with respect to $\beta_2$ might be smaller than $(\beta_2+\beta_3)^3$}. A similar condition is provided in \cite{liu13}. Additionally, we observe that the system is homogeneously stable irrespective of  the parameter values, by simply imposing that $\rho=0$, \textit{i.e.}, removing the acceleration term in the momentum equilibrium.

\begin{figure}[t!]
\begin{center}
\includegraphics[width=0.45\textwidth]{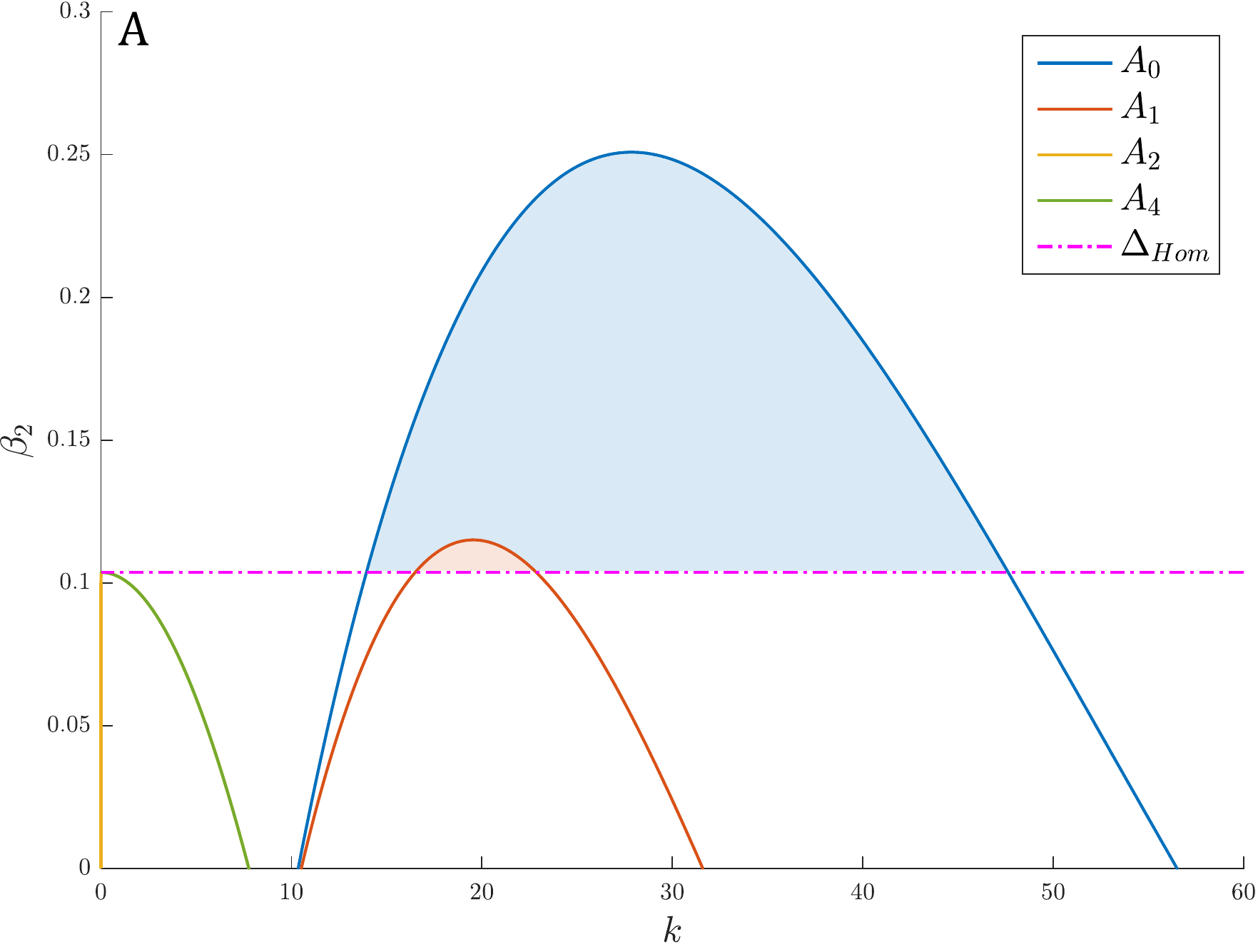}
\includegraphics[width=0.45\textwidth]{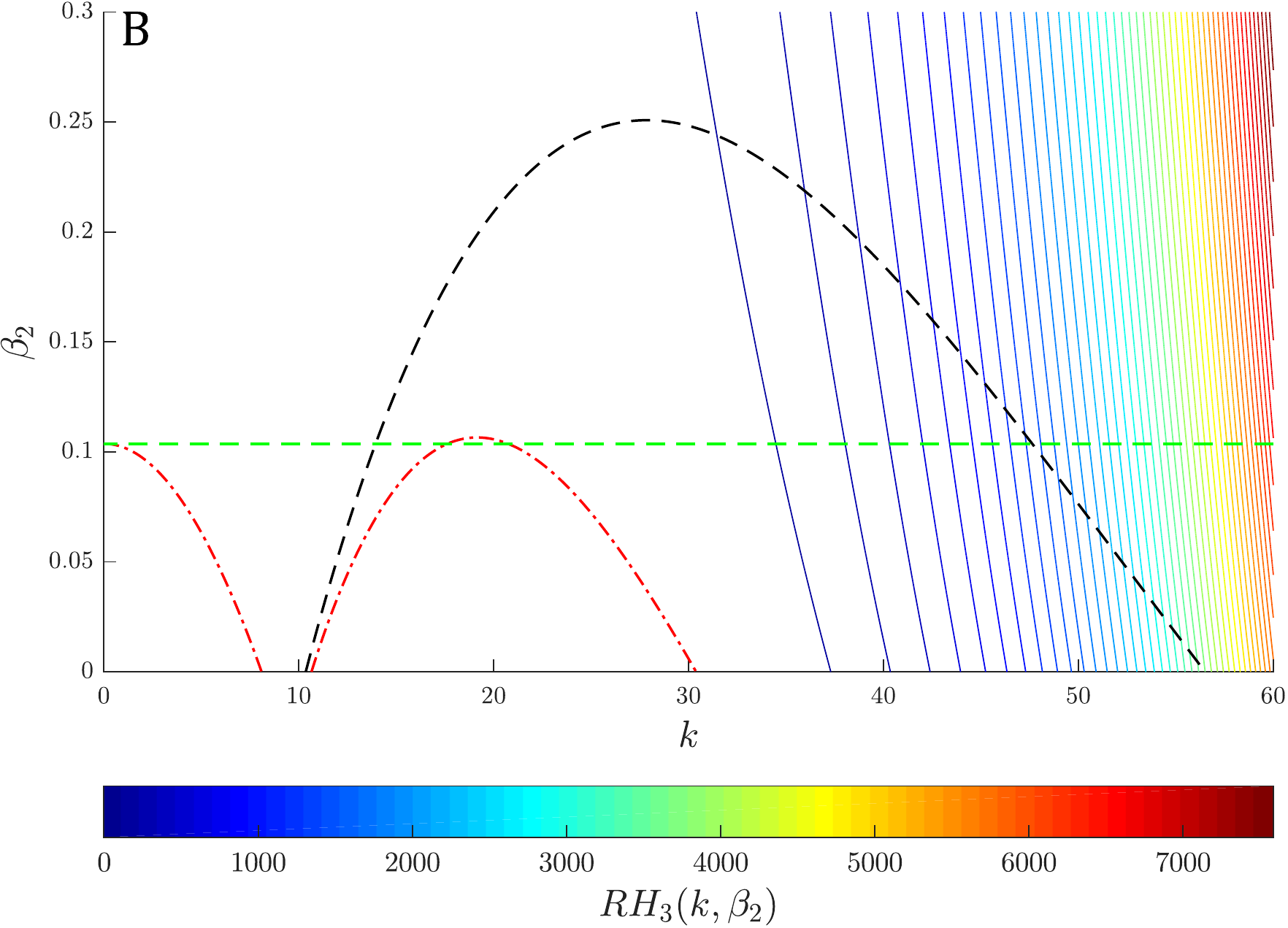}\\
\includegraphics[width=0.45\textwidth]{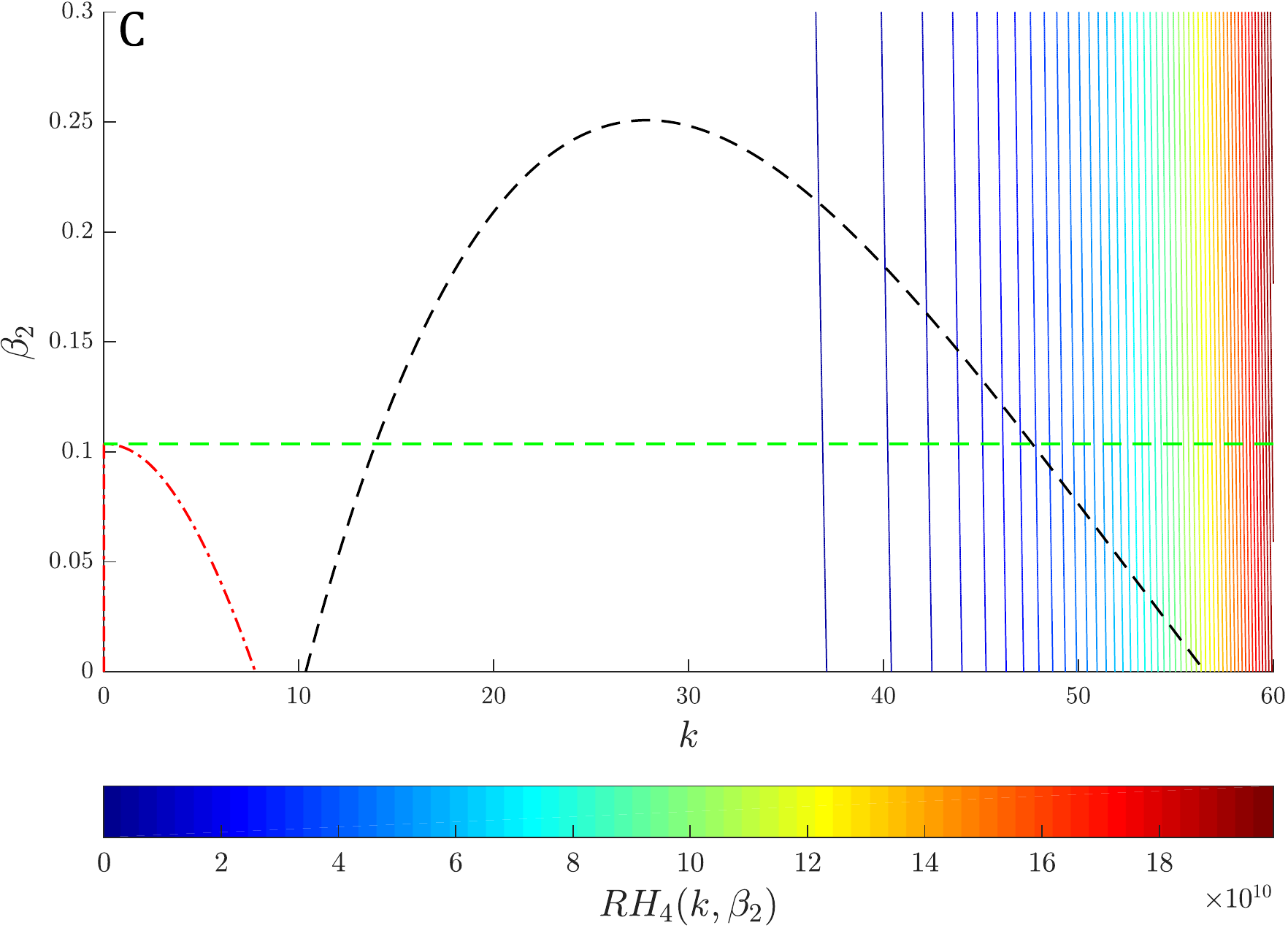}
\includegraphics[width=0.45\textwidth]{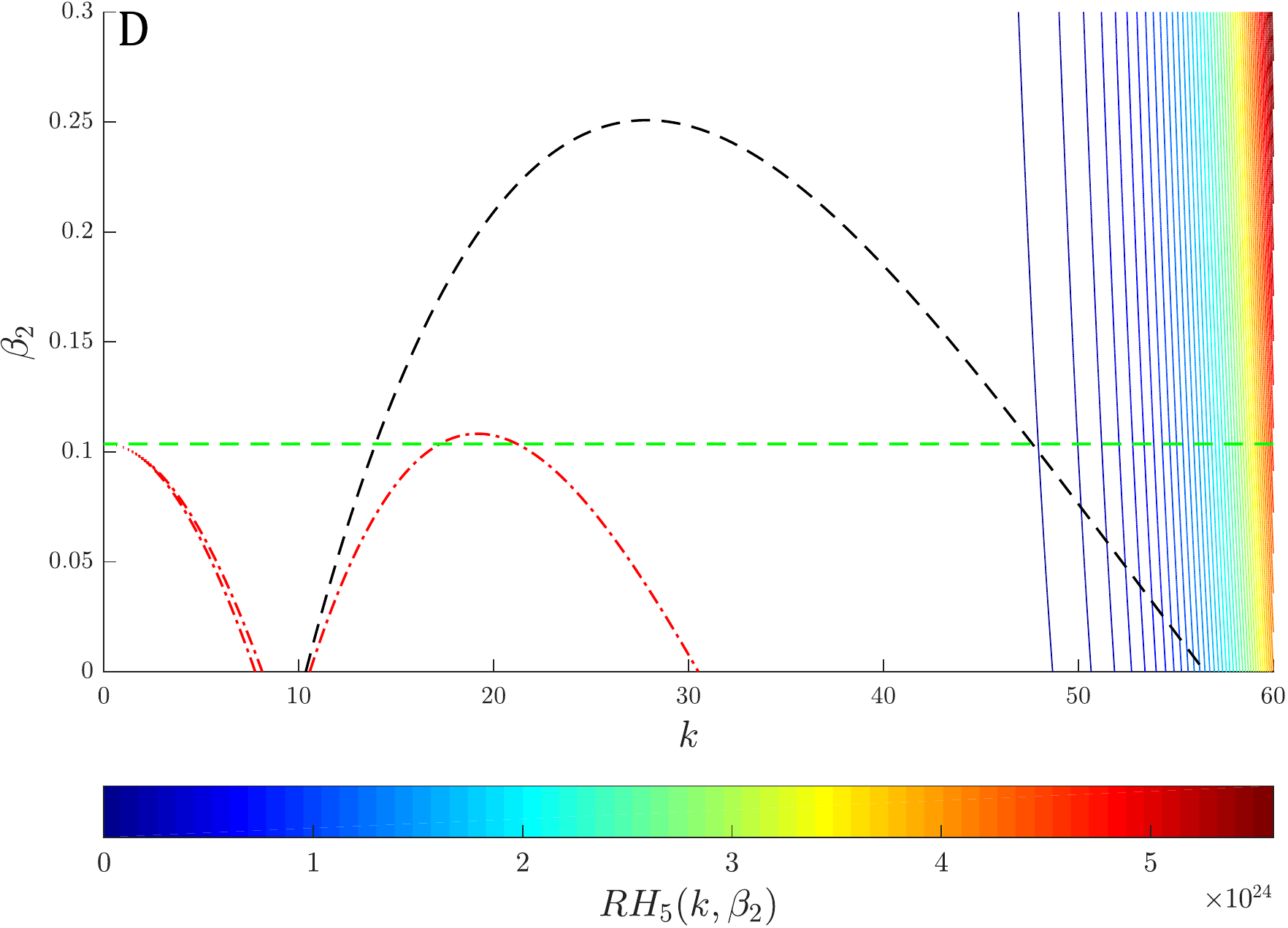}
\end{center}
\caption{Contour plots of the Routh-Hurwitz conditions for the characteristic polynomial of uncoupled system \eqref{eq1:mass}-\eqref{eq1:momentum} with $\rho=1$. (A) Null level set of $P_2$-polynomial coefficients $a_i$ defined in \eqref{eq:la-poly}. $a_3$ is strictly positive for the selected fixed parameters and so it is not presented. The magenta dot-dashed curve corresponds to the limit of condition \eqref{eq:cond-hom} with fixed $\beta_3=0.7695$. (B) Level sets (100) of condition \eqref{eq:rh-cond2}. Null levels of condition \eqref{eq:rh-cond2} (red dot-dashed);  $a_0$ (black dashed); and condition \eqref{eq:cond-hom} (green dashed) are added to locate Turing instability regions. (C-D) Similar analyses for the conditions \eqref{eq:rh-cond3} and \eqref{eq:rh-cond4} respectively. In all plots, colour-shaded regions correspond to a $k$-$\beta_2$ space that breaks the associated condition.}
\label{fig:lin_ex1}
\end{figure}

\subsection{Uncoupled system}
This scenario is reached if either $\gamma$ or $\tau$ (or both) are zero. $P_2(\phi;k^2)$ is then a fifth-order polynomial defined as in \eqref{eq:la-poly} where the terms including $\gamma$ or $\tau$ are dropped from the coefficients $A_i(k^2)$. For a polynomial of order 5, $P(\phi) = \sum_{j=0}^5 a_j \phi^j$, the Routh-Hurwitz conditions are given by
\begin{align}
\forall_j \qquad  a_j &> 0,  \label{eq:rh-cond1}  \\ 
a_3 a_4 - a_2 a_5 &> 0, \label{eq:rh-cond2} \\  
a_2 a_3 a_4 - a_2^2 a_5 - a_1 a_4^2 + a_0 a_4 a_5 &> 0, \label{eq:rh-cond3}  \\
a_0 a_2 a_3 a_4 a_5 - a_0 a_3^2 a_4^2 + a_1 a_2 a_3 a_4^2 - a_1 a_2^2 a_4 a_5 - a_1^2 a_4^3 + 2 a_0 a_1 a_4^2 a_5 - a_0^2 a_4 a_5^2 &> 0. \label{eq:rh-cond4}
\end{align}
\cblue{Condition \eqref{eq:cond-hom} indicates that $\beta_2$ and $\beta_3$ are both relevant. We decide} to perform the analysis varying $\beta_2$ and fixing all the other parameters. \cblue{The choice is justified by the influence of this parameter on the stability of the spatial homogeneous case and the direct relation of $\beta_2$ (or alternatively, $\beta_3$) on the sign of the conditions \eqref{eq:rh-cond1}-\eqref{eq:rh-cond4}}. Based on \eqref{eq:cond-hom} and the constrain on the parameters, we can be \cblue{readily deduce that $a_4$ and $a_5$ are strictly positive, and proceed to reject them since they do not lead to patterning in the system. 
The complete analysis of \eqref{eq:rh-cond1}-\eqref{eq:rh-cond4} is analytically quite involved, however some information can already be drawn by looking at 
%
conditions that violate} $a_0>0$. Beyond the tractability of the analysis, the choice of this coefficient is justified since $a_1, a_2, a_3$ can be written as affine functions of $a_0$ with positive coefficients, as long as $\rho \neq 0$. This leads, in the uncoupled system, to the property that if $a_i>0$ conditions are violated, then $a_0>0$ is inevitably unsatisfied.

\begin{figure}[!t]
\begin{center}
\includegraphics[width=0.45\textwidth]{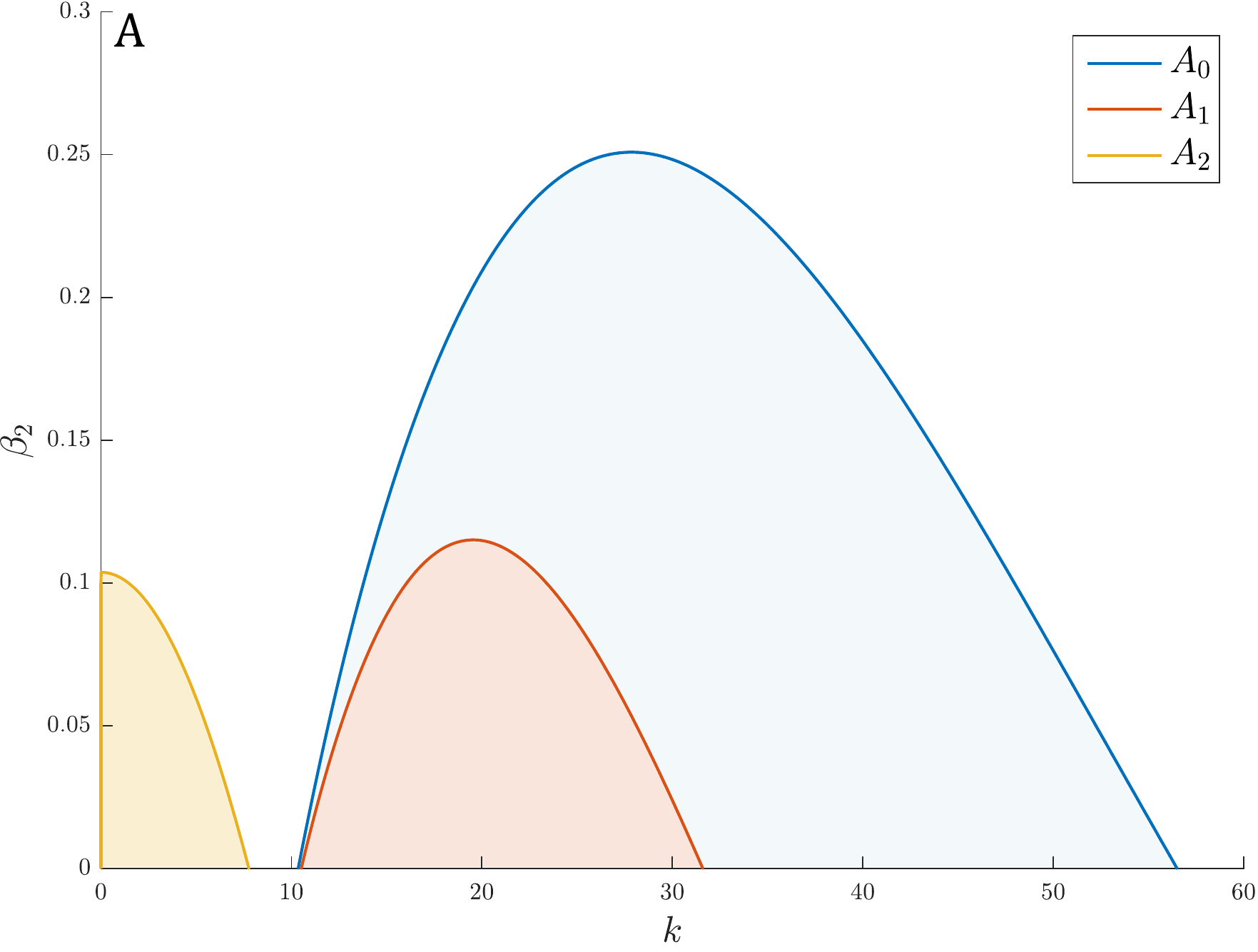}
\includegraphics[width=0.45\textwidth]{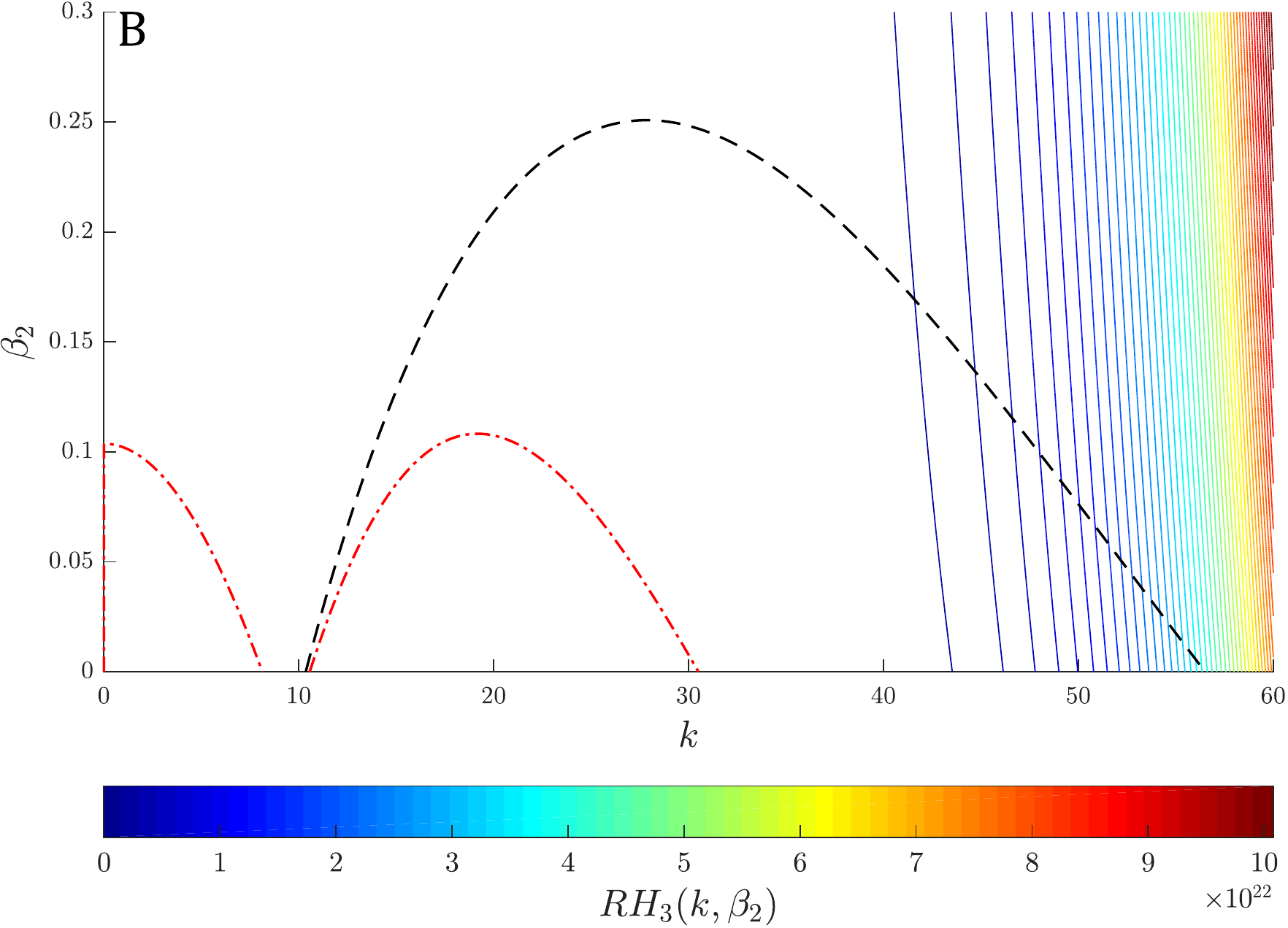}
\end{center}
\caption{Contour plots of the Routh-Hurwitz conditions for the characteristic polynomial of uncoupled system \eqref{eq1:mass}-\eqref{eq1:momentum} with $\rho=0$. (A) Null level set of $P_2$-polynomial coefficients $a_i$ defined in \eqref{eq:la-poly}. $a_3$ is strictly positive whatever the choice of the parameters and so it is not presented. (B) Level sets (100) of $2^{\text{nd}}$ condition defined in \eqref{eq:rh2-cond}. Null levels of condition \eqref{eq:rh2-cond} (red dot-dashed); and $a_0$ (black dashed) are added to locate Turing instability regions.}
\label{fig:lin_ex2}
\end{figure}

Figure~\ref{fig:lin_ex1}(A) plots the contour lines of $a_i$, $i=0,\ldots,4$. As we can observe for this specific parameter set, $a_2,a_4$ are negative only in the region below the magenta dot-dashed curve, corresponding to the limit \cblue{given by} condition \eqref{eq:cond-hom}. Consequently, these coefficients are strictly positive while \eqref{eq:cond-hom} is true. Only $a_0$ and $a_1$ present Turing instability, and the latter is just a subset of the $k$-$\beta_2$ space defined by $a_0$. Conditions \eqref{eq:rh-cond2}-\eqref{eq:rh-cond4} present a similar behaviour \cblue{as that observed for $a_0$ (Figs.~\ref{fig:lin_ex1}(B)-(D)), and therefore 
the analysis can be focused on} $a_0$ only. For the case $\rho=0$, the resulting polynomial is of order 3 and so the Routh-Hurwitz conditions are defined by
\begin{equation}
\forall_j \qquad  a_j > 0 , \qquad a_1 a_2 - a_0 a_3 > 0. \label{eq:rh2-cond} 
\end{equation}
Note that, contrary to the general case, only $a_1$ can be written as an affine function of $a_0$. As the condition \eqref{eq:cond-hom} is no longer  needed to obtain homogeneous stability, the coefficients of the affine description of $a_1$ \cblue{are not necessarily} strictly positive as before. Nevertheless, as illustrated in Figure~\ref{fig:lin_ex2}, condition $a_0>0$ is the first to be broken with respect to the value of $\beta_2$, and thus we proceed to analyse that coefficient as well.

\begin{figure}[!t]
\begin{center}
\includegraphics[width=0.325\textwidth]{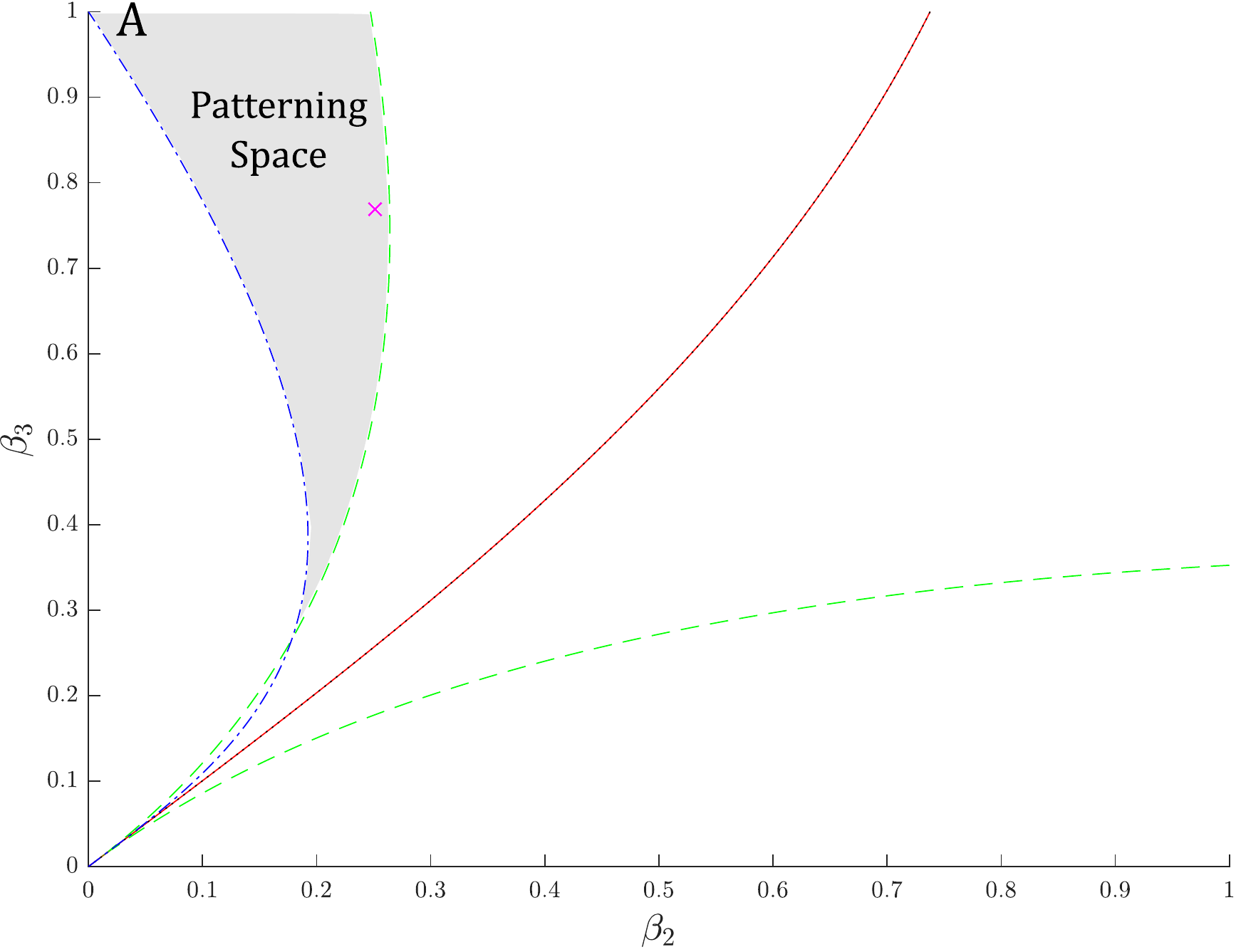}
\includegraphics[width=0.325\textwidth]{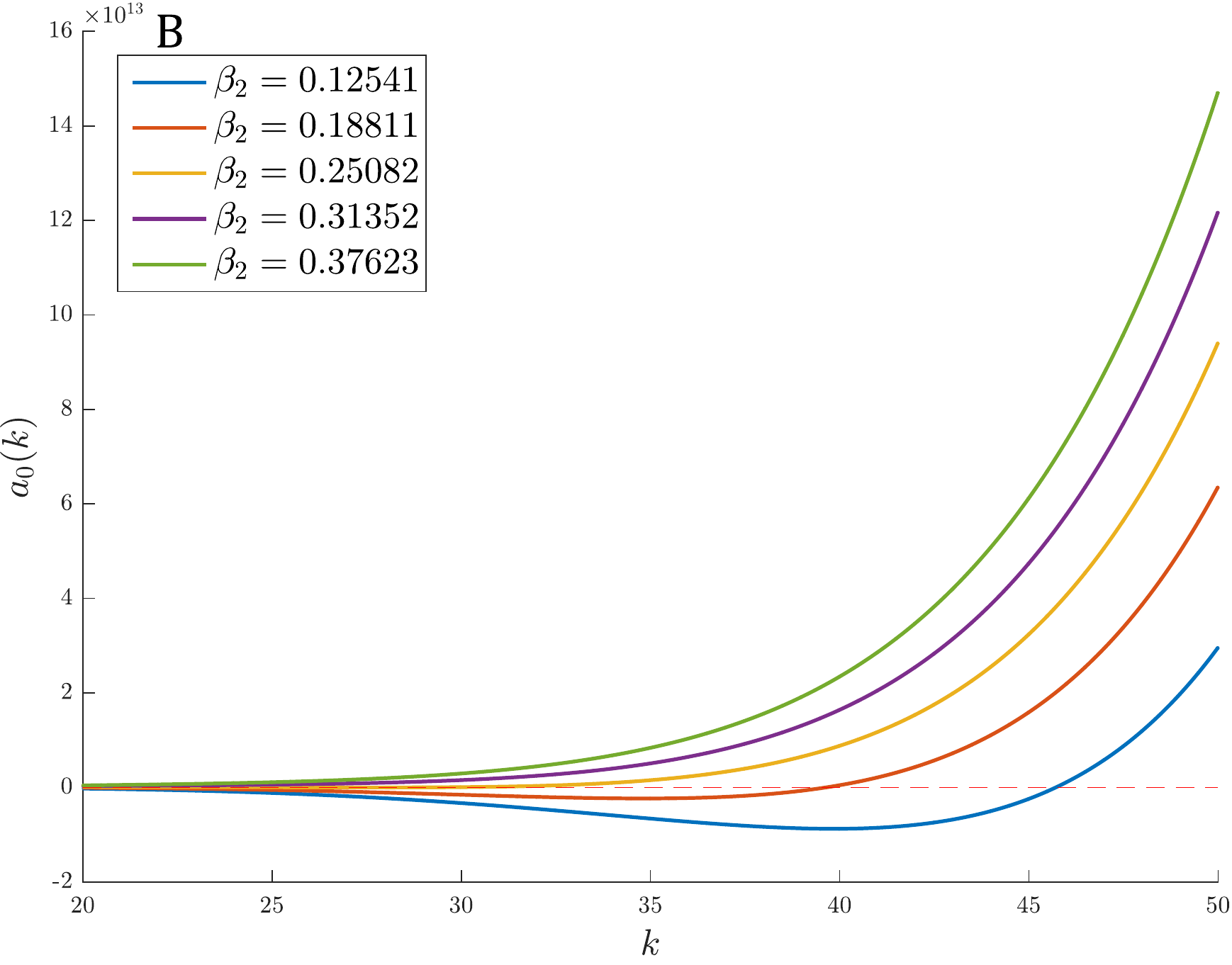}
\includegraphics[width=0.325\textwidth]{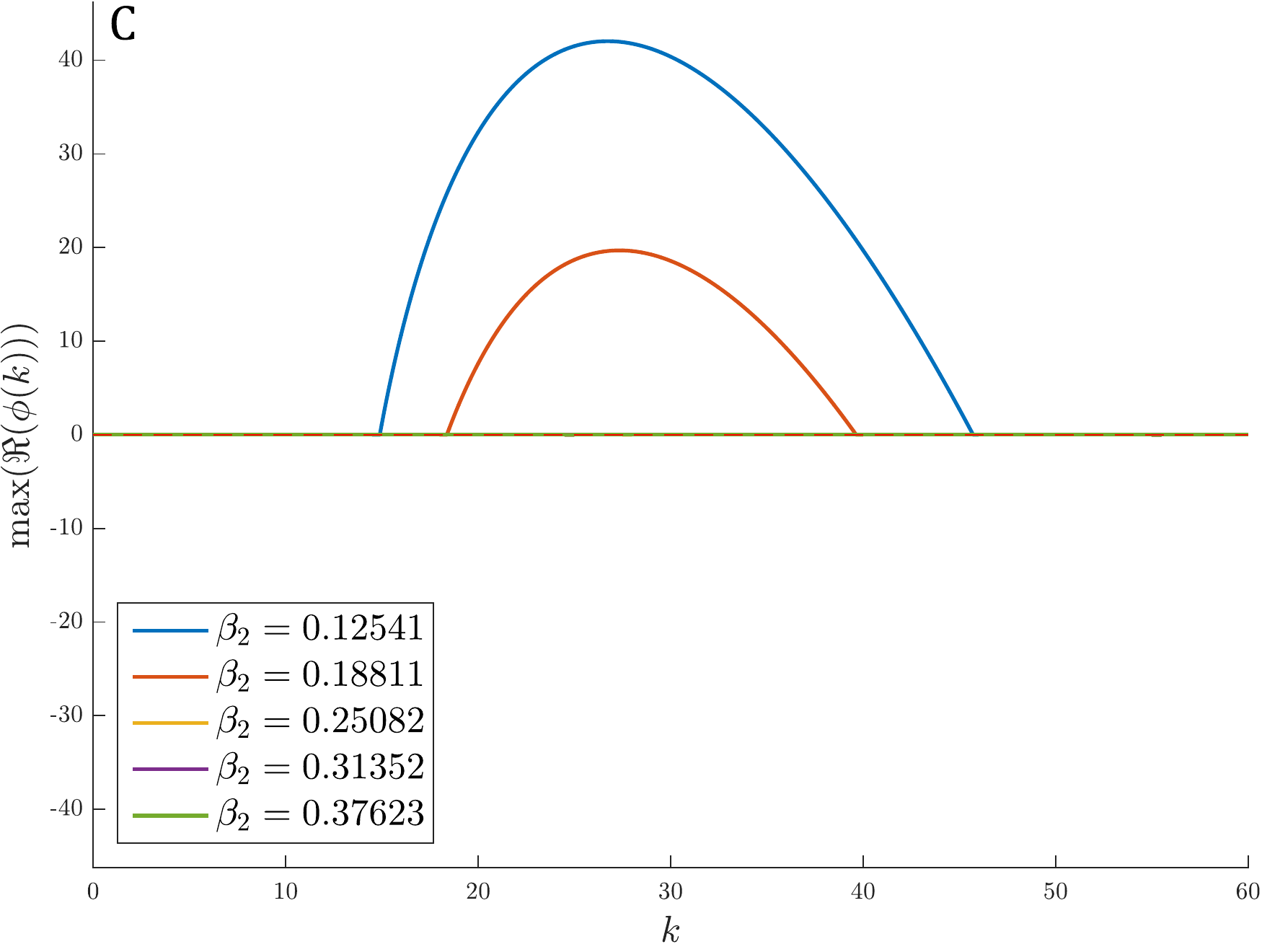}
\end{center}
\caption{Patterning space, parameter condition and dispersion relations for the uncoupled poromechano-chemical model. (A) Predicted pattering space for a selected interval in ($\beta_2$, $\beta_3$) parameter
space: boundary constructed from \eqref{eq:condUC1} (red plain); from \eqref{eq:condUC2} (green-dashed); and from \eqref{eq:cond-hom} (blue-dot-dashed). (B) Parameter coefficient condition $a_0$. Curves are drawn from the critical value $\beta_{2,c}$ (yellow) and for 25\% and 50\% increase/decrease of the parameter values. The critical parameter value is located in (A) and denoted by a magenta cross. (C) Associated dispersion relations, where the colour code is kept identical as in (B).}
\label{fig:lin_ex3}
\end{figure}

By definition, $a_0$ is a polynomial of even order with respect to $k^2$, guaranteeing that there exists at least one local extremum. Therefore, we look for the critical wave number $k_c^2>0$,  obtained by solving the equation $a_0^\prime(k_c^2)=0$, that substituting in $a_0$ will lead to the equation $a_0(k_c^2,\varphi_c)=0$, with $\varphi_c$ \cblue{being the critical parameter  to}  analyse. In the studied scenario, $a_0^{\prime}$ is a cubic polynomial with respect to $k^2$. \cblue{Therefore the following criteria may hold for $a_0$ to be negative} 
\begin{align}
D_2(\beta_3-\beta_2) - D_1(\beta_2+\beta_3)^3 &> 0, \label{eq:condUC1} \\
36\left( \AOCIII \right)^2\left( \AOCII \right)^2 - 128 \AOCIV \left( \AOCII \right)^3 &> 0. \label{eq:condUC2}
\end{align}
\cblue{Inequality \eqref{eq:condUC1} arises from  the Routh-Hurwitz conditions, and it enforces that the coefficient has a real positive part}. Combined with condition \eqref{eq:cond-hom}, it gives an interval for the ratio $(\beta_2+\beta_3)^3 / (\beta_3 - \beta_2)$  \cblue{where} Turing instabilities are reached. As the discriminant of $a_0$ is null, we look for the discriminant \eqref{eq:condUC2} of the derivative $a_0^{\prime}$ to force $k^2\in\mathbb{R}$. Figure~\ref{fig:lin_ex3} (panels (A)-(C)) presents the patterning space based on the implicit functions defined in \eqref{eq:condUC1} and \eqref{eq:condUC2} for the $(\beta_2,\beta_3)$ space. The plots  suggest that increasing the value of the production basal rate ($\beta_3$) of the inhibitor $w_2$ leads to \cblue{a larger interval of possible basal rates} ($\beta_2$) of the activator $w_1$. The condition \eqref{eq:cond-hom}, represented by the blue-dot-dashed curve in Figure~\ref{fig:lin_ex3}(A), is absent for the $\rho=0$ scenario. \cblue{This enlarges the patterning space, and therefore the presence of} acceleration in the momentum equilibrium equation in turn leads to a restriction of the Turing space.

\begin{figure}[!t]
\begin{center}
\includegraphics[width=0.45\textwidth]{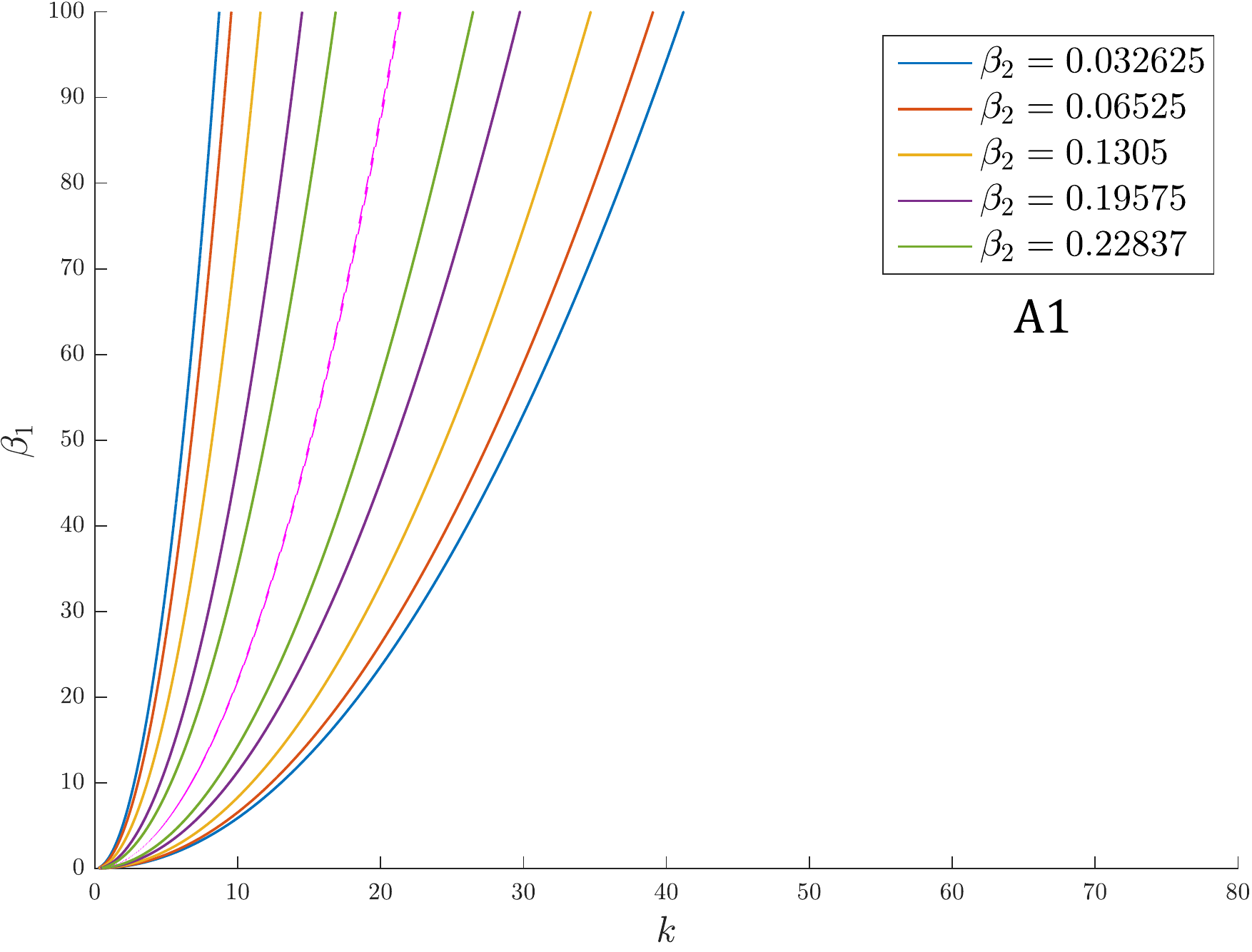}
\includegraphics[width=0.45\textwidth]{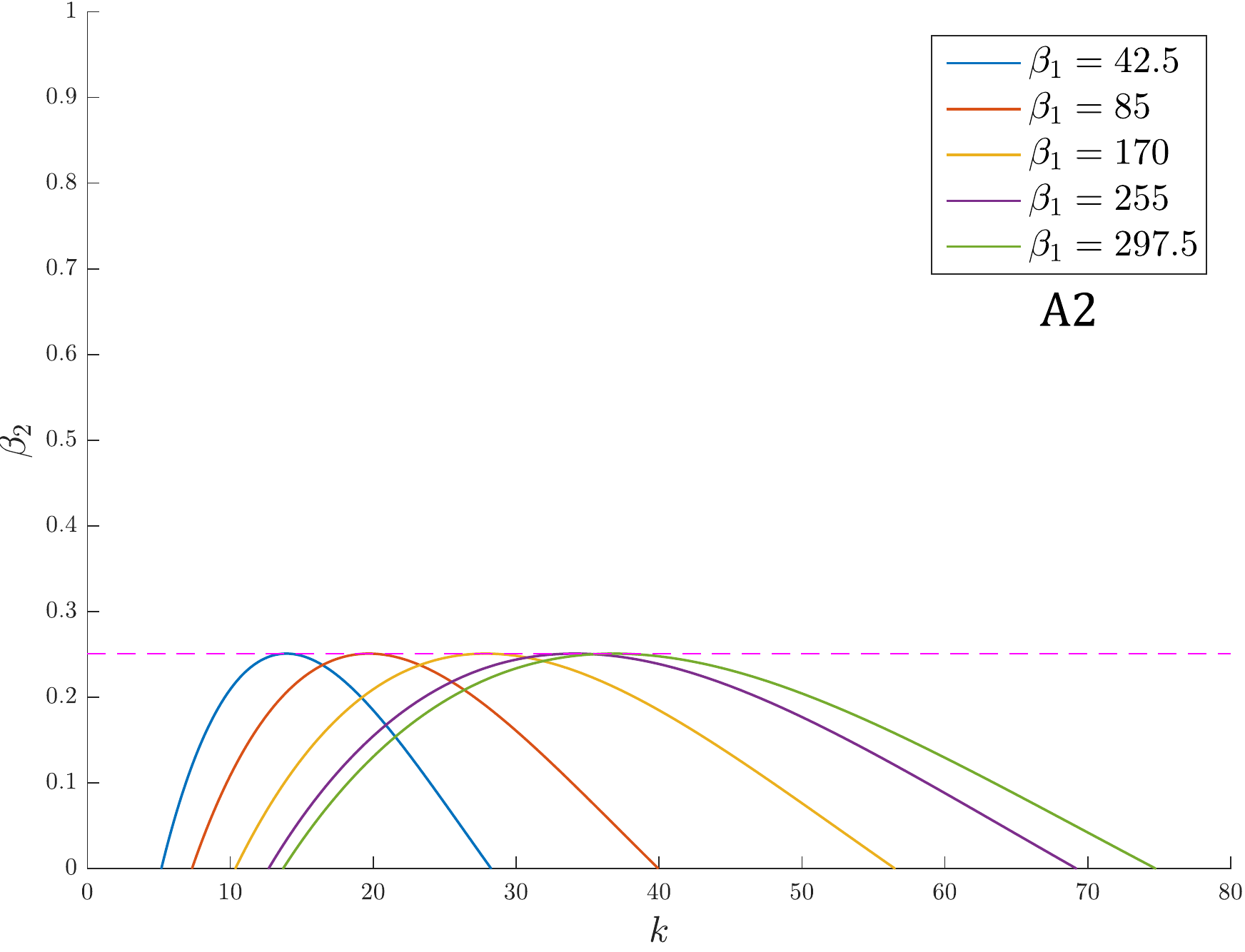}
\end{center}
\caption{\cblue{Contour plots of the Routh-Hurwitz conditions for the characteristic polynomial of the uncoupled system \eqref{eq1:mass}-\eqref{eq1:momentum}. (A1) Null level-set of $P_2$-polynomial coefficients $a_0$ defined in \eqref{eq:la-poly} for different $\beta_2$'s. (A2) Null level set of $P_2$-polynomial coefficients $a_0$ defined in \eqref{eq:la-poly} for different $\beta_1$'s. The magenta curve is used as reference for the critical $\beta_2$ value computed from $a_0$, setting $\beta_1=170$.}}
\label{fig:lin_ex4}
\end{figure}

\cblue{
Taking into account the stability of the homogeneous spatial case, the constraint \eqref{eq:condUC1} is key for generating instabilities. If \eqref{eq:condUC1} is unsatisfied, all the dispersion relation coefficients are positive, whatever the choice of the parameter values. Figure \ref{fig:lin_ex3}(A) shows that  condition \eqref{eq:condUC1} is not sufficient to onset instability, and the patterning space is smaller than the region delimited by the red curve in Figure \ref{fig:lin_ex3}(A). The question is then to check whether, after selecting a $\beta_2$ in the region delimited by the green-dashed and the red-solid curves in the patterning space, we can  have instability by varying the value of $\beta_1$ only. Figures \ref{fig:lin_ex4}(A1)-(A2) display the behaviour of the uncoupled system for different values of $\beta_1$ and $\beta_2$, fixing all other coefficients. When the system reaches the critical value for $\beta_2$ (magenta-dashed curve in Fig. \ref{fig:lin_ex4}(A1)), we see that $a_0$ tends to be strictly positive and prevents any instability from that condition (the other Routh-Hurwitz conditions are satisfied for the selected fixed parameters). Nevertheless, $\beta_1$ still has an influence by extending the space scale where instability can occur (see Fig. \ref{fig:lin_ex4}(A2)), without affecting the value of the critical value of $\beta_2$ (all the null level-sets do not exceed the reference magenta-dashed curve in Fig. \ref{fig:lin_ex4}(A2)). This confirms that $\beta_2$ (or eventually $\beta_3$) are appropriate parameters to analyse the stability of the system in the uncoupled scenario.
}

\subsection{Coupled system - Null production/degradation rates}
%
%

\begin{figure}[!t]
\begin{center}
\includegraphics[width=0.375\textwidth]{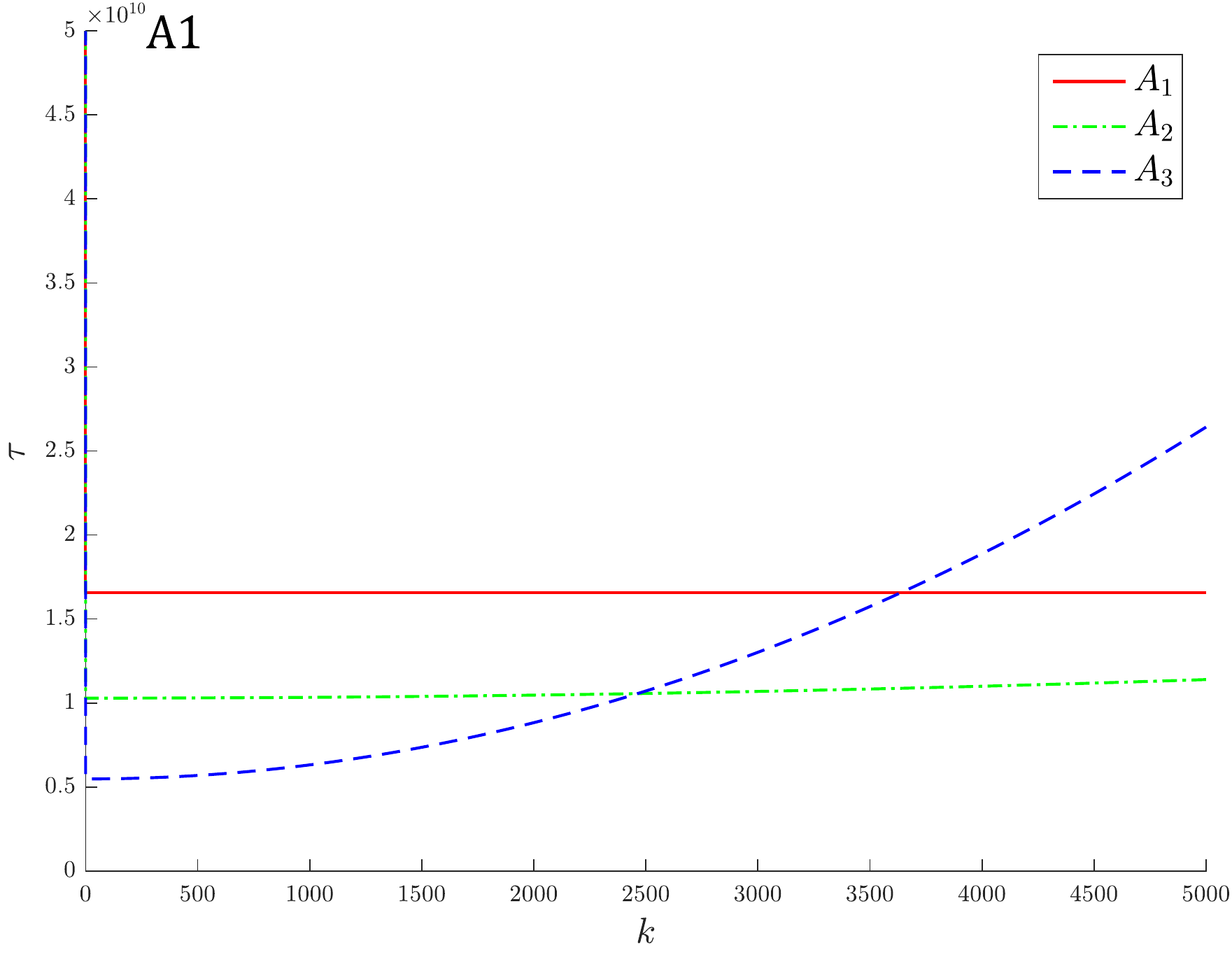}
\includegraphics[width=0.375\textwidth]{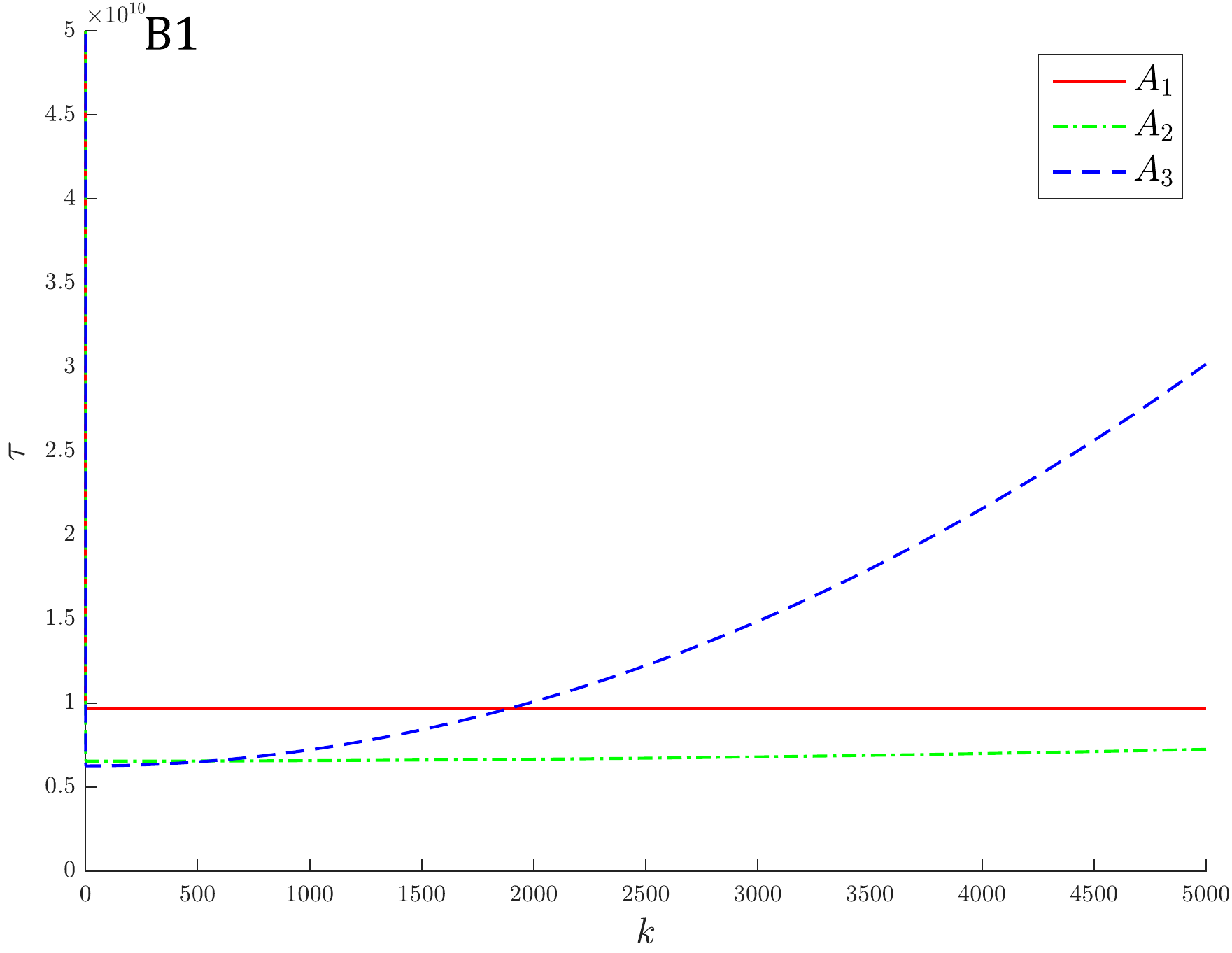}\\
\includegraphics[width=0.375\textwidth]{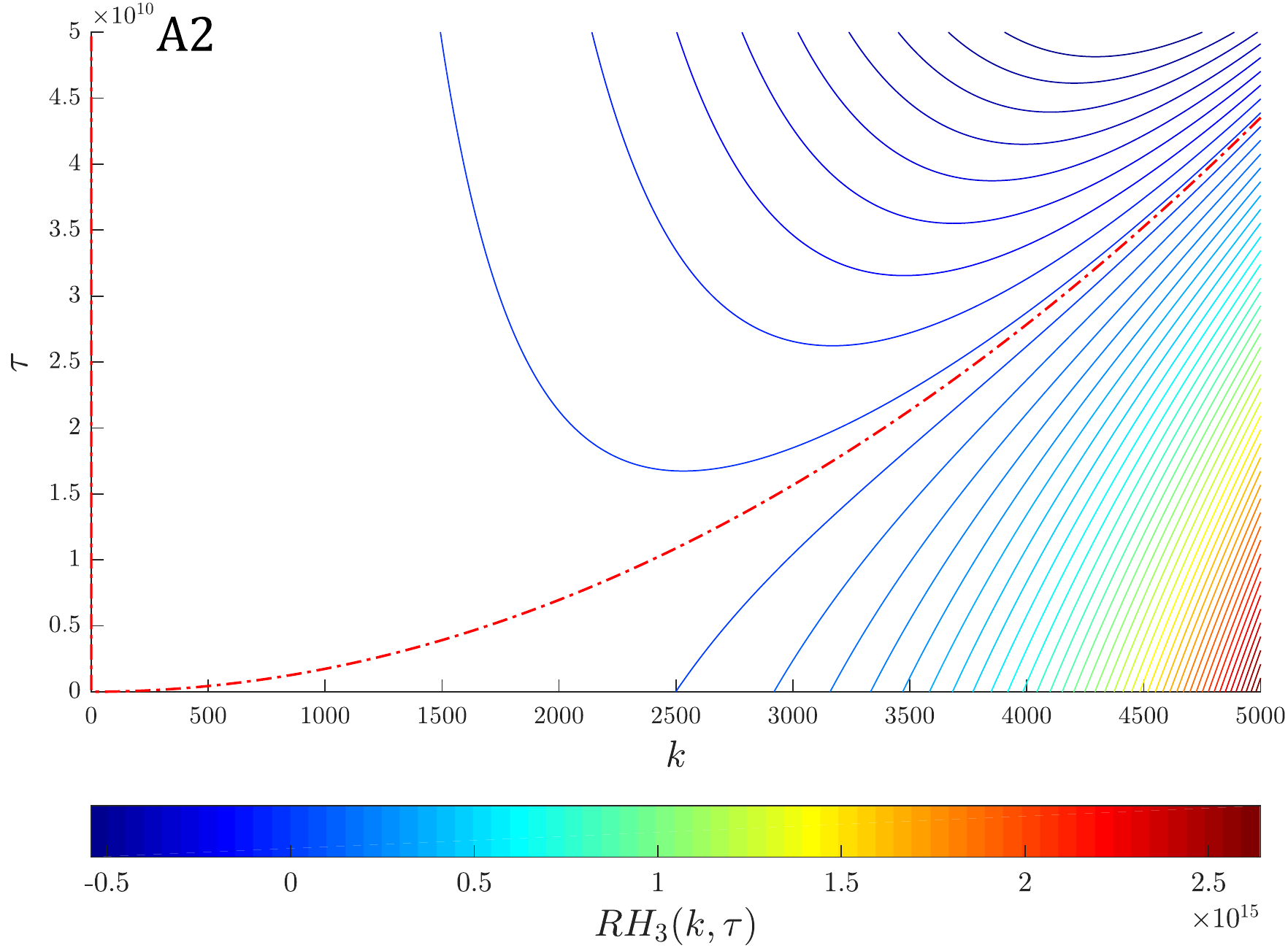}
\includegraphics[width=0.375\textwidth]{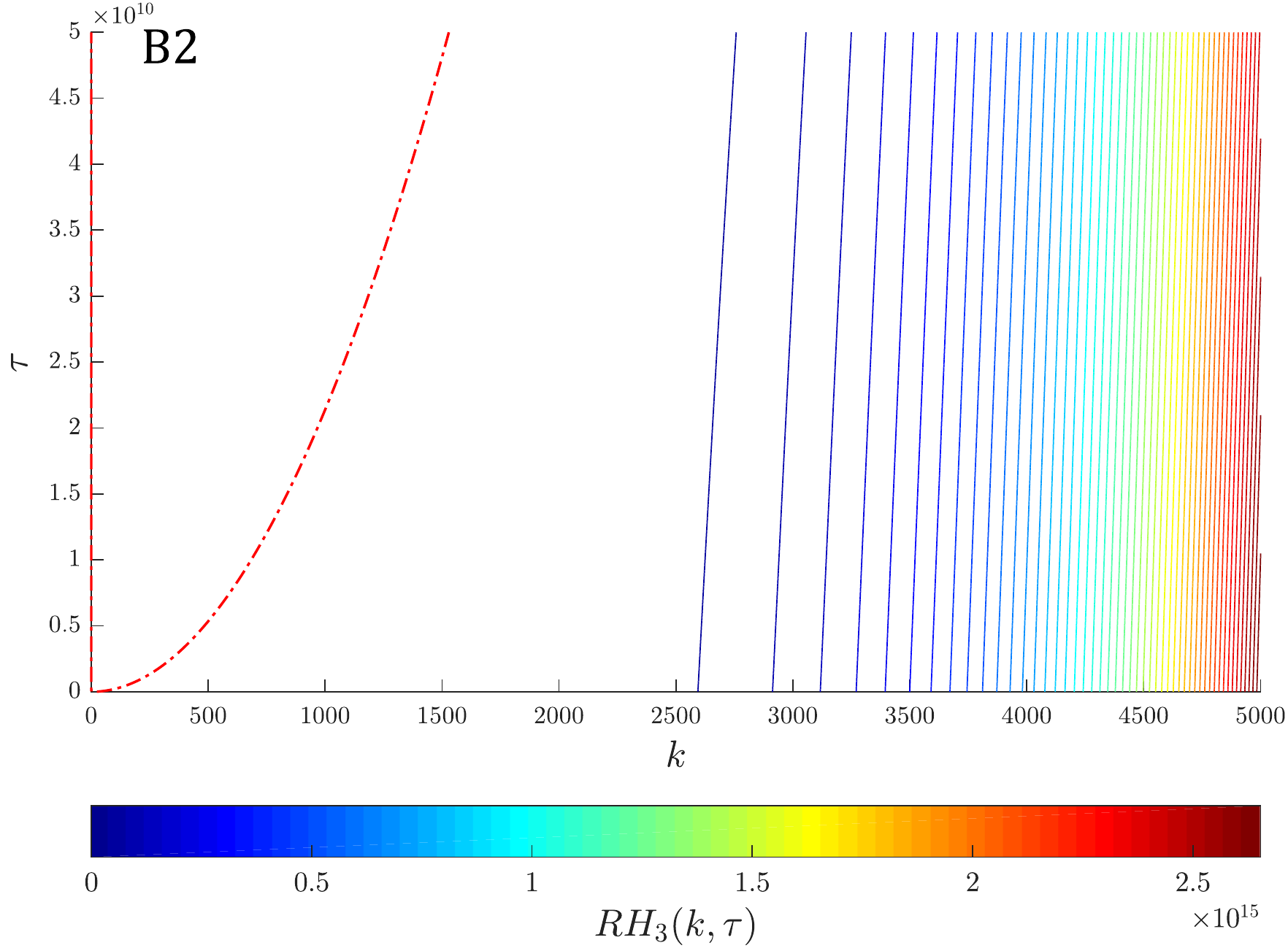}\\
\includegraphics[width=0.375\textwidth]{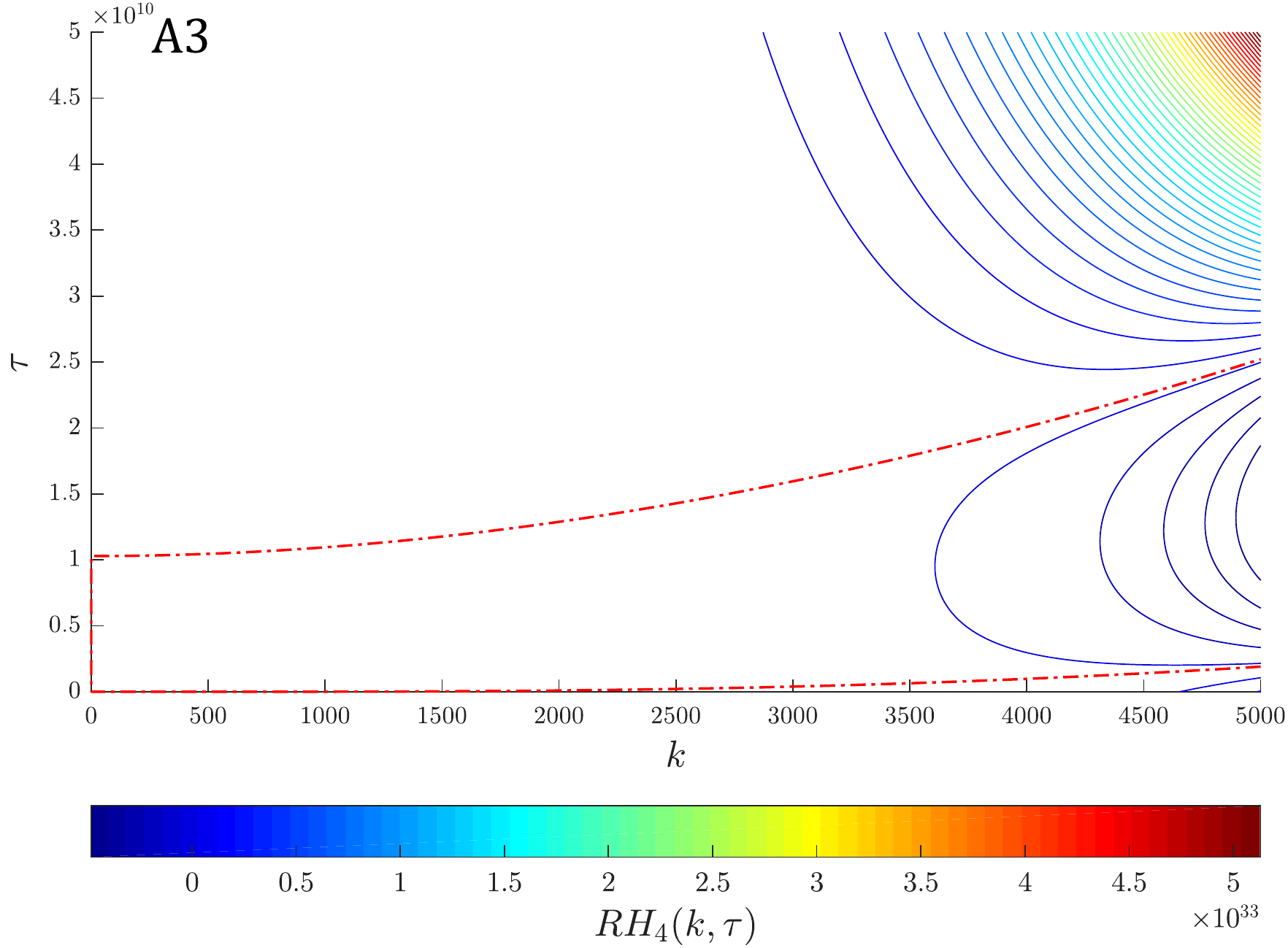}
\includegraphics[width=0.375\textwidth]{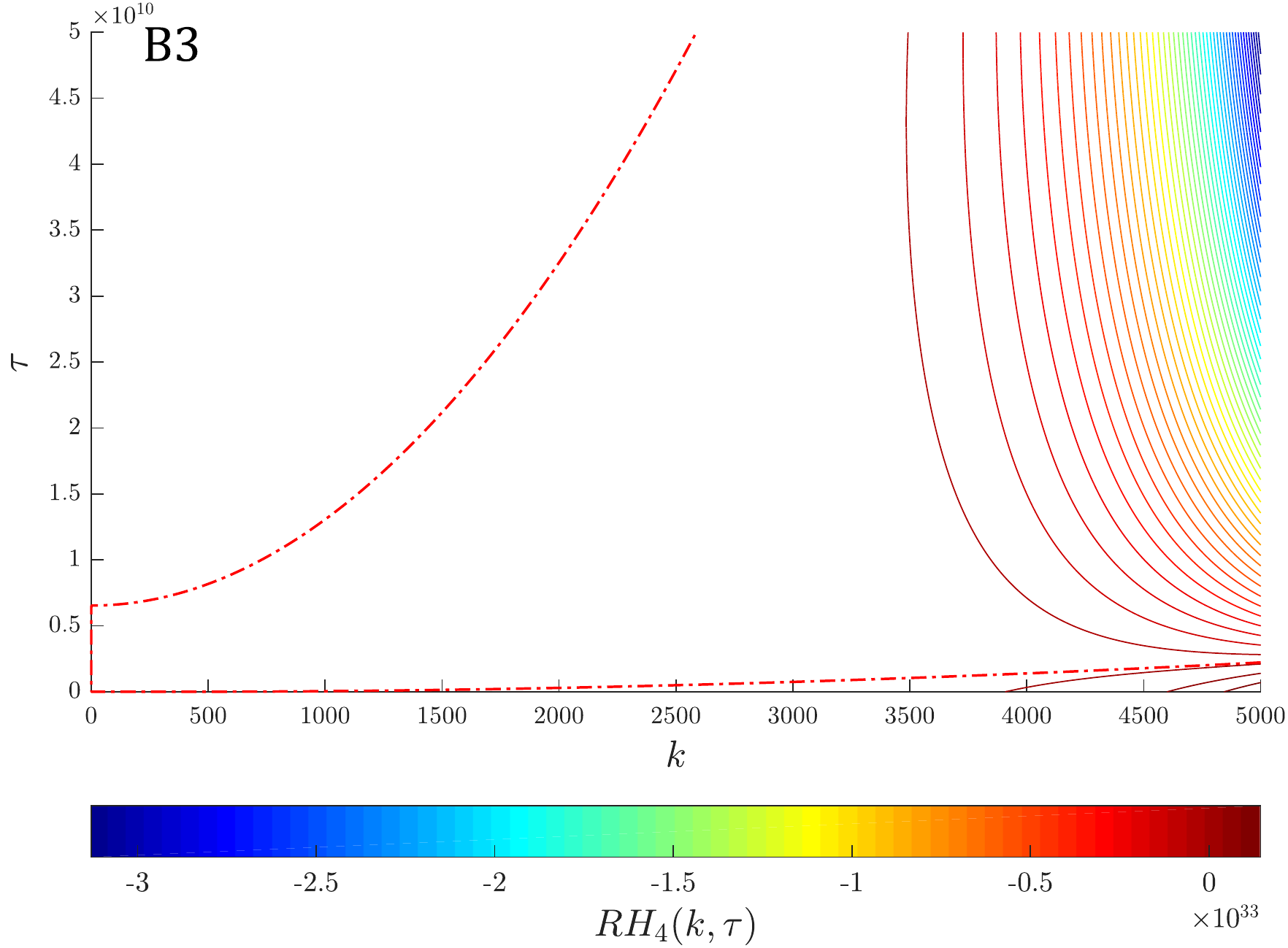}\\
\includegraphics[width=0.375\textwidth]{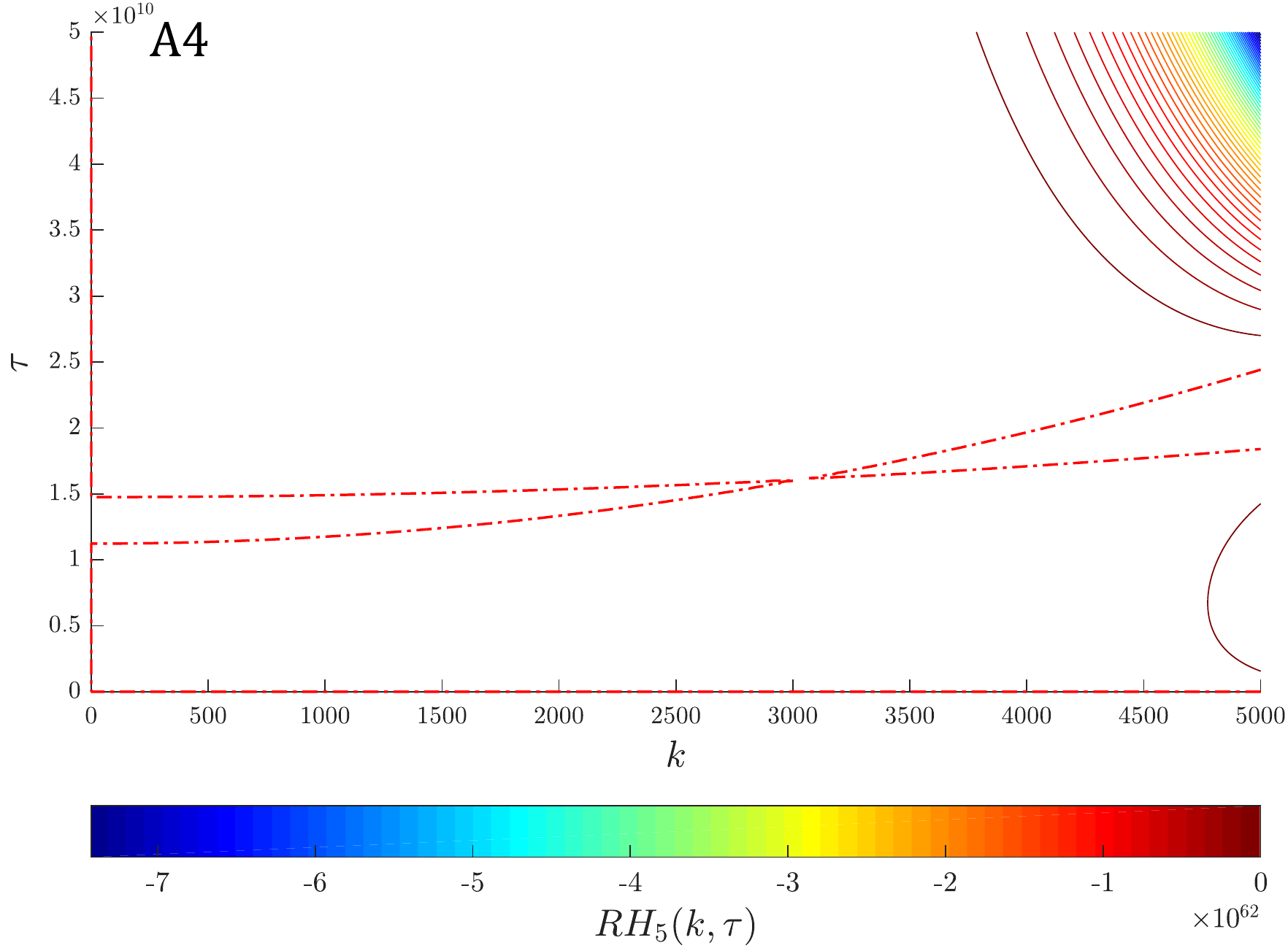}
\includegraphics[width=0.375\textwidth]{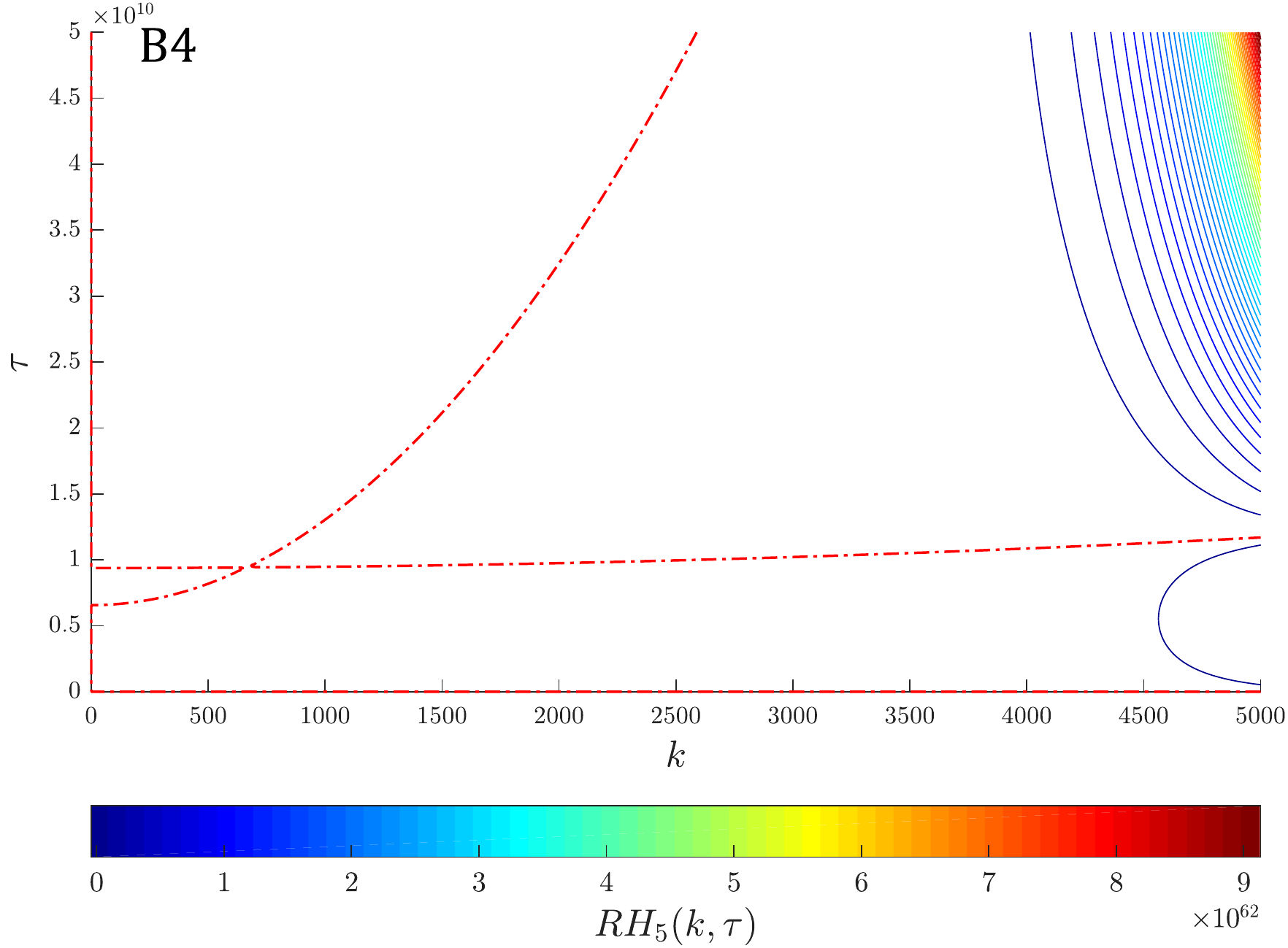}
\end{center}
\caption{Contour plots of the Routh-Hurwitz conditions for the characteristic polynomial of system \eqref{eq1:mass}-\eqref{eq1:momentum} for $\beta_1=0$ and $\rho=1$. (A1) Null level set of $P_2$-polynomial coefficients $a_i^{(1)}$ for the coupling term $\btheta^{(1)}$. (A2) Level sets (100) of condition \eqref{eq:rh-cond2} with null level of the condition \eqref{eq:rh-cond2} (red dot-dashed) for the coupling term $\btheta^{(1)}$. (A3-A4) Similar analysis for the conditions \eqref{eq:rh-cond3} and \eqref{eq:rh-cond4} respectively. (B1-B4) Similar analysis for the coupling term $\btheta^{(2)}$.}
\label{fig:lin_ex5}
\end{figure}

In contrast with classical reaction-diffusion systems, due to the coupling with the poroelastic deformations we can perfectly encounter cases where production/degradation rates are missing. For the first case of $\beta_1=0$ (corresponding to a pure convection-diffusion chemical system) the characteristic polynomial for $\rho \neq 0$ is still of order 5 with the new coefficients
{\small
\begin{align*}
&\hspace{2.0cm} A_5^{(1)}(k^2) = A_5(k^2), \qquad A_4^{(1)}(k^2) = \AIVCI k^2, \\
A_3^{(1)}(k^2) &= \AIIICII k^4 + \left[ \AIIICIbO \right] k^2, \\
A_2^{(1)}(k^2) &= \AIICIII k^6 + \Bigg[ \AIICIIbO \\
&\qquad + \gamma \left( \frac{\kappa}{\eta}\left( (\beta_2 + \beta_3 )\theta_1^{(i)} +\frac{\beta_3}{( \beta_2 + \beta_3 )^2}\theta_2^{(i)} \right) + c_0\left( (\beta_2 + \beta_3 )\theta_1^{(i)}D_2 +\frac{\beta_3}{( \beta_2 + \beta_3 )^2}\theta_2^{(i)}D_1 \right) \right) \Bigg] k^4, \\
A_1^{(1)}(k^2) &= \left[ \AICIIIbis \right] k^6, \\
A_0^{(1)}(k^2) &= \frac{\kappa}{\eta}(2\mu+\lambda) \AOCIV k^8,
\end{align*}
}
where each $\btheta^{(i)}$ (for $i=1,2$) is specified as 
\begin{equation}
\btheta^{(1)} = \begin{pmatrix}
\theta^{(1)}_1 \\
\theta^{(1)}_2
\end{pmatrix} = -\tau(1,1)^{\intercal}, \qquad \btheta^{(2)} = - 2\tau(1,0)^{\intercal}.
\label{eq:theta}
\end{equation}
\cblue{The spatially homogeneous case of $k=0$ remains stable irrespective of the parameter values, by definition. Also note that} for a general $k$, all the coefficients that do not include the coupling parameter $\gamma$ are strictly positive, and consequently they do not influence the stability of the coupled system. \cblue{As there is no restriction on the parameter values} for the homogeneous case, any parameter can be chosen as the critical one. In order to analyse the effect of poromechanics on the chemical system, we regard $\tau$ as the parameter of interest.

Figure~\ref{fig:lin_ex5} presents contour plots of the Routh-Hurwitz conditions with respect to the wave number $k$ and the parameter $\tau$, for both $\btheta$ defined in \eqref{eq:theta}. From Figures~\ref{fig:lin_ex5}(A1) and \ref{fig:lin_ex5}(B1), we observe that  $a_3$ is the first coefficient to break the inequality condition with respect to $\tau$, and this occurs at a low wave number. Along $k$, the parabolic shape of the null levels shows how, depending on the size of the system, any of the three coefficients can break the Routh-Hurwitz inequality. Consequently, the coupled system presents complex instabilities and makes it difficult to choose only one coefficient to analyse the full patterning space. Nevertheless, patterns are reachable for \cblue{large values of $\tau$, irrespective of the} wave number (at least in the presented interval). This is contrary to the uncoupled case, where all the $a_i$'s are strictly positive. Figures~\ref{fig:lin_ex5}(A2)-(A4) and \ref{fig:lin_ex5}(B2)-(B4) depict the sign of the conditions \eqref{eq:rh-cond2}-\eqref{eq:rh-cond4} for $\btheta^{(1)}$ and $\btheta^{(2)}$, respectively. We readily see that the instability region starts already at a value of $\tau$ ($\sim10^{5}$) lower than that provided by the $a_i$'s. This \cblue{emphasises further the effects of the poromechanical} coupling into pattern formation.

\begin{figure}[!t]
\begin{center}
\includegraphics[width=0.45\textwidth]{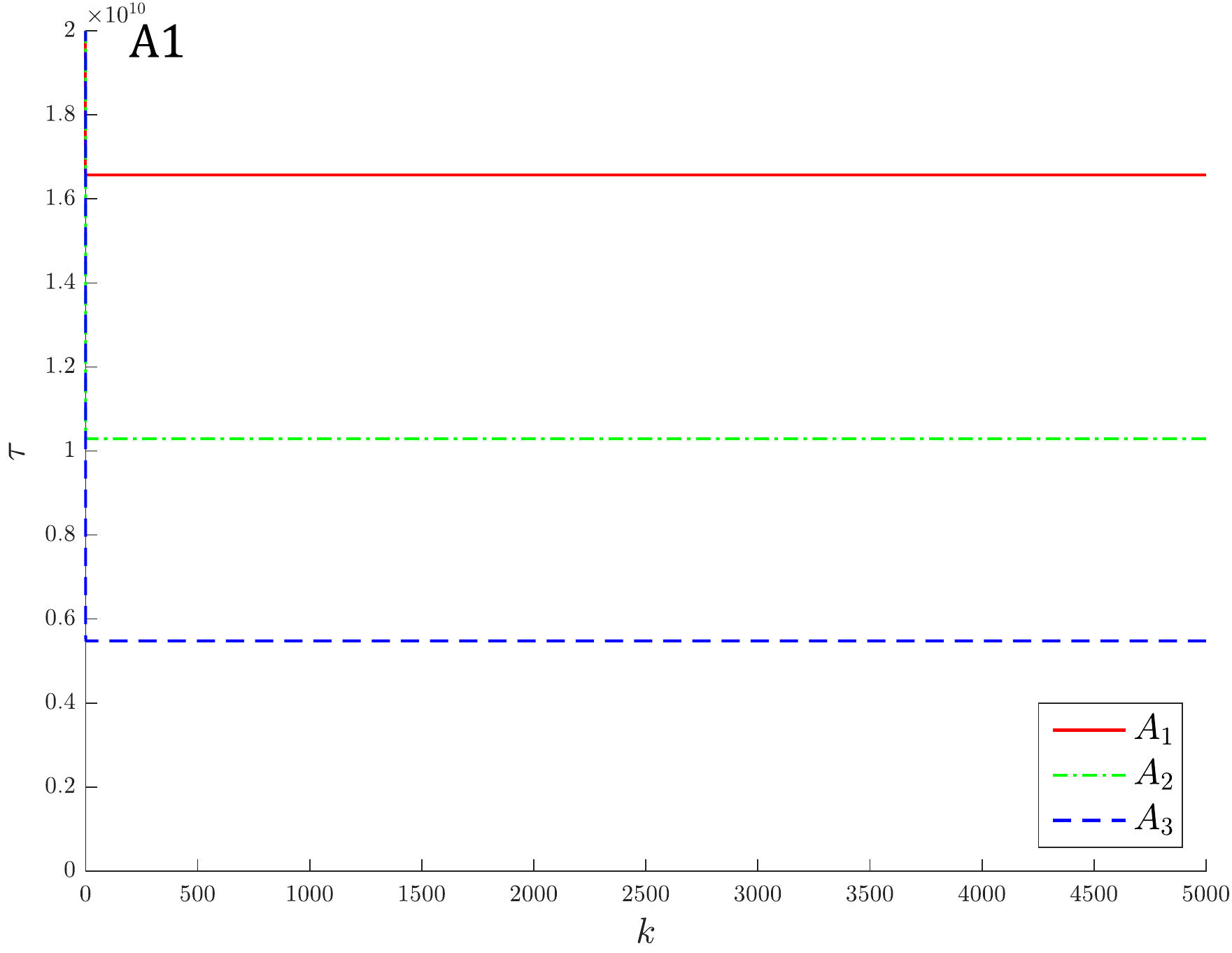}
\includegraphics[width=0.45\textwidth]{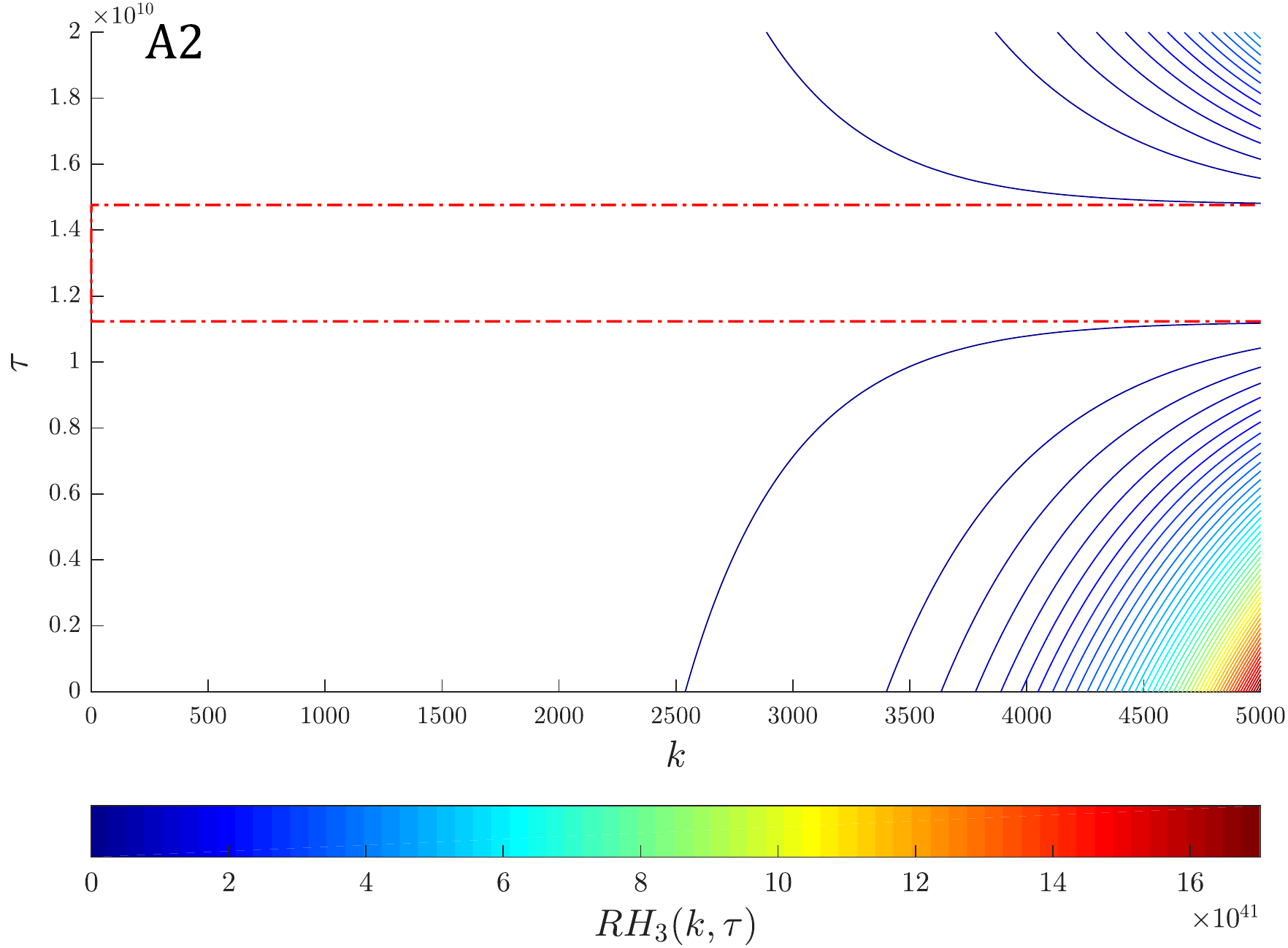}\\
\includegraphics[width=0.45\textwidth]{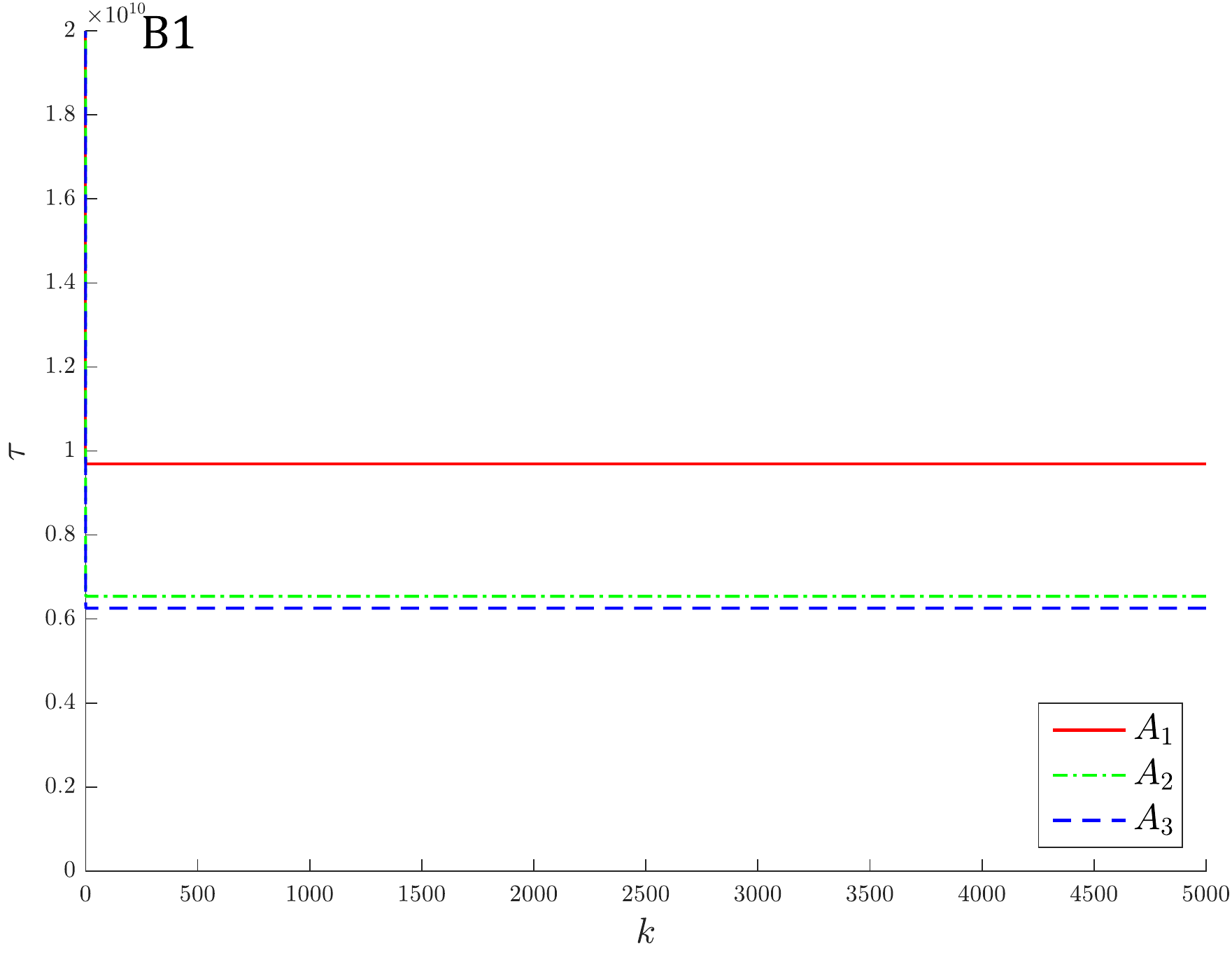}
\includegraphics[width=0.45\textwidth]{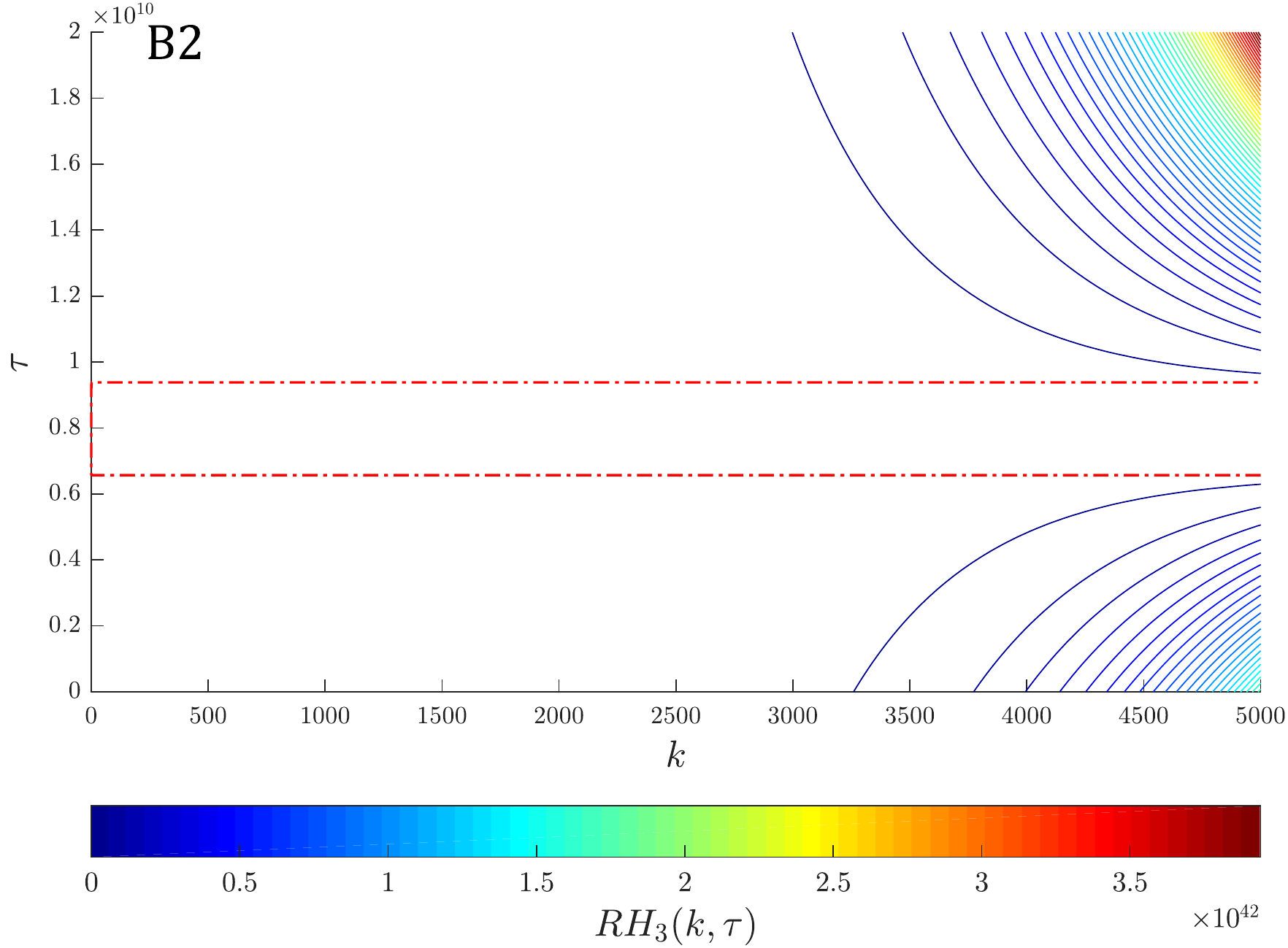}
\end{center}
\caption{Contour plots of the Routh-Hurwitz conditions for the characteristic polynomial of system \eqref{eq1:mass}-\eqref{eq1:momentum} for $\beta_1=0$ and $\rho=0$. (A1) Null level set of $P_2$-polynomial coefficients $a_i^{(1)}$ for the coupling term $\btheta^{(1)}$. (A2) Level sets (100) of condition \eqref{eq:rh2-cond} with associated null level (red dot-dashed) for the coupling term $\btheta^{(1)}$. (B1-B2) Similar analysis for the coupling term $\btheta^{(2)}$.}
\label{fig:lin_ex6}
\end{figure}

The choice of the coupling function $\btheta$ has a non-intuitive influence on the Routh-Hurwitz conditions. While the area of the instability region for conditions \eqref{eq:rh-cond2} and \eqref{eq:rh-cond4} (\textit{cf.} respectively Figs.~\ref{fig:lin_ex5}(A2), (B2) and Figs.~\ref{fig:lin_ex5}(A4), (B4)) \cblue{is decreased, it is in turn increased} for condition \eqref{eq:rh-cond3} (\textit{cf.}  Figs.~\ref{fig:lin_ex5}(A3), (B3)). Proposing a patterning space based uniquely on the 
parameters can still be quite difficult due to the complexity of the Routh-Hurwitz conditions and by the evolution of  the inequality constrain through the wave number. Nonetheless, we can still show how the poromechanics coupled to a pure convection-diffusion system may generate patterns. In the absence of acceleration, the analysis of the Routh-Hurwitz conditions is quite similar (see Figure~\ref{fig:lin_ex6}) to the general case. We observe again that for a specific $\tau$, the conditions can break for any coupling function $\btheta$. Contrary to the case of $\rho=1$, here the defect on the criteria occurs almost at one specific value of $\tau$ for any wave number $k$ (see Figs.~\ref{fig:lin_ex6}(A1), (B1)). The coupling function seems to move closer to each other the levels of the different coefficients $a_i$, and it also decreases the interval of $\tau$ that leads to breaking the second condition \eqref{eq:rh2-cond}.

Secondly, we impose either the basal rate of the activator, $\beta_2$, or the inhibitor, $\beta_3$, to be zero. Again the characteristic polynomial $P_2(\phi,k^2)$ for $\rho\neq0$ is still of order 5 with the new coefficients given by 
{\small
\begin{align*}
 A_5^{(2)}(k^2) &= A_5^{(3)}(k^2) = A_5(k^2), \quad
A_4^{(2,3)}(k^2) = \AIVCI  k^2 - \Bigg\{\begin{matrix} \rho c_0 \beta_1( 1 - \beta_3^2 ) \\ -\rho c_0 \beta_1( 1 + \beta_2^2 )
\end{matrix}, \\
A_3^{(2,3)}(k^2) &= \AIIICII  k^4 \\
&\quad + \Bigg\{\begin{matrix}
\left[ \AIIICIbD + \Psi(\beta_3;\tau) \right] k^2 + c_0\rho\beta_1^2\beta_3^2 \\ 
\left[ \AIIICIbT + \Psi(\beta_2;\tau) \right] k^2 + c_0\rho\beta_1^2\beta_2^2
\end{matrix}, \\
A_2^{(2,3)} &= \AIICIII k^6 \\ 
&\quad + \Bigg\{\begin{matrix}
\left[ \AIICIIbD + \Psi(\beta_3;\tau) \right] k^4 \\ 
\left[ \AIICIIbT + \Psi(\beta_2;\tau) \right] k^4
\end{matrix} \\
&\quad + \Bigg\{\begin{matrix}
\left[ \AIICIbD + \Psi(\beta_3;\tau) \right] k^2 \\
\left[ \AIICIbT + \Psi(\beta_2;\tau) \right] k^2
\end{matrix}, \\
A_1^{(2,3)} &= \Bigg\{\begin{matrix}\left[ \AICIII + \Psi(\beta_3;\tau) \right] k^6 \\
\left[ \AICIII + \Psi(\beta_2;\tau) \right] k^6
\end{matrix}
\\
&\quad - \Bigg\{\begin{matrix}
\left[ \AICIIbD + \Psi(\beta_3;\tau) \right]k^4 \\ 
- \left[ \AICIIbT - \Psi(\beta_2;\tau) \right]k^4
\end{matrix}  \\
&\quad + \Bigg\{\begin{matrix}
( c_0( 2\mu + \lambda ) + \alpha^2 )\beta_1^2\beta_3^2 k^2 \\
( c_0( 2\mu + \lambda ) + \alpha^2 )\beta_1^2\beta_2^2 k^2 \\
\end{matrix}, \\
A_0^{(2,3)} &= \Bigg\{\begin{matrix}
\frac{\kappa}{\eta}(2\mu + \lambda)k^4\left[ \AOCIV k^4 - \beta_1(D_2 - D_1\beta_3^2)k^2 + \beta_1^2\beta_3^2 \right] \\
\frac{\kappa}{\eta}(2\mu + \lambda)k^4\left[ \AOCIV k^4 + \beta_1(D_2 + D_1\beta_2^2)k^2 + \beta_1^2\beta_2^2 \right]
\end{matrix},
\end{align*}
}
where $\Psi$ is a generic function that summarises the appropriate coupling term defined precisely in \eqref{eq:la-poly} and the upper (resp. lower) line in braces exhibits the coefficient description for $\beta_2=0$ (resp. $\beta_3=0$). We observe from the new set of coefficients that the stability of the system behaves differently with respect to \cblue{a zero basal rate}. In an uncoupled scenario, the case $\beta_3=0$ implies that the $a_i$'s terms are strictly positive whatever the value of the coefficients. This indicates that removing the basal production of the inhibitor prevents any patterning (in the uncoupled scenario), and this occurs for any $\rho$ and any wave number $k$. For $\beta_2=0$, the system can enter an instability region if and only if 
\begin{equation*}
0 < \beta_3 < \sqrt{\min\left(1,\frac{D_2}{D_1}\right)}.
\end{equation*}

\begin{figure}[!t]
\begin{center}
\includegraphics[width=0.45\textwidth]{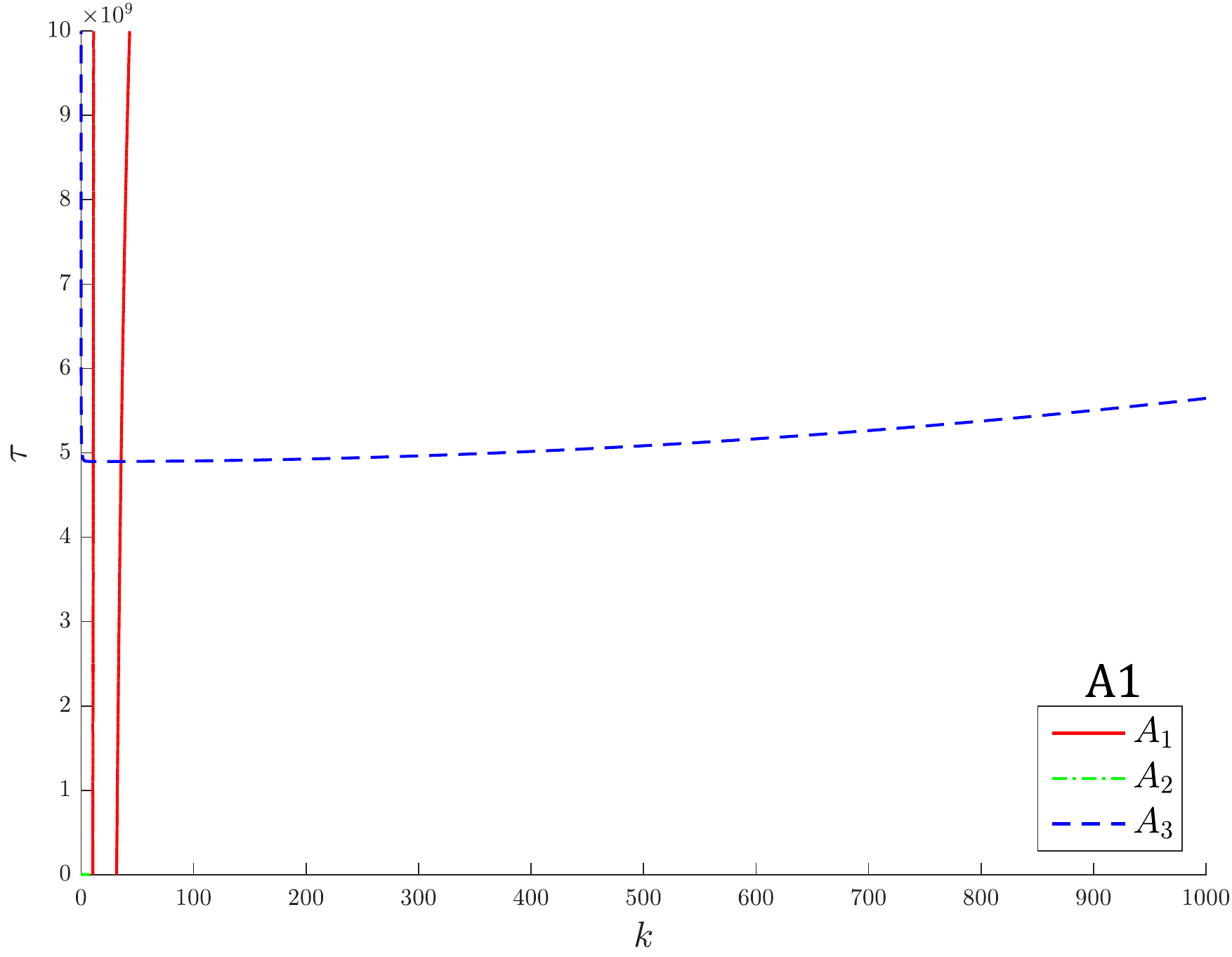}
\includegraphics[width=0.45\textwidth]{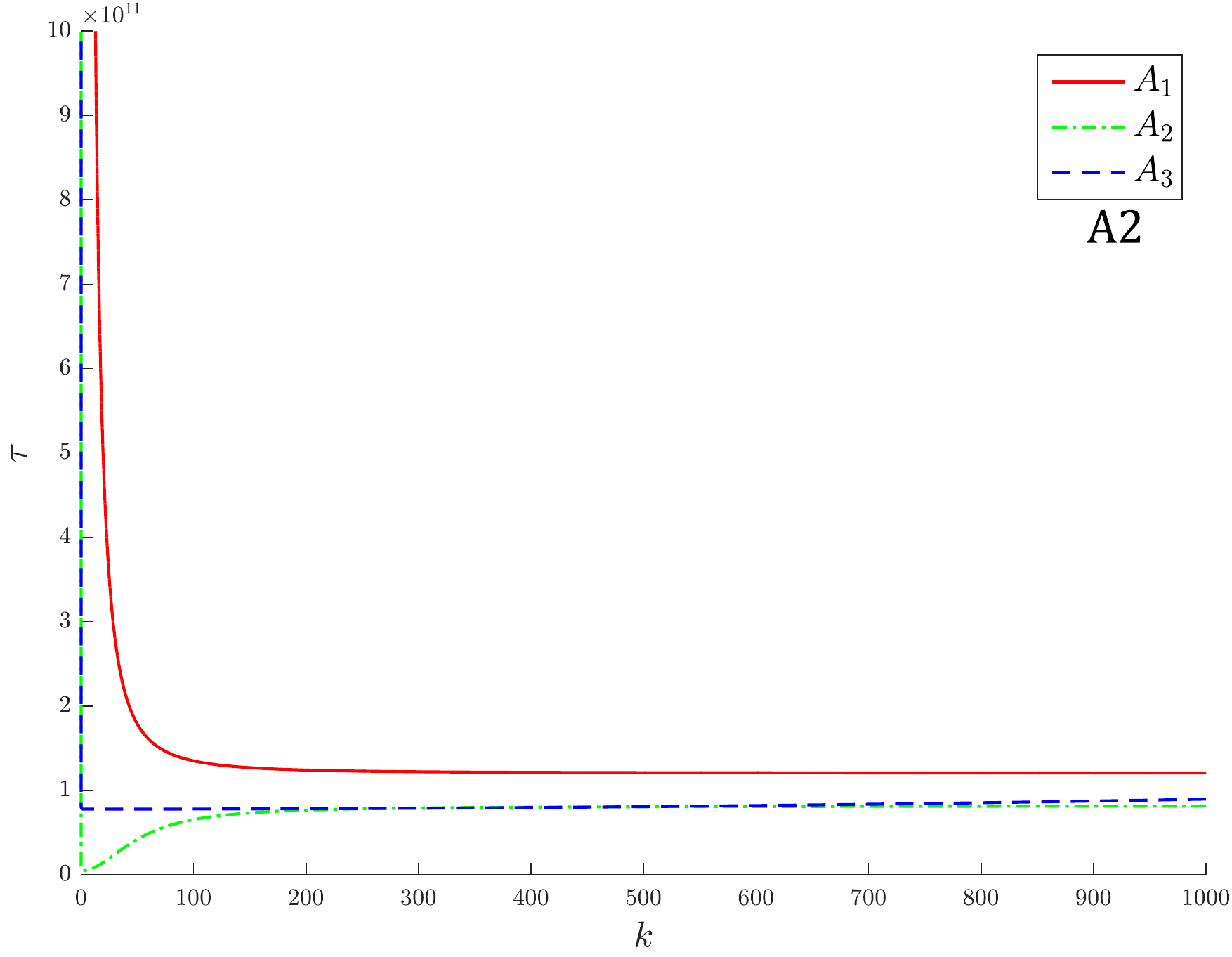}\\
\includegraphics[width=0.45\textwidth]{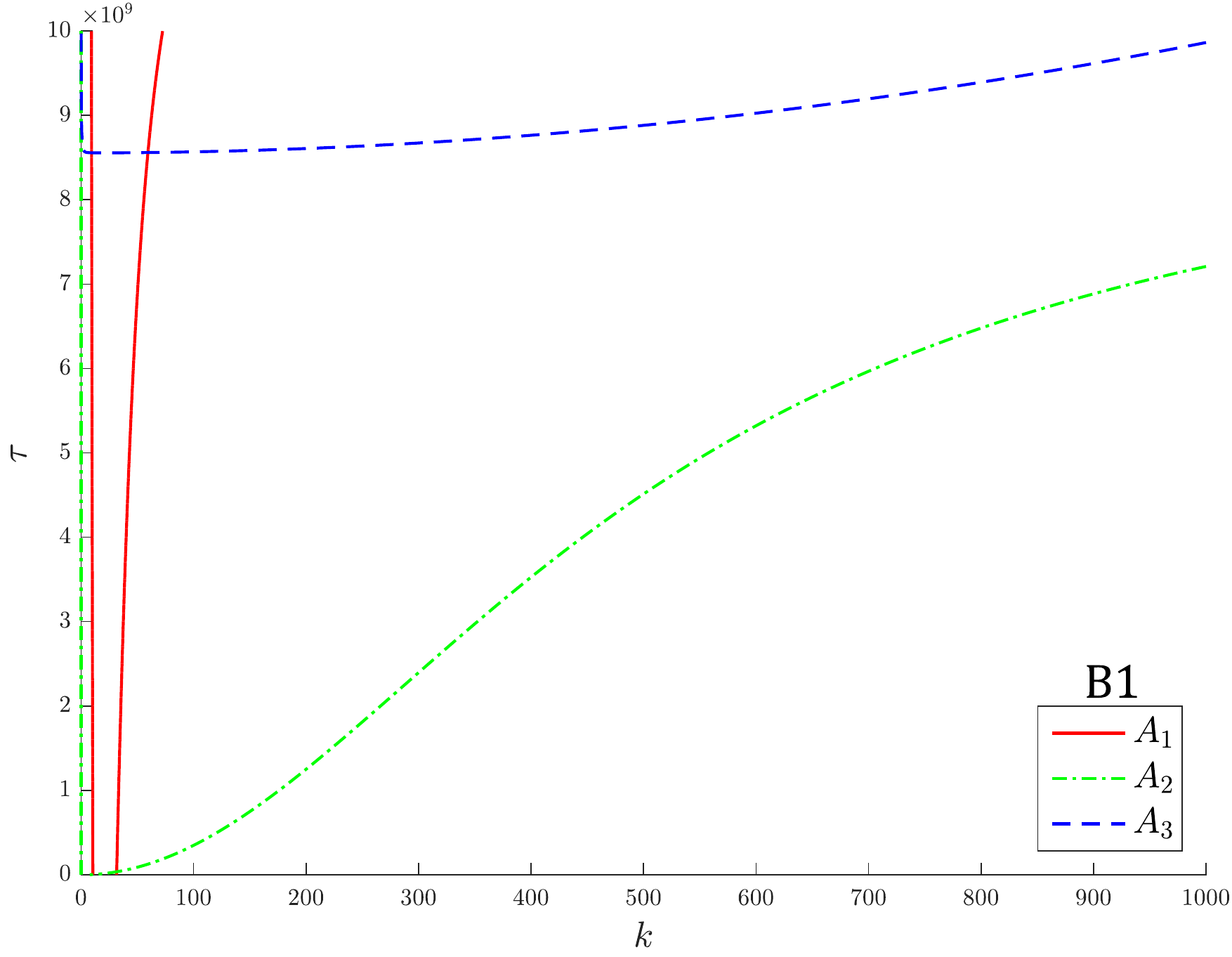}
\includegraphics[width=0.45\textwidth]{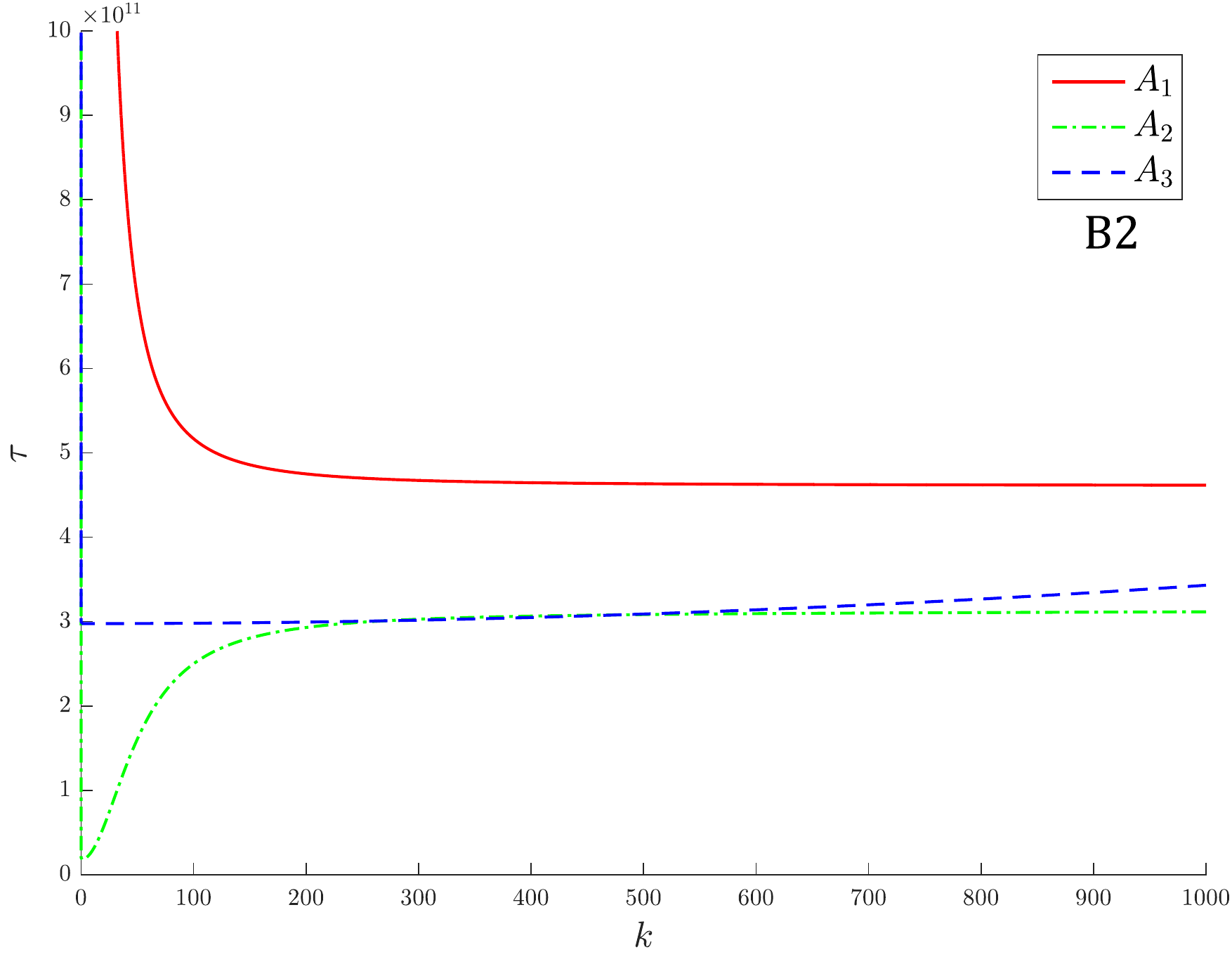}
\end{center}
\caption{Null level set of the characteristic polynomial coefficients $a_i^{(2,3)}$ defined in \eqref{eq:la-poly} for different null basal rate and coupling term $\btheta$. (A1) $\beta_2=0$ with $\btheta^{(1)}$. (A2) $\beta_3=0$ with $\btheta^{(1)}$. (B1-B2) Similar analysis for $\btheta^{(2)}$. In all plots \cblue{we use} $\rho=1$.}
\label{fig:lin_ex7}
\end{figure}

Coupling the convection-reaction-diffusion system to the poromechanics enables both scenarios to reach an instability for some parameter values. We analyse again \cblue{the coupled} system focusing on $\tau$. Figure~\ref{fig:lin_ex7} presents the null level set of the coefficients for both $\btheta$ (\textit{cf.}  \eqref{eq:theta}) and  $\rho=1$. The case $\beta_2=0$ (see Figs.~\ref{fig:lin_ex7}(A1), (B1)) is significantly \cblue{affected} by the choice of $\btheta$. Going from a linear $\btheta^{(1)}$ to a nonlinear $\btheta^{(2)}$ coupling function, the instability region is clearly modified, especially for the coefficient $a_2$, \cblue{leading to a larger area}. Furthermore, starting from \cblue{a given} wave number $k$, the system presents instability whatever the choice of the parameter $\tau$ (\textit{e.g.}, red null level in Fig.~\ref{fig:lin_ex7}(A1)). In $\beta_3=0$, instabilities can be produced using a large value of $\tau$ (see Figs.~\ref{fig:lin_ex7}(A2), (B2)), implying that coupling the diffusion system to the \cblue{poromechanics bypasses the intrinsic stability of an uncoupled system}.

\subsection{General coupled system}

\begin{figure}[!t]
\begin{center}
\includegraphics[width=0.45\textwidth]{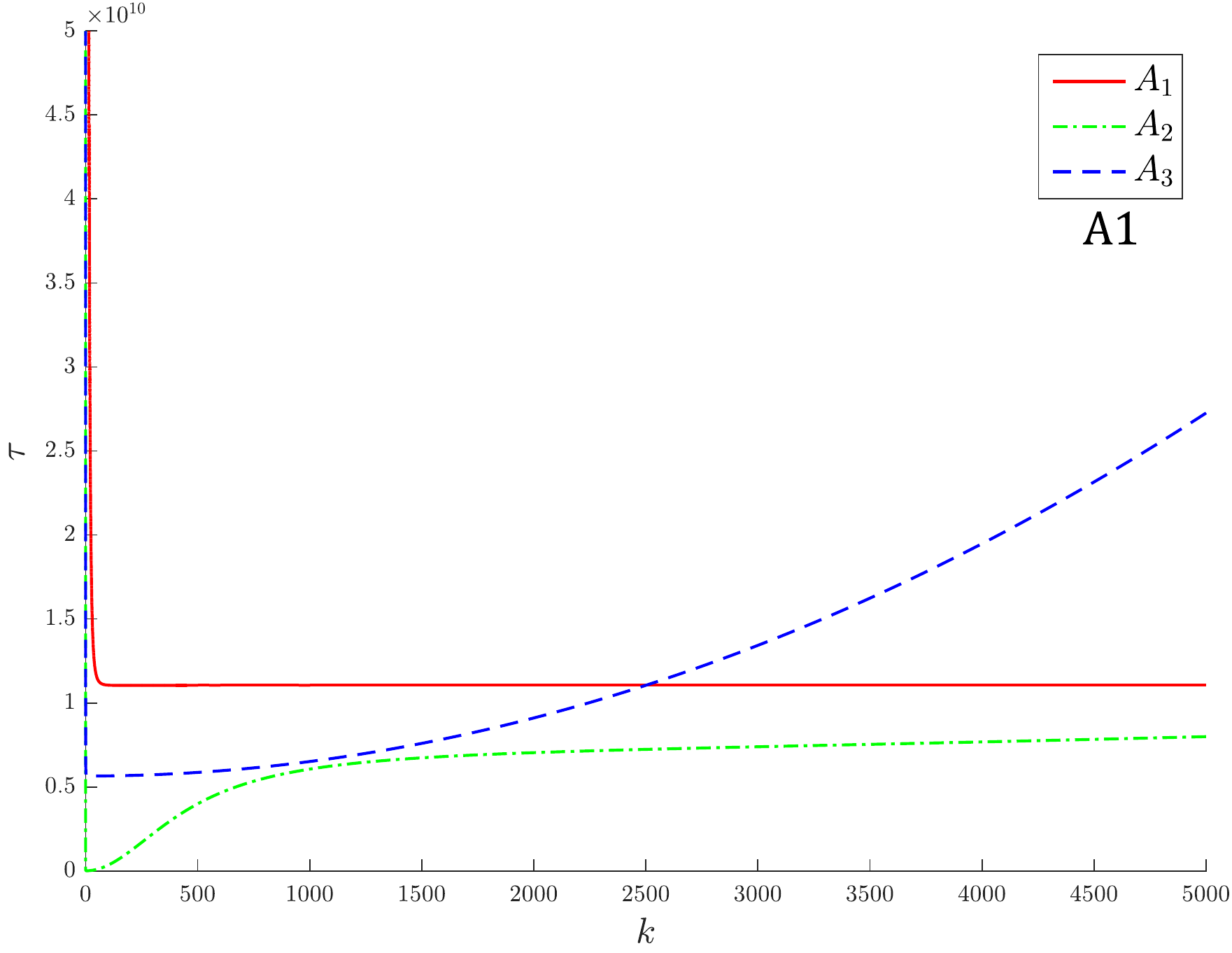}
\includegraphics[width=0.45\textwidth]{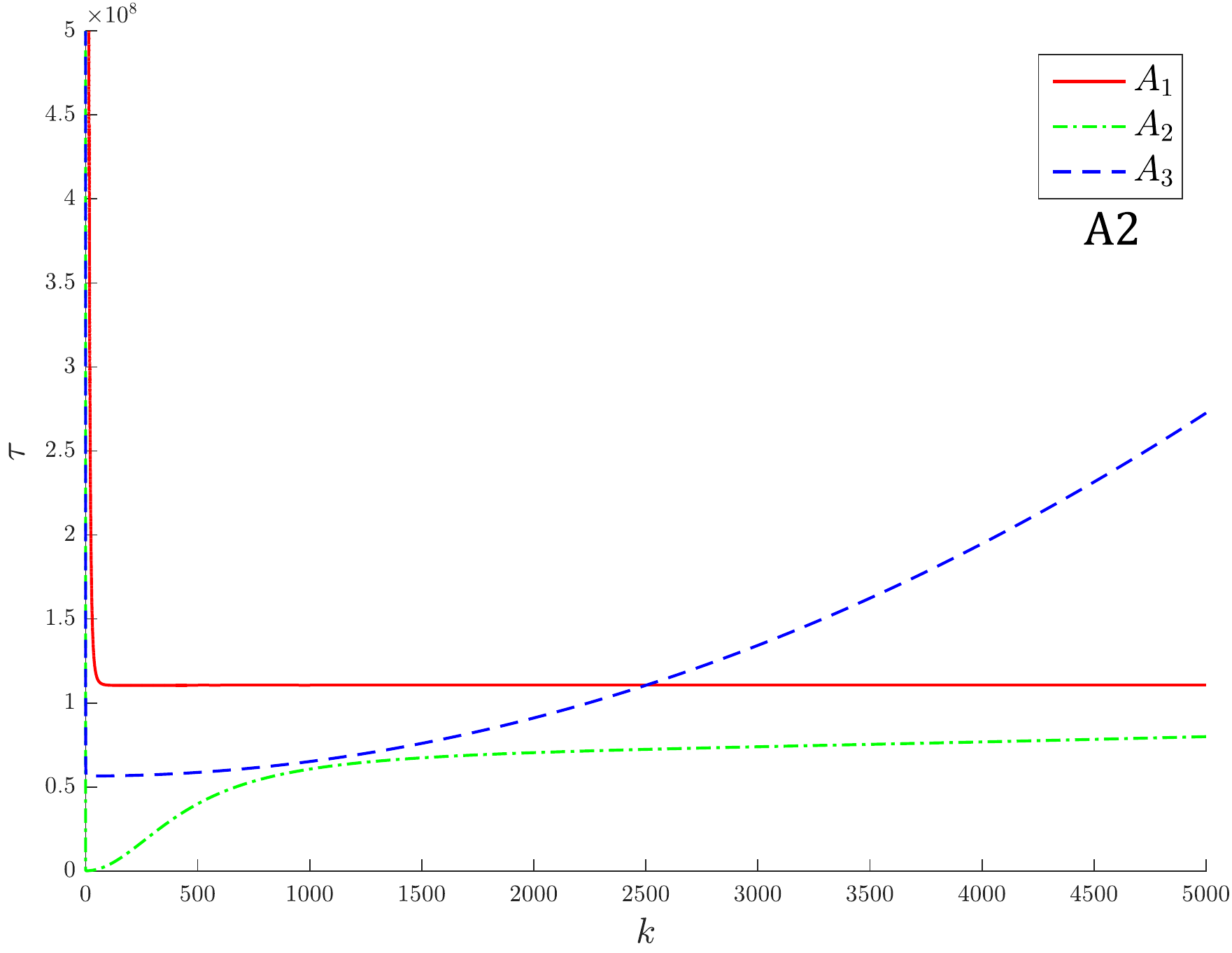}\\
\includegraphics[width=0.45\textwidth]{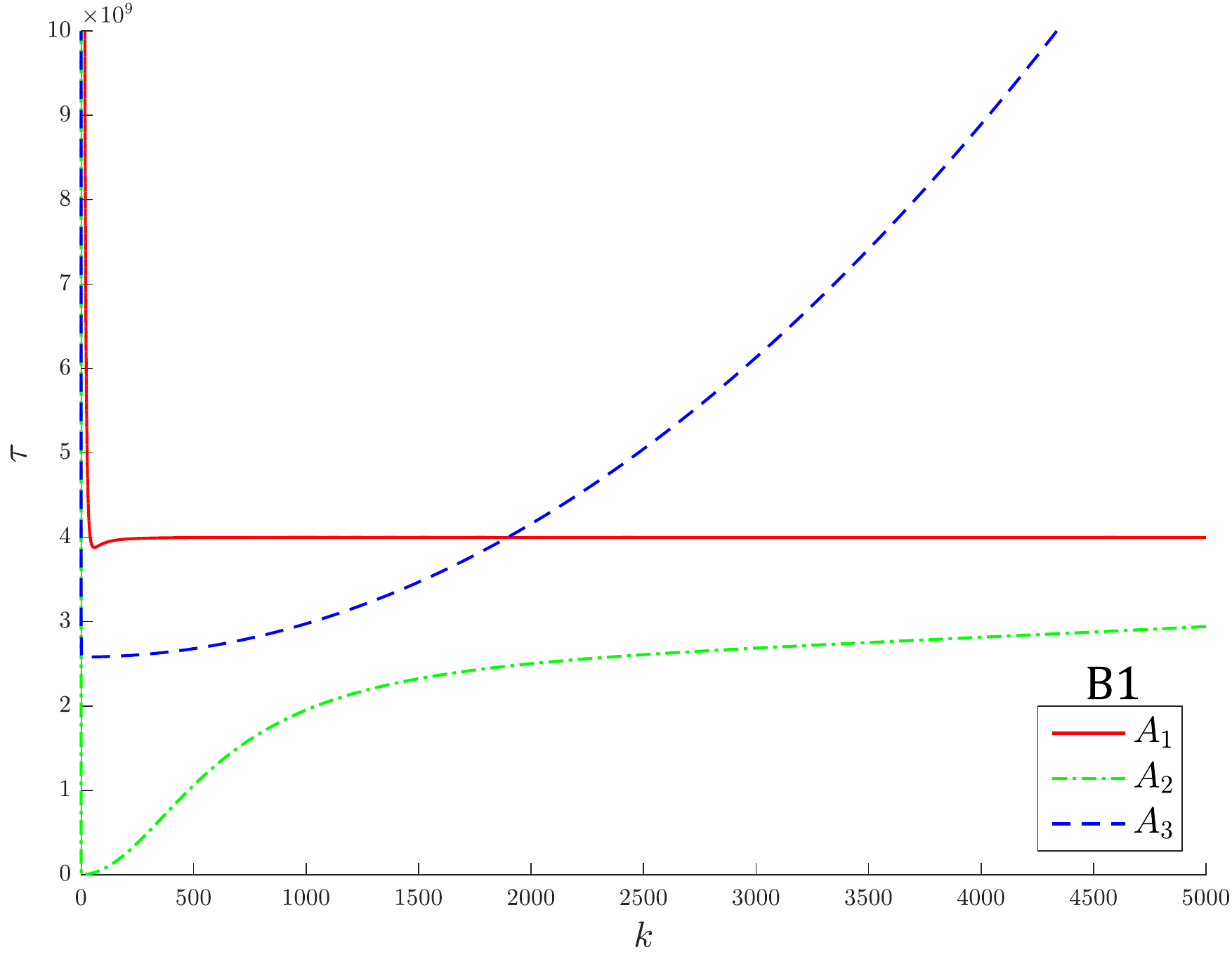}
\includegraphics[width=0.45\textwidth]{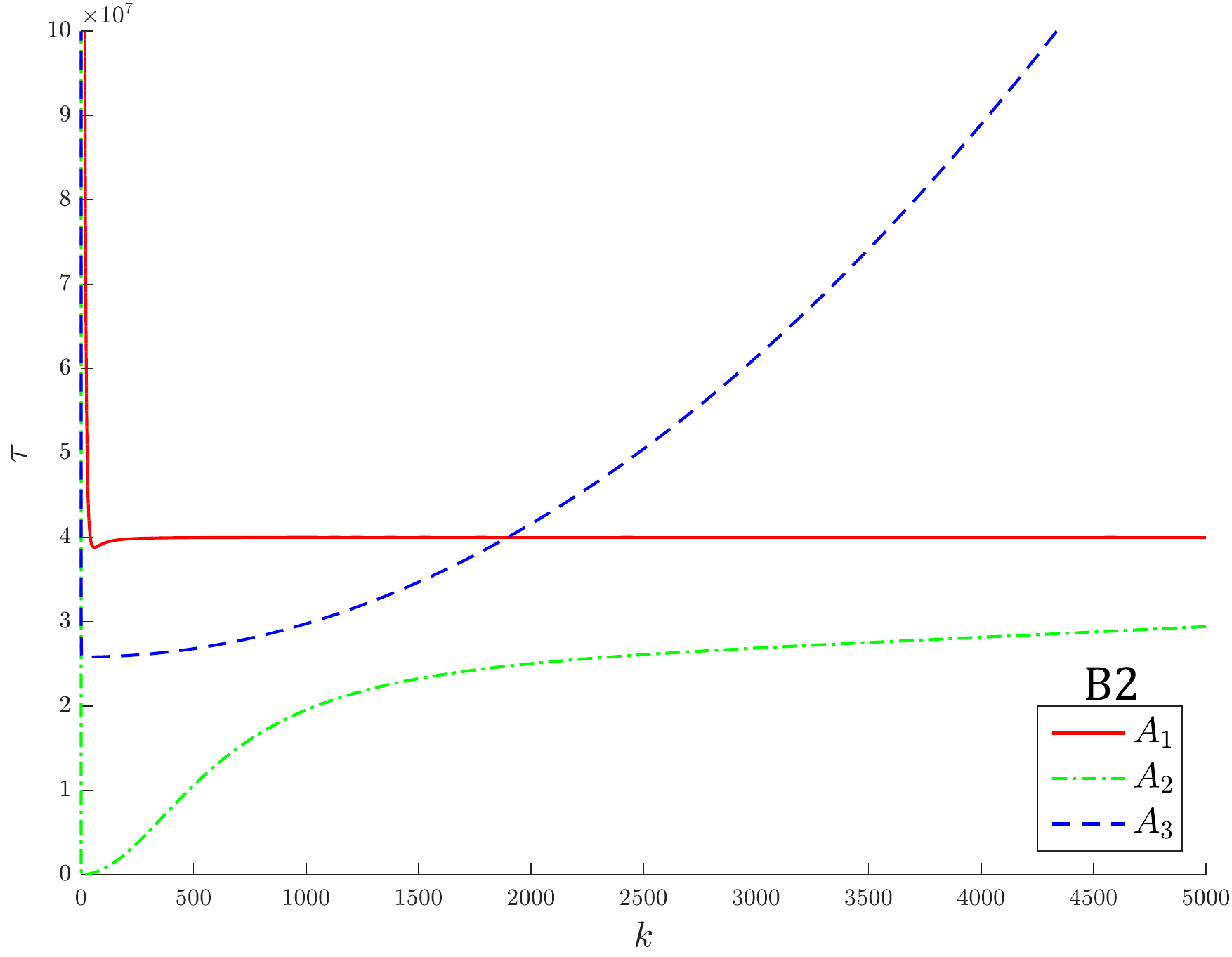}
\end{center}
\caption{Null level set of the characteristic polynomial coefficients $a_i$ with respect to $\tau$ defined in \eqref{eq:la-poly} for the general coupled system and for different coupling term $\btheta$. (A1) $\gamma=10^{-4}$ with $\btheta^{(1)}$. (A2) $\gamma=10^{-2}$ with $\btheta^{(1)}$. (B1-B2) Similar analysis for $\btheta^{(2)}$. In all plots \cblue{we use $\rho=1$, $\beta_2=0.6319$}.}
\label{fig:lin_ex8}
\end{figure}

\begin{figure}[!t]
\begin{center}
\includegraphics[width=0.325\textwidth]{Aif_ktau_DSet1}
\includegraphics[width=0.325\textwidth]{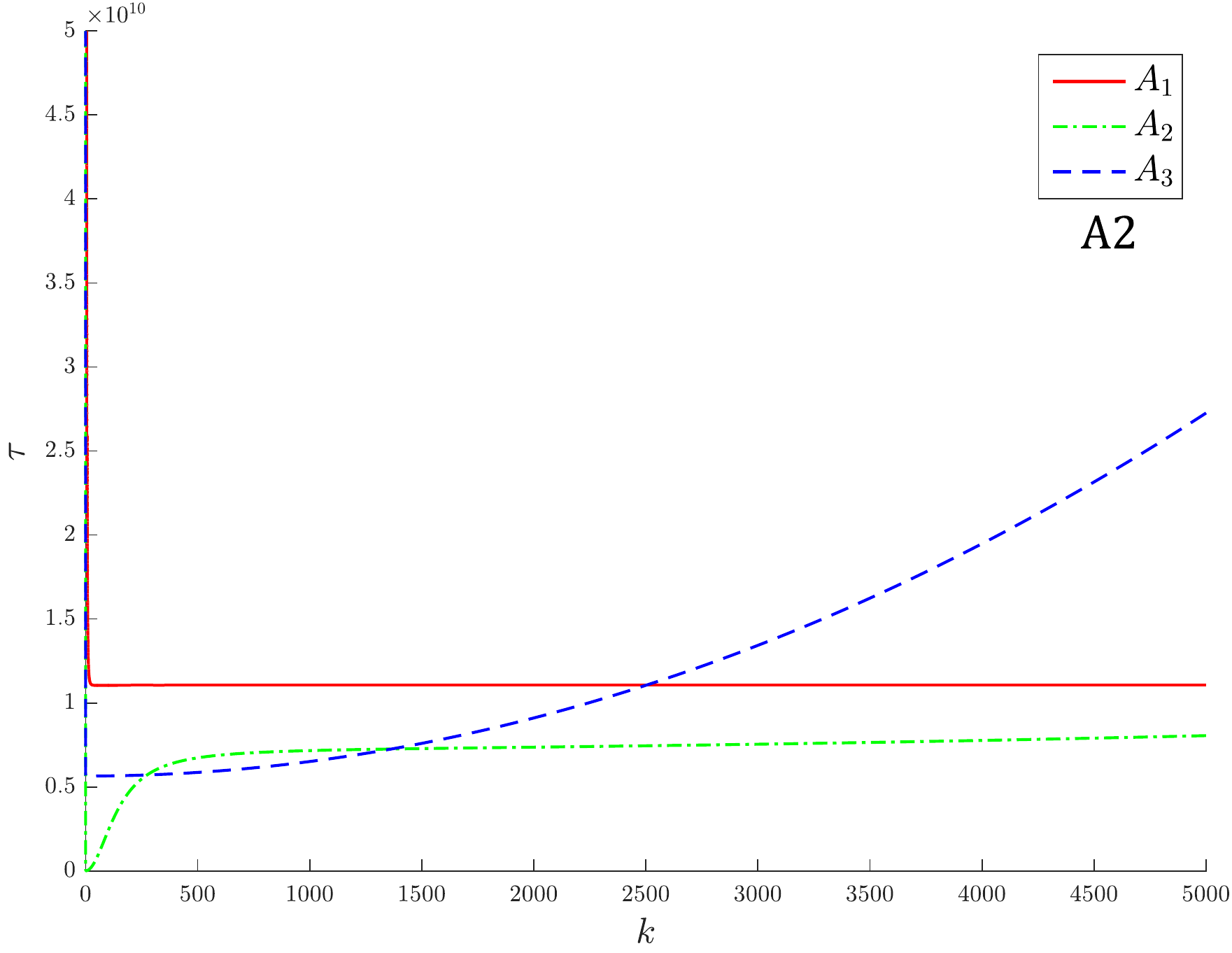}
\includegraphics[width=0.325\textwidth]{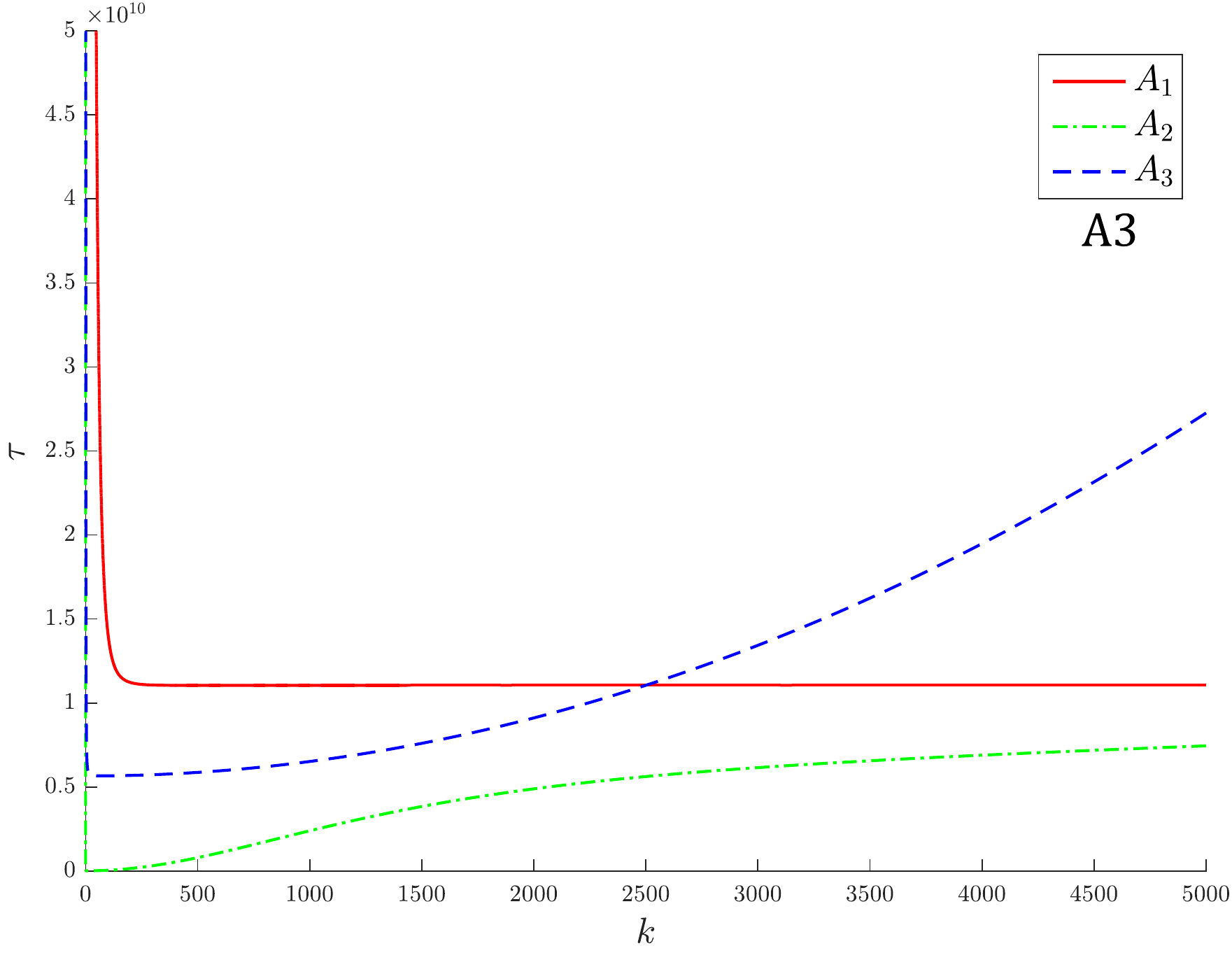}\\
\includegraphics[width=0.325\textwidth]{Aif_ktau_DSet1_MSetY}
\includegraphics[width=0.325\textwidth]{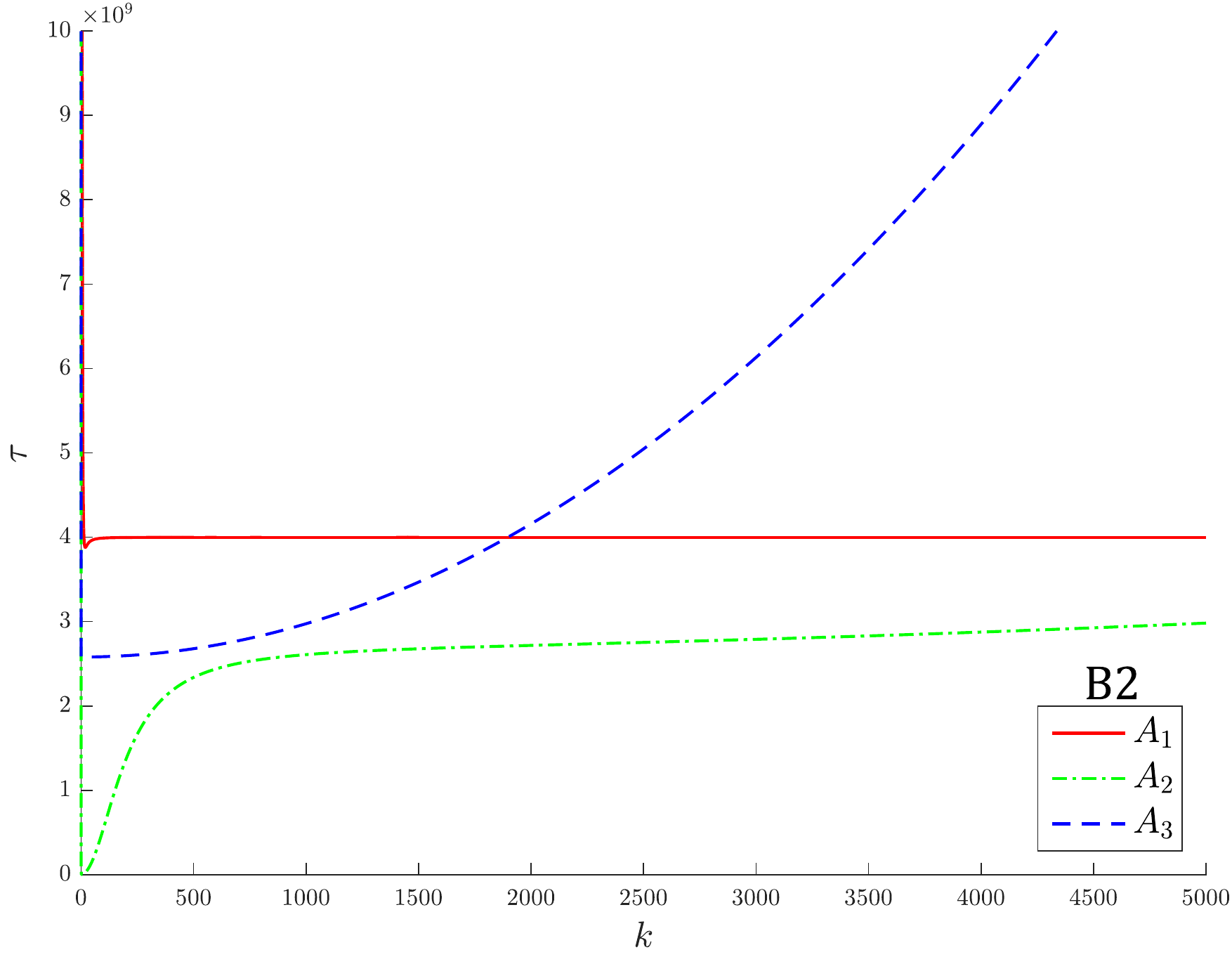}
\includegraphics[width=0.325\textwidth]{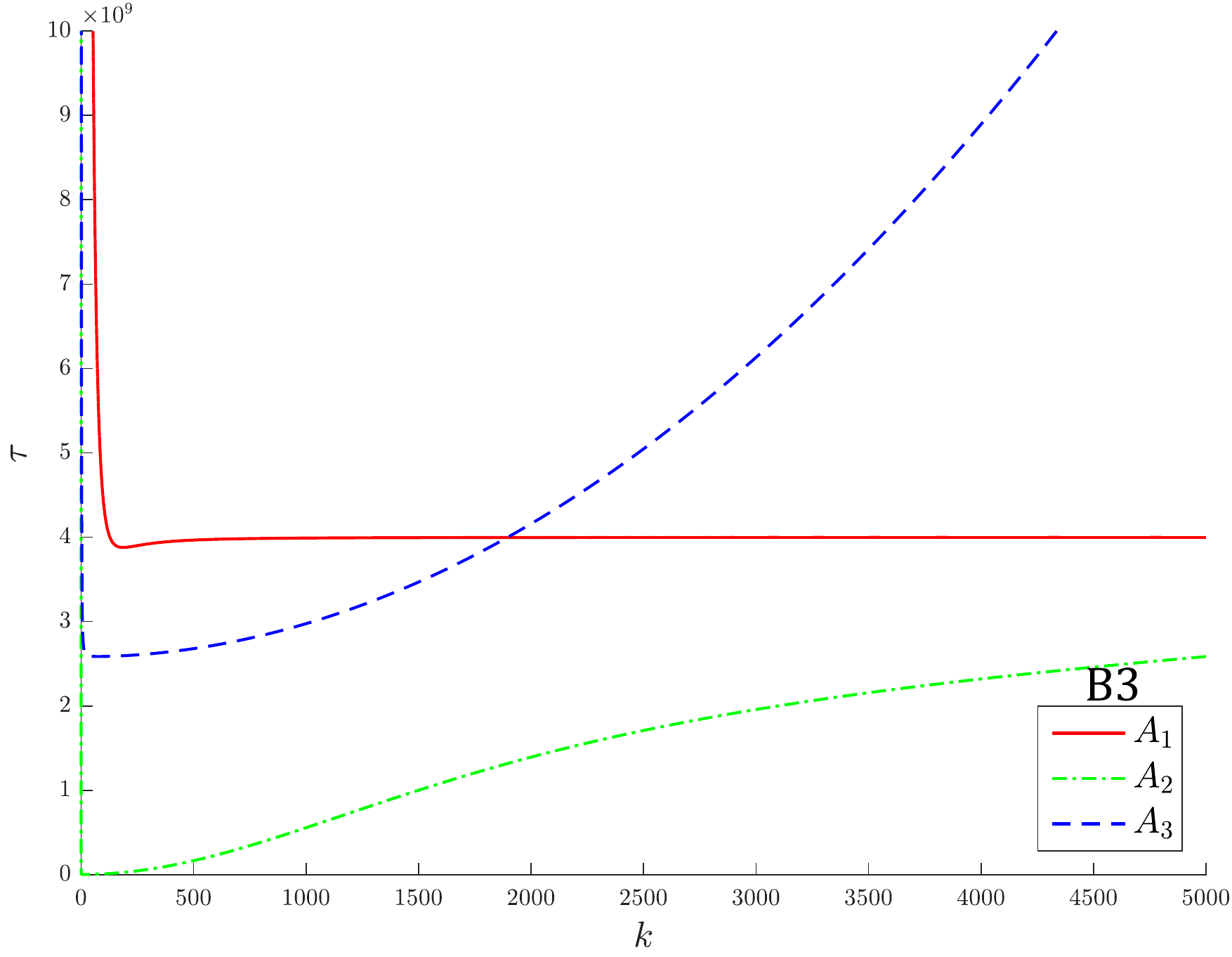}
\end{center}
\caption{\cblue{Null level set of the characteristic polynomial coefficients $a_i$ with respect to $\tau$ defined in \eqref{eq:la-poly} for the general coupled system and for different coupling term $\btheta$. (A1) $\beta_1=170$ with $\btheta^{(1)}$. (A2) $\beta_1=17$ with $\btheta^{(1)}$. (A3) $\beta_1=1700$ with $\btheta^{(1)}$. (B1-B3) Similar analysis for $\btheta^{(2)}$. In all plots we use $\rho=1$ and $\gamma = 10^{-4}$.}}
\label{fig:lin_ex9}
\end{figure}

\cblue{In a more general scenario, we look at how the strength of the coupling between the poromechanical and chemical systems modifies the} linear stability properties. Figure~\ref{fig:lin_ex8} \cblue{presents}, for different values of \cblue{the coupling parameter} $\gamma$, the null levels of the dispersion relation coefficients. \cblue{In order to reduce the complexity of the stability analysis, we fix $\beta_2=0.6319$ letting the other parameter values unchanged. This modification forces the coefficient $a_0$ to be strictly positive, reducing the complex pattern generated by the Routh-Hurwitz conditions (cf. \eqref{eq:rh-cond1}-\eqref{eq:rh-cond4}).} Increasing $\gamma$ \cblue{reduces the critical} value of $\tau$ needed 
to reach instability (compare with Figs.~\ref{fig:lin_ex8}(A1), (A2)) without affecting the pattern generated by the null level-set. \cblue{Additionally, the source terms in the modified Schnackenberg model depend also on $\beta_1$. In order to compare the effectiveness of both parameters in stabilising or destabilising the system, we present in Figure \ref{fig:lin_ex9} the null level-set of the coefficients $a_i$ for different values of $\beta_1$. In that case, we see that the critical value for $\tau$ is affected by increasing $\beta_1$ (see Fig.~\eqref{fig:lin_ex9}(A1)-(A3)). Moreover, and contrary to what occurs with $\gamma$, it seems to affect significantly the pattern of the null level-sets (especially the one associated with the coefficient $a_2$) and therefore also the spatial scale where instabilities can occur}. Analogous conclusions can be drawn even using a nonlinear coupling function $r^{(2)}(\bw)$. As in the previous scenario, the different coefficients present a large interval where the Routh-Hurwitz conditions are not satisfied. In summary, the coupled system \cblue{is able to produce} non-trivial patterns at very different length scales.

\begin{figure}[!t]
\begin{center}
\includegraphics[width=0.45\textwidth]{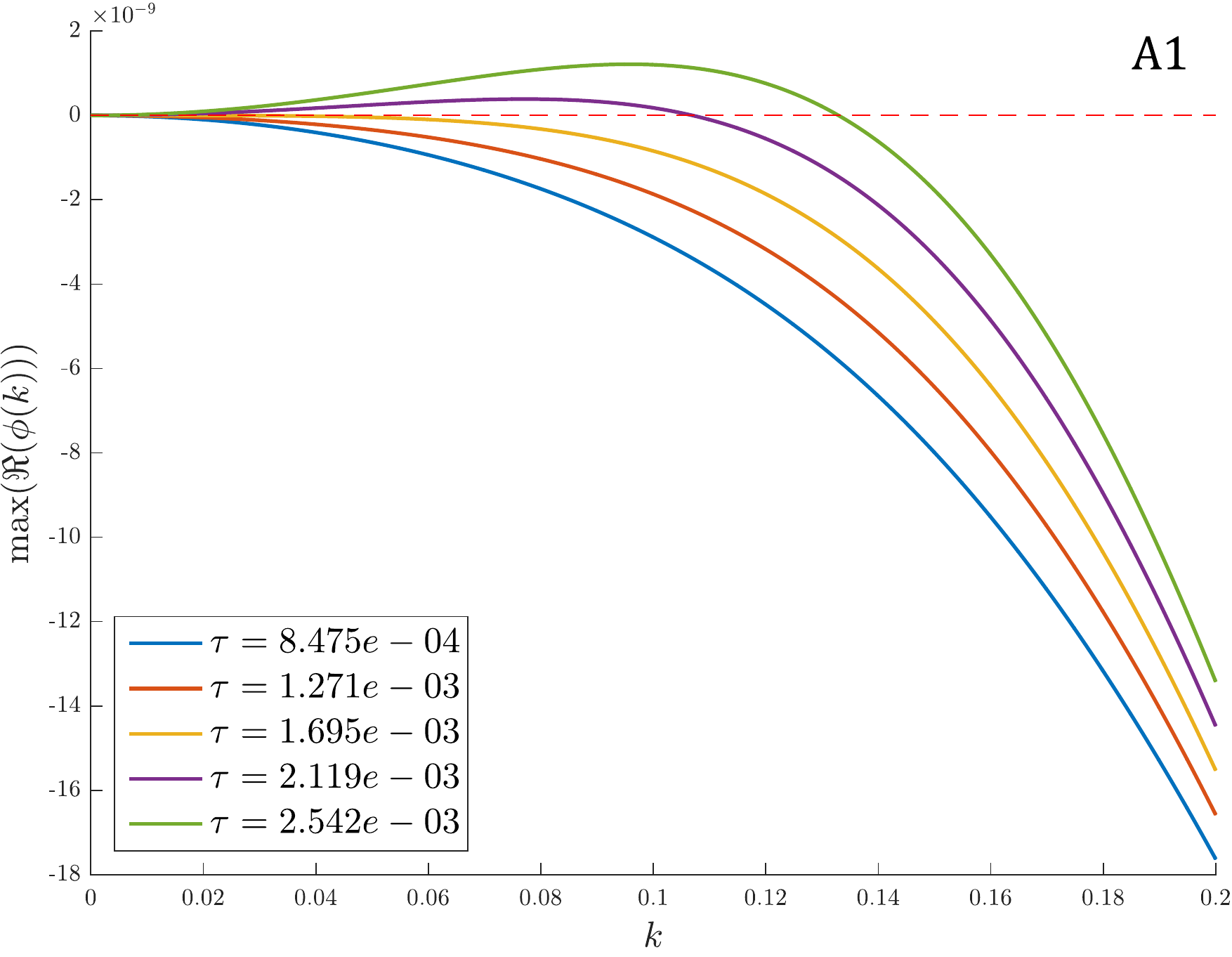}
\includegraphics[width=0.45\textwidth]{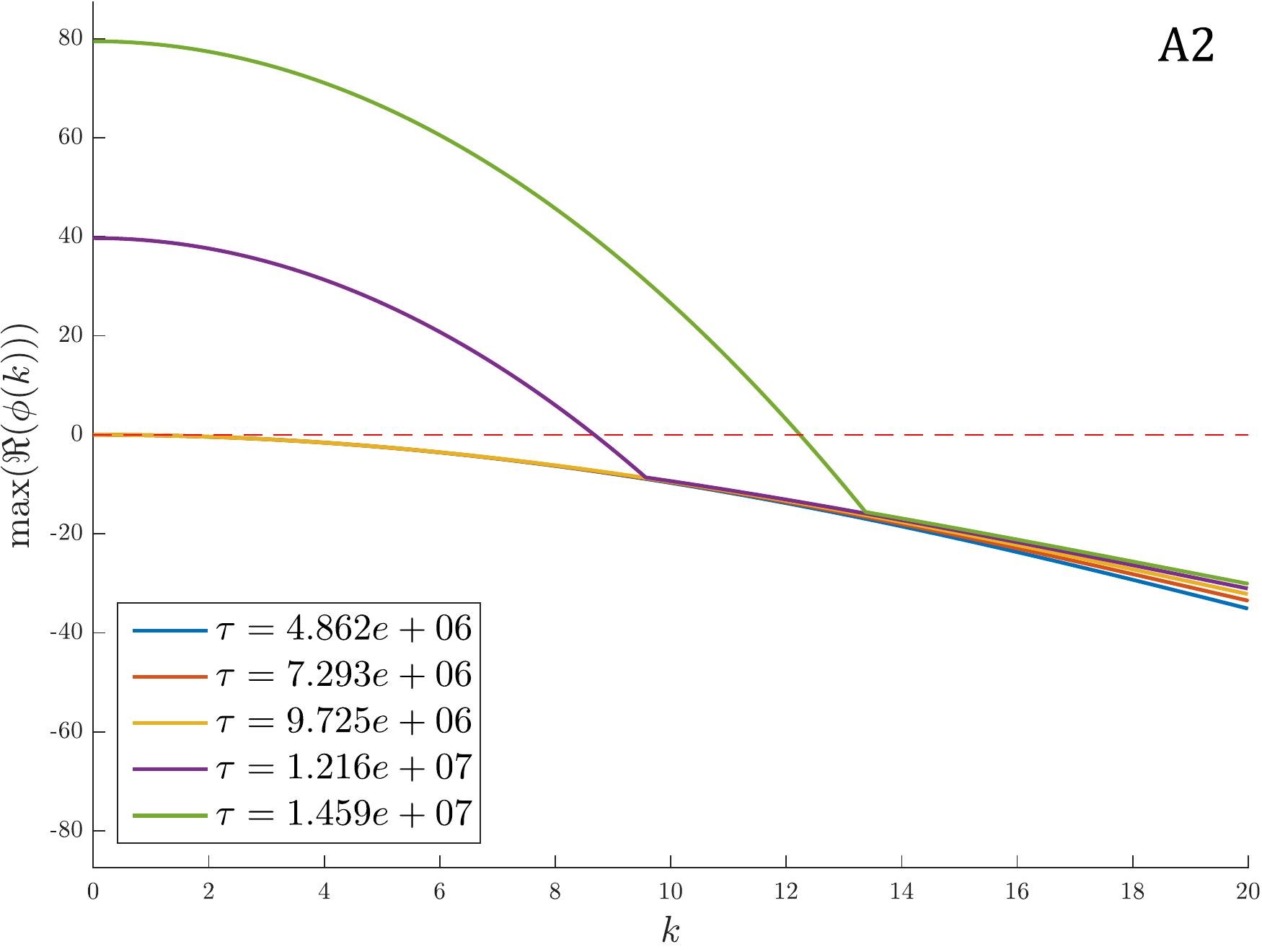}\\
\includegraphics[width=0.45\textwidth]{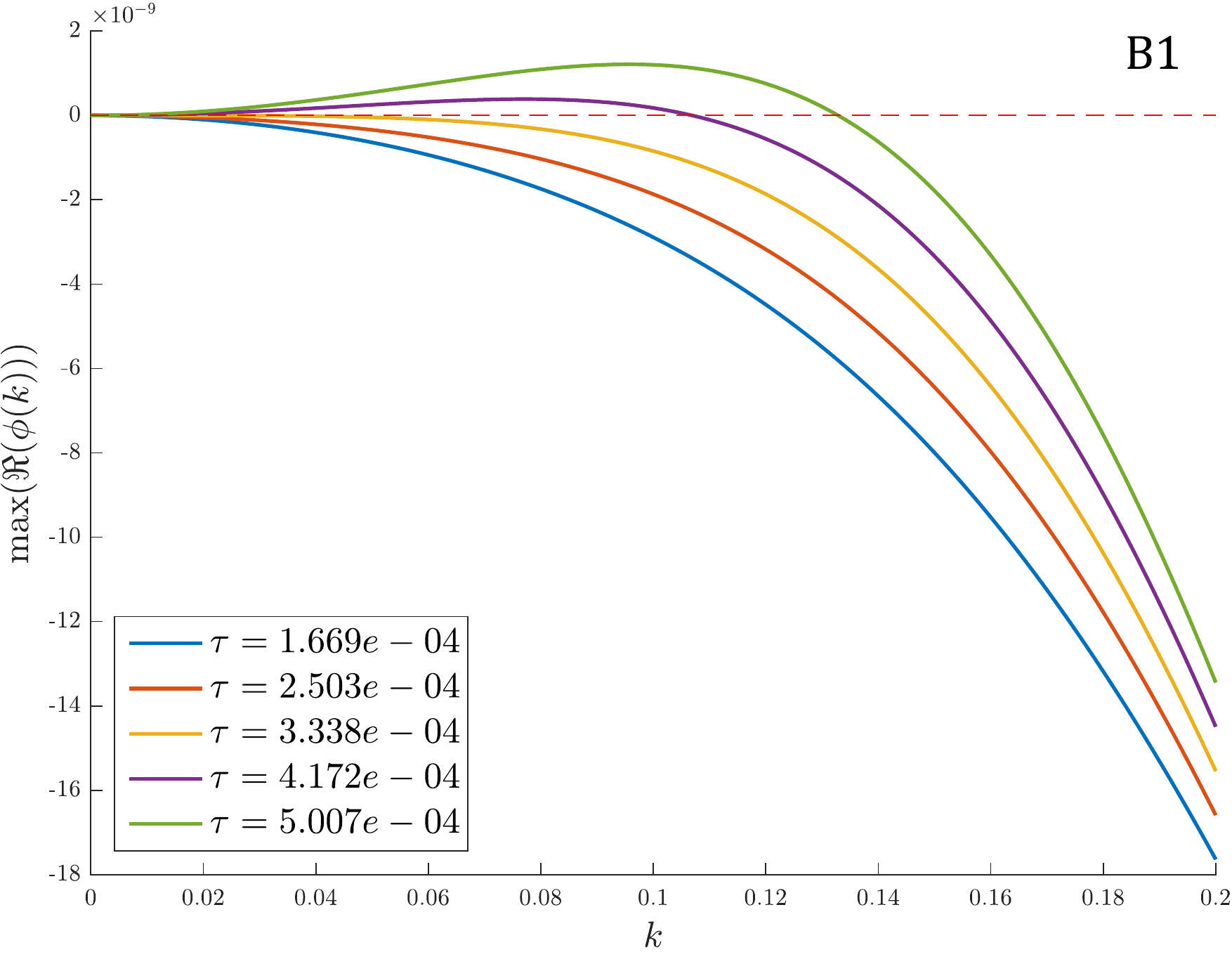}
\includegraphics[width=0.45\textwidth]{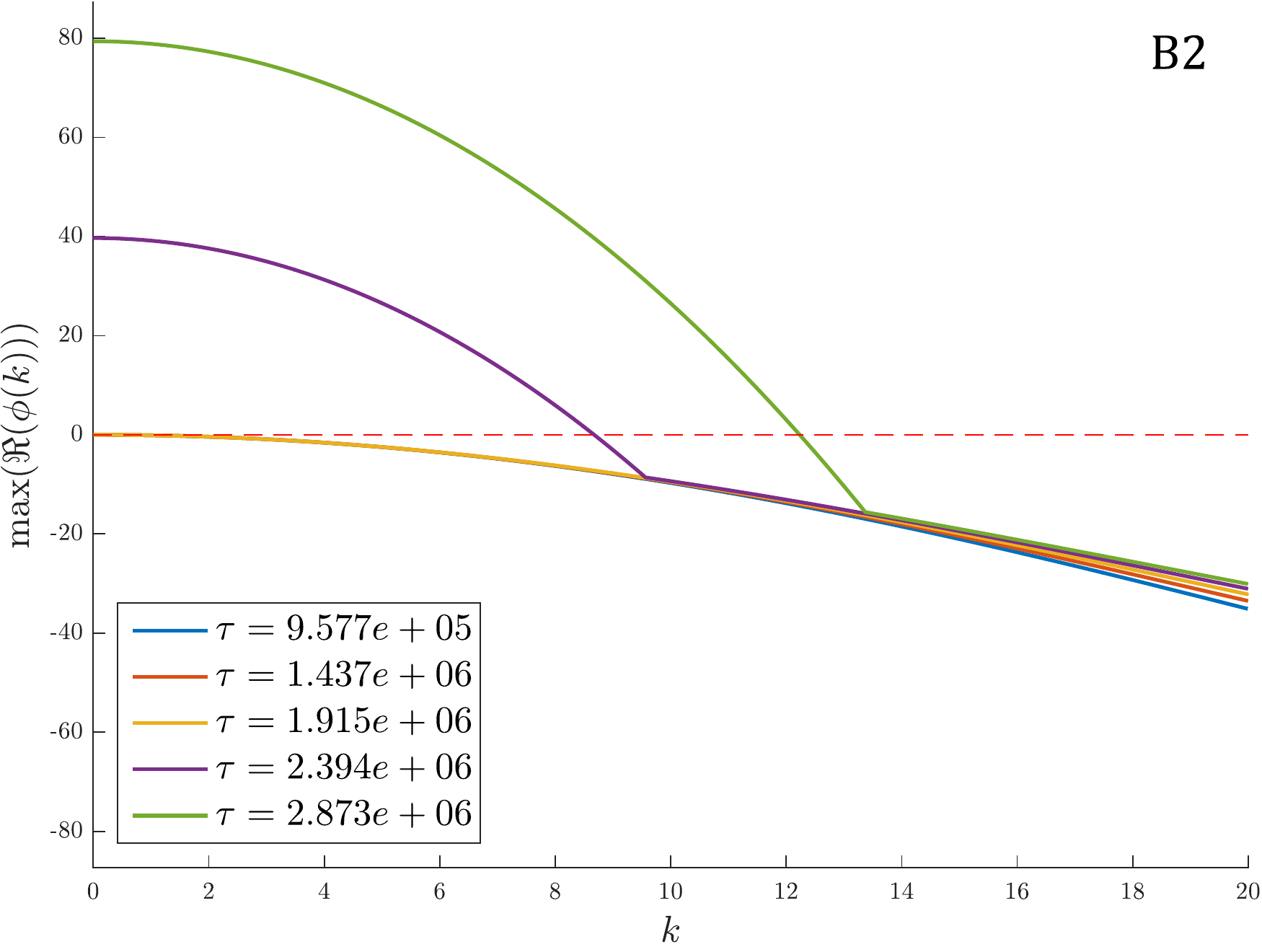}
\end{center}
\caption{\cblue{Dispersion relations for the general coupled poromechano-chemical model for different coupling terms $\btheta$. (A1) Dispersion relation associated to $\rho=1$ with $\btheta^{(1)}$. (A2) Dispersion relation associated to $\rho=0$ with $\btheta^{(1)}$. (B1-B2) Similar analysis with $\btheta^{(2)}$. In all plots, curves are drawn from the critical value $\tau_c$ (yellow) and for 25\% and 50\% increase/decrease of the parameter values, additionally we use $\beta_2=0.6319$ and $\gamma = 10^{-4}$.}}
\label{fig:lin_ex10}
\end{figure}	

\cblue{Another way to show pattern generation is through dispersion relations related to the parameter $\tau$, see Figure \ref{fig:lin_ex10}. These relations are computed using the implicit functions obtained from the coefficient $a_2$ when $\rho=0$, and from condition \eqref{eq:rh-cond4}, otherwise. The choice was made as they are the first coefficients to break the Routh-Hurwitz conditions with respect to $\tau$ in the explored space scale. We observe that if acceleration is present, a relatively small value of $\tau$  is sufficient (at a very low space scale) to induce instability, as seen in Figure \ref{fig:lin_ex10}(A1) (or (B1)). When $\rho=0$, the system can still reach instability but it needs a larger value of $\tau$, which shows again the effect of inertia on the system stability. Comparing Figures \ref{fig:lin_ex10}(A1)-(A2) with \ref{fig:lin_ex10}(B1)-(B2) we see that, for the selected parameters, the choice of the coupling function $\btheta$ affects the critical value of $\tau$, but not the pattern of the dispersion relations.}

\begin{figure}[!t]
\begin{center}
\includegraphics[width=0.325\textwidth]{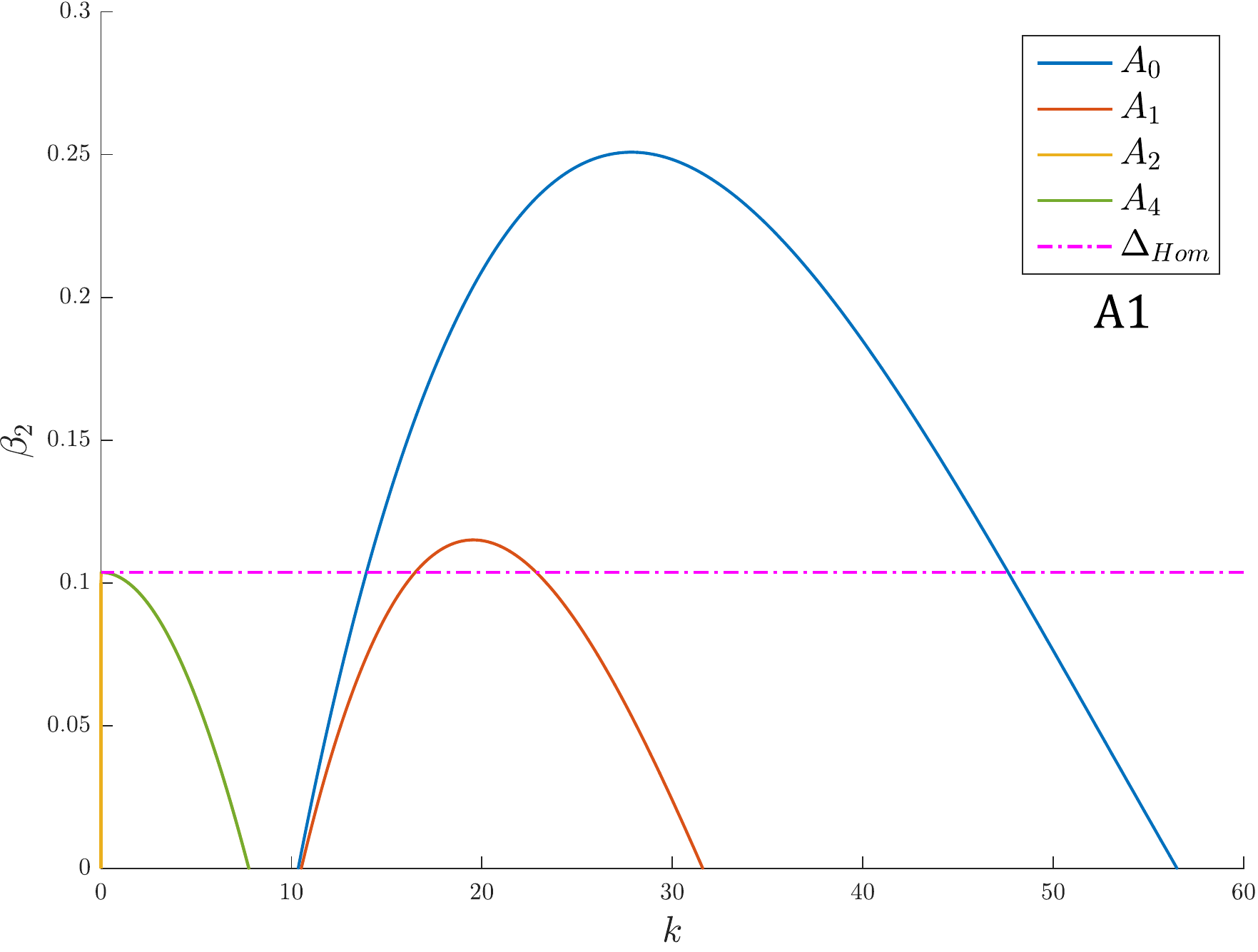}
\includegraphics[width=0.325\textwidth]{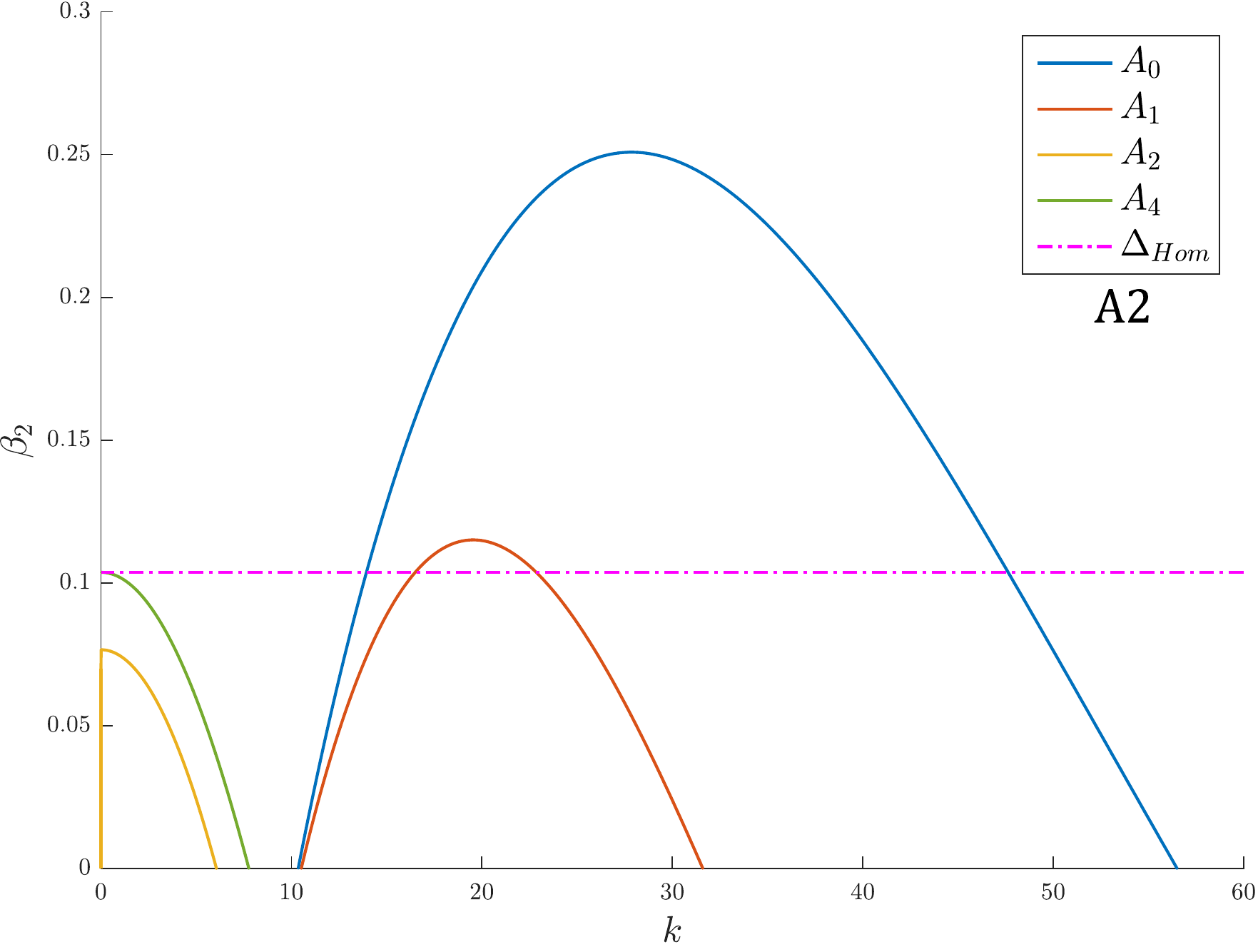}
\includegraphics[width=0.325\textwidth]{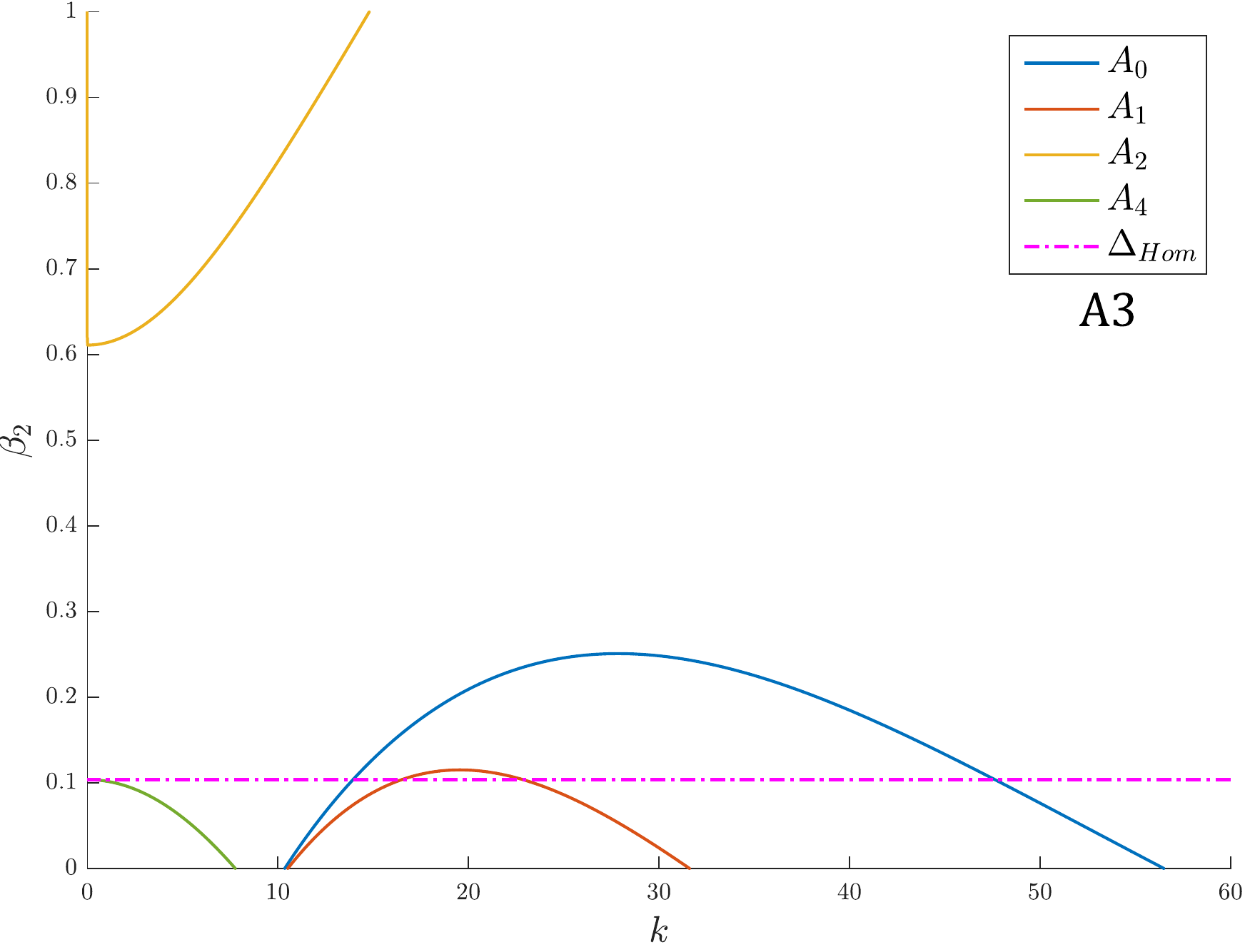}\\
\includegraphics[width=0.325\textwidth]{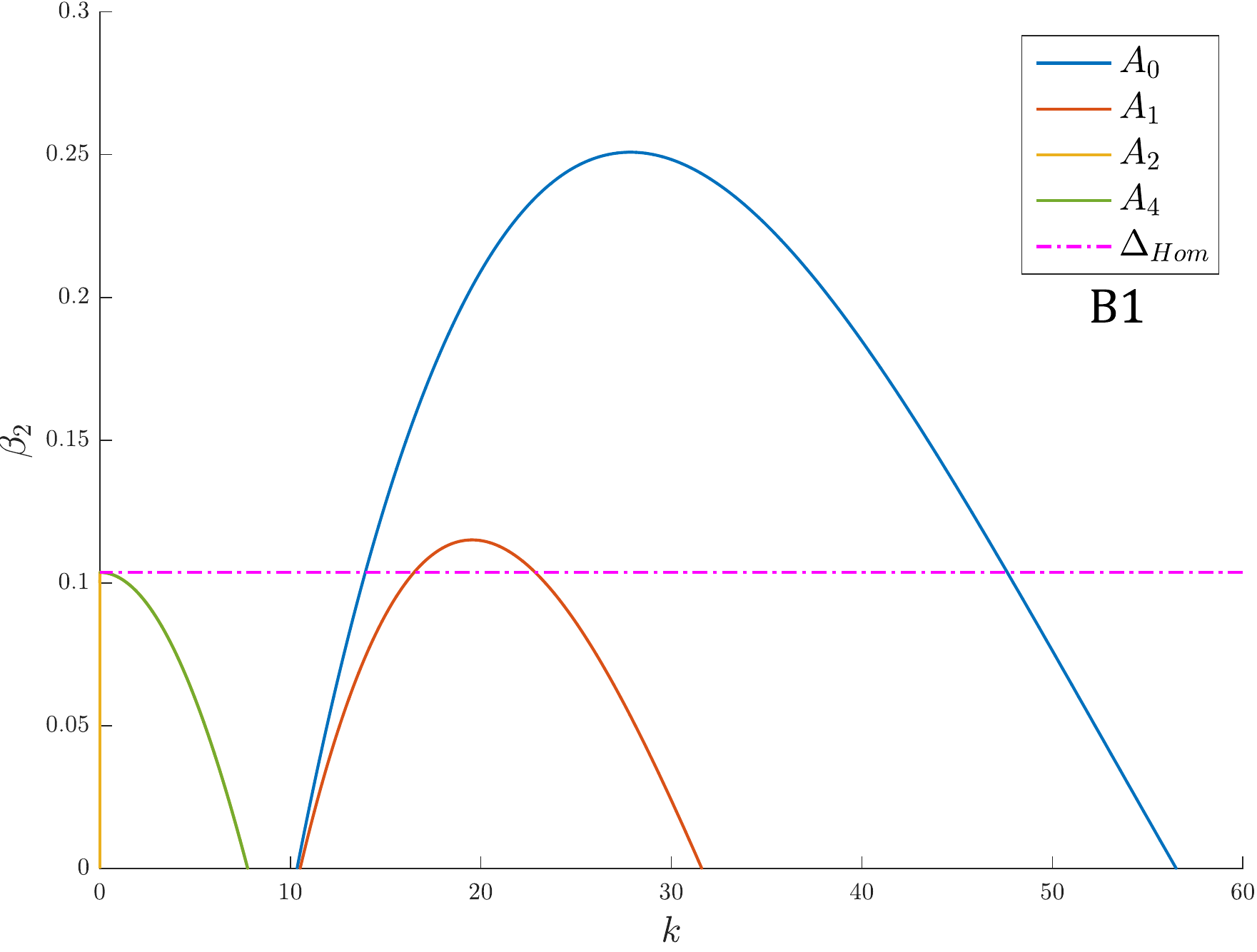}
\includegraphics[width=0.325\textwidth]{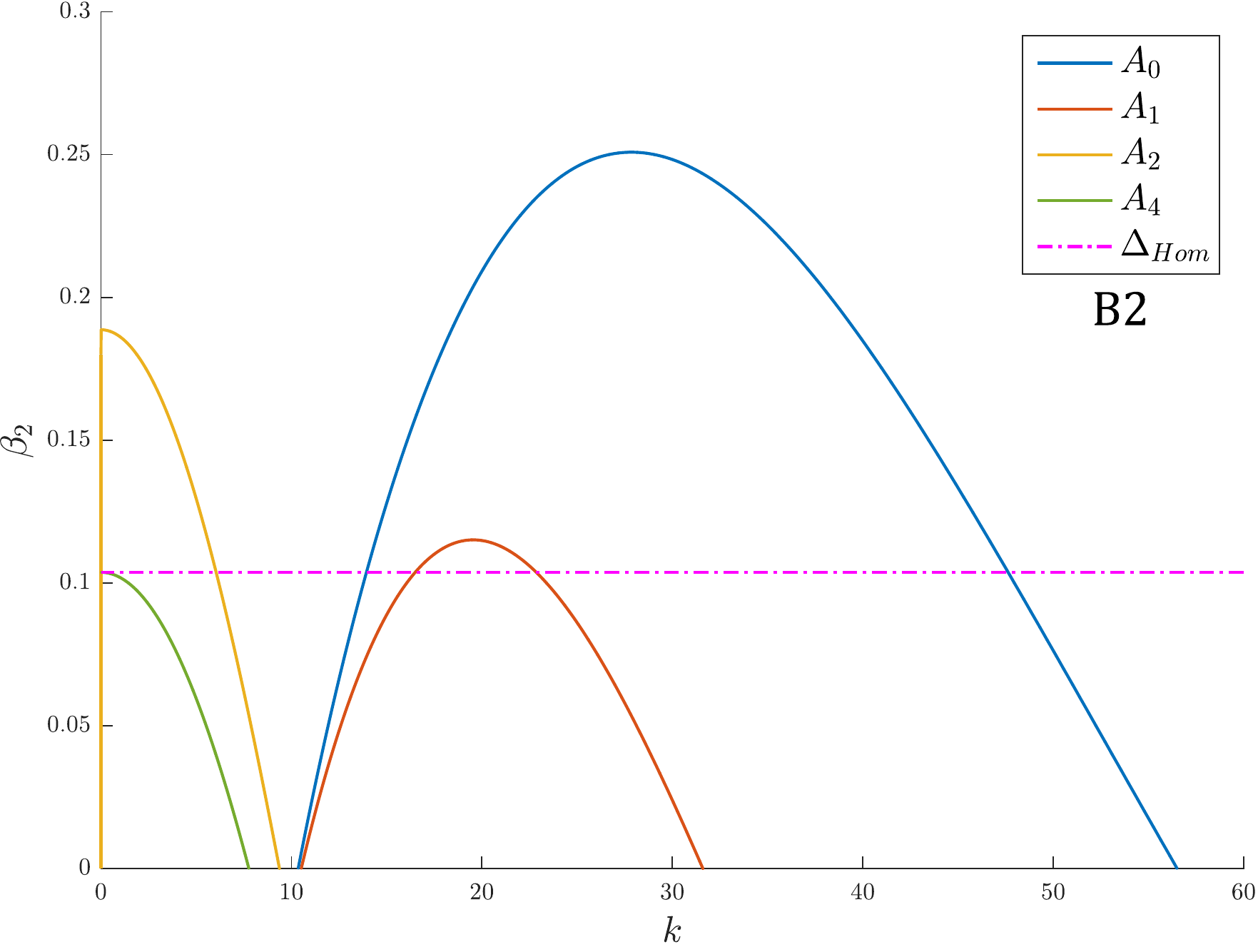}
\includegraphics[width=0.325\textwidth]{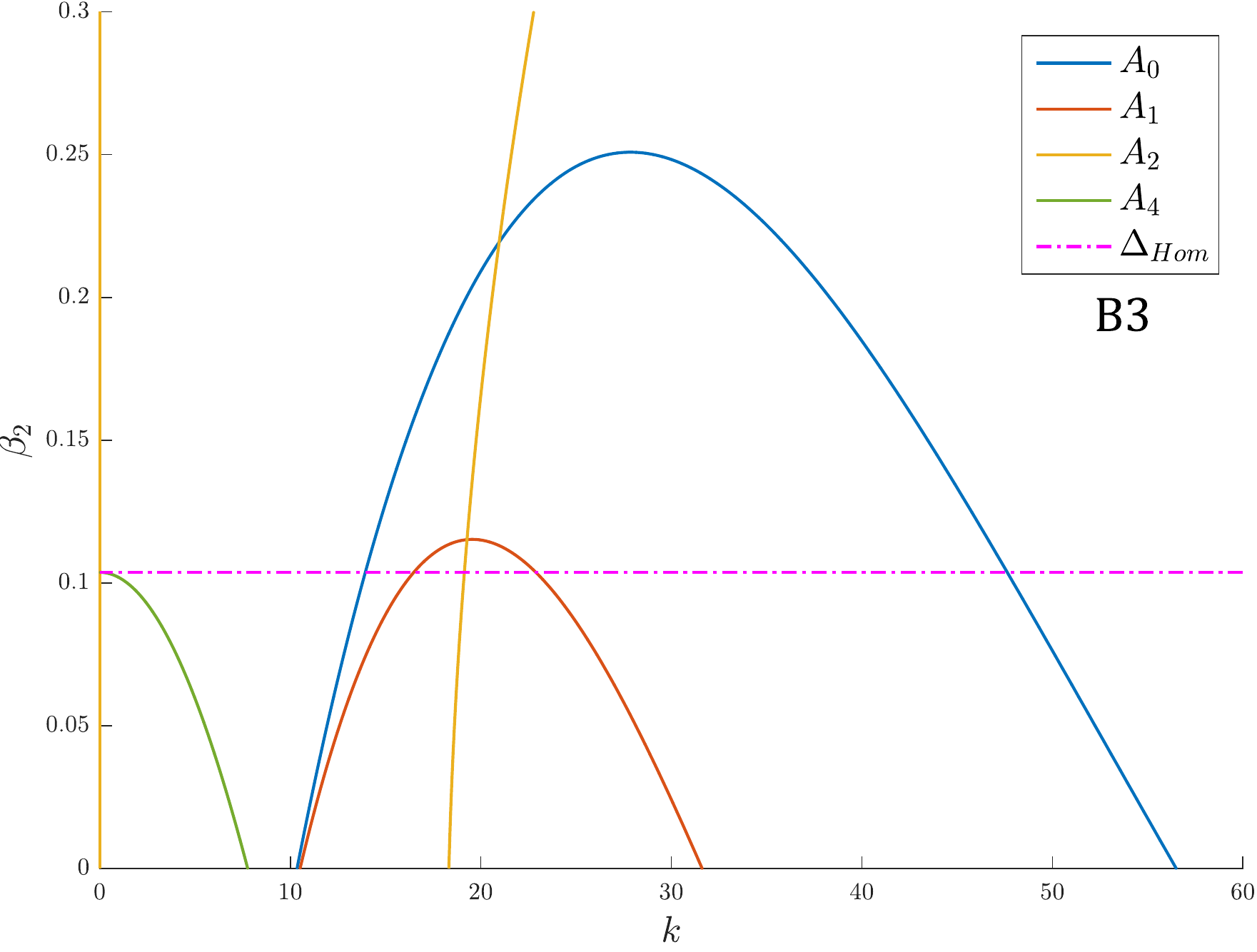}
\end{center}
\caption{Null level set of the characteristic polynomial coefficients $a_i$ with respect to $\beta_2$ defined in \eqref{eq:la-poly} for the general coupled system and for different coupling term $\btheta$. (A1) $\tau=100$, $\gamma=10^{-4}$,  with $\btheta^{(1)}$. (A2) $\tau=10^4$, $\gamma=10^{-2}$ with $\btheta^{(1)}$. (A3) $\tau=10^5$, $\gamma=10^{-2}$ with $\btheta^{(1)}$. (B1-B3) Similar analysis with $\btheta^{(2)}$. In all plots \cblue{we use} $\rho=1$.}
\label{fig:lin_ex11}
\end{figure}

To conclude this section, we analyse the sign of \cblue{each $a_i$} against the value of $\beta_2$, and compare with the results from 
the uncoupled system studied above. Figure~\ref{fig:lin_ex11} displays the null level set \cblue{when increasing the coupling strength (realised by augmenting $\tau$ and $\gamma$ simultaneously). From Figures} \ref{fig:lin_ex11}(A1) to \ref{fig:lin_ex11}(A3), we observe that $a_2$ is particularly affected by the choice of these two parameters. First, reducing the region of instability (see Fig.~\ref{fig:lin_ex11}(A2)) we obtain, in a strongly coupled system, two instability regions delimited by the red ($a_0$) and green ($a_2$) level sets that overlap with respect to the wave number $k$, but not with respect to $\beta_2$. Consequently, and depending on the strength of the coupling between reaction-diffusion and 
poroelastic effects, we can discriminate different values of $\beta_2$ that will produce distinct patterns at desired specific scales. This is more clearly seen for the case with $\btheta=\btheta^{(2)}$, where the interval of $\beta_2$ that leads to instability increases rapidly with the augmentation of $\tau$ and $\gamma$, and tends \cblue{to finally overlap with} the region delimited by $a_0$ (see red level in Fig.~\ref{fig:lin_ex11}(B3)).

	\section{Numerical method and implementation}\label{sec:FE}
	
In \cite{verma01} we propose a discretisation in space using a mixed finite element method. 
In its lowest-order form, the method consists of piecewise bilinear elements enriched with 
bubbles for the \cblue{displacements (cubic for triangles, and quartic for tetrahedra)}, piecewise linear and continuous approximations for the fluid pressure and for the chemical solutes; and piecewise constant approximation for total pressure. 
For sake of completeness we outline here the spatio-temporal method. 

The time discretisation is achieved by a backward Euler 
scheme and an implicit centred difference method for the first and second order time derivatives, respectively. Denoting 
$$\delta_t X^{n+1}:=\frac{X^{n+1}-X^{n}}{\Delta t}, \qquad \delta_{tt} X^{n+1}:= \frac{X^{n+1}-2X^n+X^{n-1}}{\Delta t^2},$$ 
the fully discrete method  reads: From initial data \cblue{$\bu_h^{s,0},p_h^{f,0},\psi_h^{0}, w_{1,h}^0, w_{2,h}^0$} 
(which will be projections of the continuous initial conditions of each field) and
for $n=1,\ldots$, find $\bu_h^{s,n+1}\in\bV_h,p_h^{f,n+1}\in Q_h,\psi_h^{n+1}\in Z_h,w_{1,h}^{n+1}\in W_h,w_{2,h}^{n+1}\in W_h$ such that
\begin{alignat}{5}
\widetilde{a}_1(\bu_h^{s,n+1},\bv_h^s)  \;\,&+&   a_1(\bu_h^{s,n+1},\bv_h^s)   &&                 &\;+&\; b_1(\bv_h^s,\psi_h^{n+1})     &=&\;F_{r_h^{n+1}}(\bv_h^s)&\; \;\forall \bv_h^s\in\bV_h, \label{weak-u-h}\\
\tilde{a}_2(p_h^{f,n+1},q_h^f)   &\;+&               &&      a_2(p_h^{f,n+1},q_h^f)   &\;-&\;   {b}_2(q_h^f,\delta_t\psi_h^{n+1})  &=&\; G_{\ell^{n+1}}(q_h^f) &\; \; \forall q_h^f\in Q_h, \label{weak-p-h}\\
&&b_1(\bu_h^{s,n+1},\phi_h)  &\;+\;& b_2(p_h^{f,n+1},\phi_h)&\;-&\; a_3(\psi_h^{n+1},\phi_h) &=&\; 0 &\; \; \forall\phi_h\in Z_h, \label{weak-psi-h}    \\
\tilde{a}_4(w_{1,h}^{n+1},s_{1,h}) &\;+ &\; a_4(w_{1,h}^{n+1},s_{1,h}) &&& & & = &\; J_{f_h^{n+1}}(s_{1,h}) &\; \; \forall s_{1,h} \in W_h, \label{weak-w1-h}\\
\tilde{a}_5(w_{2,h}^{n+1},s_{2,h}) &\;+ &\; a_5(w_{2,h}^{n+1},s_{2,h}) &&&&& =&\; J_{g_h^{n+1}}(s_{2,h}) &\; \; \forall s_{2,h}\in W_h\cblue{,} \label{weak-w2-h}
\end{alignat}
with the bilinear forms and linear functionals defined as 
\begin{gather}
\widetilde{a}_1(\bu_h^{s,n+1},\bv_h^s):=\int_{\Omega}\delta_{tt}\bu^{s,n+1}_h\cdot\bv_h ,\;\;a_1(\bu^{s,n+1},\bv^s) := 2\mu \int_{\Omega} \beps(\bu_h^{s,n+1}):\beps(\bv_h^s),\nonumber\\
\;\;b_1(\bv_h^s,\phi_h):= -\int_{\Omega}\phi_h\vdiv \bv_h^s,\;\;\, b_2(p_h^{f,n+1},\phi_h):= \frac{\alpha}{\lambda} \int_{\Omega} p_h^{f,n+1}\phi_h,\;\;a_3(\psi_h^{n+1},\phi_h):=  \frac{1}{\lambda} \int_{\Omega} \psi_h^{n+1}\phi_h,\nonumber\\
\tilde{a}_2(p_h^{f,n+1},q_h^f)  := \biggl(c_0+\frac{\alpha^2}{\lambda}\biggr)
\int_{\Omega}  \delta_t p_h^{f,n+1} q_h^f, \quad a_2(p_h^{f,n+1},q_h^f)  := \frac{1}{\eta}\int_{\Omega}\kappa \nabla p_h^{f,n+1}\cdot\nabla q_h^f,\nonumber\\
\tilde{a}_4(w_{1,h}^{n+1},s_1) : = \delta_t \int_\Omega  w_{1,h}^{n+1}s_{1,h}, \quad a_4(w_{1,h}^{n+1},s_{1,h}) : =  \int_\Omega D_1(\bx)\nabla w_{1,h}^{n+1}\cdot \nabla s_{1,h}, \label{eq:forms}\\
\tilde{a}_5(w_{2,h}^{n+1},s_{2,h})  : = \delta_t\int_\Omega  w_{2,h}^{n+1}s_{2,h}, \quad a_5(w_{2,h}^{n+1},s_{2,h})  : = \int_\Omega D_2(\bx)\nabla w_{2,h}^{n+1}\cdot \nabla s_{2,h},\nonumber \\
F_{r_h^{n+1}}(\bv_h^s) := \rho\int_{\Omega} \bb^{n+1} \cdot\bv_h^s
+\tau\int_{\Omega} r_h^{n+1}  \bk\otimes\bk :\beps(\bv_h^s), \nonumber \\
G_{\ell^{n+1}}(q_h^f) :=  \int_{\Omega} \ell^{n+1}  q_h^f, \quad
J_{f^{n+1}}(s_{1,h}) : = \int_{\Omega} f^{n+1} s_{1,h},
\quad J_{g^{n+1}}(s_{2,h}) : = \int_{\Omega} g^{n+1} s_{2,h}\cblue{,}
\nonumber
\end{gather}
and where for the treatment of the convection-diffusion-reaction problem, we have proceeded as in \cite{murphy19}. The only nonlinearities
reside in the terms $ F_{r_h^{n+1}}(\bv_h^s),\;J_{f_h^{n+1}}(s_{1,h})$, and $J_{g_h^{n+1}}(s_{2,h}).$

With the aim to rewrite the Galerkin scheme \eqref{weak-u-h}-\eqref{weak-w2-h} as a matrix equation, we write the unknowns $\bu^s_h, \psi_h, p_h, w_{1,h}$ and $w_{2,h}$ in terms of their  basis functions:
$$\bu^s_h=\sum_{j=1}^{N_1}{U}_j\boldsymbol{\varphi}_j,\quad \psi_h=\sum_{j=1}^{N_2}\Phi_j{\widehat{\varphi}}_j,\quad p_h=\sum_{j=1}^{N_3}P_j{\widetilde{\varphi}}_j,\quad w_{1,h}=\sum_{j=1}^{N_4}W_{1,j}{\overline{\varphi}}_j,\quad
w_{2,h}=\sum_{j=1}^{N_4}W_{2,j}{\overline{\varphi}}_j,$$
and substituting back into \eqref{weak-u-h}-\eqref{weak-w2-h} we obtain 
\begin{align*}
&\sum_{j=1}^{N_1}(\boldsymbol{\varphi}_j,\boldsymbol{\varphi}_i)U_j^{n+1}+   \Delta t^2\sum_{j=1}^{N_1}a_1(\boldsymbol{\varphi}_j,\boldsymbol{\varphi}_i)U_j^{n+1}                   \;+\; \Delta t^2 \sum_{j=1}^{N_2}b_1(\boldsymbol{\varphi}_i,\widehat{\varphi}_j)\Phi_j^{n+1}     =\Delta t^2\sum_{j=1}^{N_1}F_{r_h^{n+1}}(\boldsymbol{\varphi}_i)\\
&\qquad\qquad\qquad\qquad\qquad\qquad+\sum_{j=1}^{N_1}(\boldsymbol{\varphi}_j,\boldsymbol{\varphi}_i)U_j^{n}-\sum_{j=1}^{N_1}(\boldsymbol{\varphi}_j,\boldsymbol{\varphi}_i)U_j^{n-1}\quad i=1,...,N_1, \\
&\biggl(c_0+\frac{\alpha^2}{\lambda}\biggr)\sum_{j=1}^{N_3}(\widetilde{\varphi}_j,\widetilde{\varphi}_i)P_j^{n+1}+\Delta t\sum_{j=1}^{N_3}{a}_2(\widetilde{\varphi}_j,\widetilde{\varphi}_i)P_j^{n+1}                    -\sum_{j=1}^{N_2}{b}_2(\widetilde{\varphi}_i,\widehat{\varphi}_j)\Phi_j^{n+1} =\Delta t \sum_{j=1}^{N_3} G(\widetilde{\varphi}_i)\\
&\qquad\qquad\qquad\qquad\qquad\qquad+\biggl(c_0+\frac{\alpha^2}{\lambda}\biggr)\sum_{j=1}^{N_3}(\widetilde{\varphi}_j,\widetilde{\varphi}_i)P_j^{n}-\sum_{j=1}^{N_2}{b}_2(\widetilde{\varphi}_i,\widehat{\varphi}_j)\Phi_j^{n}\quad i=1,...,N_3, \\
&\sum_{j=1}^{N_1}b_1(\boldsymbol{\varphi}_j,\widehat{\varphi}_i)U_j^{n+1}+\sum_{j=1}^{N_3}b_2(\widetilde{\varphi}_j,\widehat{\varphi}_i)P_j^{n+1} -\sum_{j=1}^{N_2} a_3(\widehat{\varphi}_j,\widehat{\varphi}_i)\Phi_j^{n+1} = 0 \quad i=1,...,N_2, \\
&\sum_{j=1}^{N_4}(\overline{\varphi}_j,\overline{\varphi}_i)W_j^{n+1}-\Delta t\sum_{j=1}^{N_4}{a}_4(\overline{\varphi}_j,\overline{\varphi}_i)W_j^{n+1}= \; \Delta t\sum_{j=1}^{N_4}J_{f_h^{n+1}}(\overline{\varphi}_i) +\sum_{j=1}^{N_4}(\overline{\varphi}_j,\overline{\varphi}_i)W_j^{n} \quad i=1,...,N_4,\\
&\sum_{j=1}^{N_4}(\overline{\varphi}_j,\overline{\varphi}_i)W_j^{n+1}-\Delta t\sum_{j=1}^{N_4}{a}_5(\overline{\varphi}_j,\overline{\varphi}_i)W_j^{n+1}= \; \Delta t\sum_{j=1}^{N_4}J_{g_h^{n+1}}(\overline{\varphi}_i) +\sum_{j=1}^{N_4}(\overline{\varphi}_j,\overline{\varphi}_i)W_j^{n} \quad i=1,...,N_4.
\end{align*}
Then, we can organise the above system in terms of matrices and vectors: 
\begin{gather*}
\widetilde{A}_1\in \mathbf{R}^{N_1\times N_1},\quad 
 {A}_1\in \mathbf{R}^{N_1\times N_1}, \quad 
 {B}_1\in \mathbf{R}^{N_2\times N_1}, \quad 
 {B}_2\in \mathbf{R}^{N_3\times N_2}, \quad 
 \widetilde{A}_2\in \mathbf{R}^{N_3\times N_3},\quad {A}_2\in \mathbf{R}^{N_3\times N_3}, \\
 {A}_3\in \mathbf{R}^{N_2\times N_2}, \quad 
 \widetilde{A}_4\in \mathbf{R}^{N_4\times N_4}, \quad
  {A}_4\in \mathbf{R}^{N_4\times N_4}, \quad  
 \widetilde{A}_5\in \mathbf{R}^{N_4\times N_4},\quad {A}_5\in \mathbf{R}^{N_4\times N_4}, \\
 F\in \mathbf{R}^{N_1}, \quad 
 G\in \mathbf{R}^{N_3}, \quad 
 J_1\in \mathbf{R}^{N_4},\quad 
 J_2\in \mathbf{R}^{N_4},\end{gather*}
 such that 
$$\widetilde{a}_{1,ij}=(\boldsymbol{\varphi}_j,\boldsymbol{\varphi}_i),\qquad {a}_{1,ij}=\Delta t^2 a_1(\boldsymbol{\varphi}_j,\boldsymbol{\varphi}_i)\quad i,j=1,...,N_1,\qquad a_{3,ij}=a_3(\widehat{\varphi}_j,\widehat{\varphi}_i), \quad i,j=1,...,N_2,$$
$$\widetilde{a}_{2,ij}=\left(c_0+\frac{\alpha}{\lambda}\right)(\widetilde{\varphi}_j,\widetilde{\varphi}_i),\qquad {a}_{2,ij}=\Delta t\,a_2(\widetilde{\varphi}_j,\widetilde{\varphi}_i)\quad i,j=1,...,N_3,$$
$$ b_{1,ij}=b_1(\boldsymbol{\varphi}_j,\widehat{\varphi}_i)\quad i=1,...,N_2,\quad j=1,...,N_1,\qquad b_{2,ij}=b_2(\widetilde{\varphi},\widehat{\varphi}_i)\quad i=1,...,N_2,\quad j=1,...,N_3,$$
$$\widetilde{a}_{4,ij}=(\overline{\varphi}_j,\overline{\varphi}_i),\qquad {a}_{4,ij}=\Delta t\,a_4(\overline{\varphi}_j,\overline{\varphi}_i)\quad i,j=1,...,N_4,$$
$$\widetilde{a}_{5,ij}=(\overline{\varphi}_j,\overline{\varphi}_i),\qquad {a}_{5,ij}=\Delta t\, a_5(\overline{\varphi}_j,\overline{\varphi}_i)\quad i,j=1,...,N_4,\qquad F_i=\Delta t^2F_{r^{n+1}}(\boldsymbol{\varphi}_i)\quad i=1,...,N_1,$$
$$G_i=\Delta tG(\widetilde{\varphi}_i)\quad i=1,...,N_3,\qquad J_{1,i}=\Delta t J_{f_h^{n+1}}(\overline{\varphi}_i),\qquad J_{2,i}=\Delta tJ_{f_h^{n+1}}(\overline{\varphi}_i)\quad i=1,...,N_4,$$ and then, denoting 
$$
\mathbf{A}:=\left[\begin{array}{ccccc}
\widetilde{A}_1+A_1&\Delta t^2{B}^{\mathrm{T}}_1&\mathrm{O}&F\\
\mathrm{O}&\widetilde{A}_2+A_2&-{B}^{\mathrm{T}}_2&\mathrm{O}\\
{B}_1&{B}_2&-A_3&\mathrm{O}&\\
\mathrm{O}&\mathrm{O}&\mathrm{O}&\widetilde{A}_4+A_4+\widetilde{A}_5+A_5
\end{array}\right],\quad \mathbf{X}:=
\left[\begin{array}{c}
U_j\\
P_j\\
\Phi_j\\
(W_{1,j},W_{2,j})
\end{array}\right]$$
$$\mathbf{H}:=\left[\begin{array}{c}
2\widetilde{A}_1U_j^n-\widetilde{A}_1U_j^{n-1}\\
{G}+\widetilde{A}_2P^n_j-B_2\Phi^n_j\\
\mathrm{O}\\
J_1+\widetilde{A}_4W^{n}_{1,j}+J_2+\widetilde{A}_5W^{n}_{1,j}\\
\end{array}\right],
$$
the fully-discrete matrix problem for \eqref{eq:forms} reads   $$\mathbf{A}\mathbf{X}^{n+1}=\mathbf{H}^n,$$
which will be used for the development of the numerical tests.

In addition, a Newton method with exact Jacobian is derived for the solution of \eqref{weak-u-h}-\eqref{weak-w2-h} at each time step. Then, regarding both chemical species in a single vector $\bW$,
the tangent algebraic systems to be solved at each Newton step (for a given time step) adopt the following form
\begin{alignat*}{5}
&\; \widehat{A}_{11} \delta\bU^{s}_{k+1}  && &\;+&\; \widehat{B}_{13}\delta\bPsi_{k+1} &\;+ \widehat{F}_1\delta\bW_{k+1}  & = &\; \bR_{1,k},\\
 & && \widehat{A}_{22} \delta\bP^{f}_{k+1}  &\;-&\; \widehat{B}_{23}\delta\bPsi_{k+1} & & = &\; \bR_{2,k},\\
 & \widehat{B}_{13}^\intercal \delta\bU^s_{k+1} &\;+\;& \widehat{B}_{32} \delta\bP^{f}_{k+1} &\;-&\; \widehat{A}_{33}\delta\bPsi_{k+1} && = &\; \bR_{3,k},\\
\;& \widehat{C}_{42}^k\delta\bU^{s}_{k+1} && &&&\;+ (\widehat{J}_4 + \widehat{C}_{41}^k + \widehat{A}_{41}) \delta\bW_{k+1}  & =&\; \bR_{4,k},
\end{alignat*}
where $\delta(\cdot)_{k+1}$ represent the vector of nodal values for the incremental unknowns that are premultiplied by the respective
elementary matrices constructed with the bilinear forms in \eqref{eq:forms} or their linearisation; that is, the matrix $\widehat{F}_1$ is induced by the linearisation of $F_{r_h^{n+1}}(\cdot)$,
$\widehat{A}_{11}$ by $a_1(\cdot,\cdot)+\tilde{a}_1(\cdot,\cdot)$, $\widehat{B}_{13}$ and $\widehat{B}_{13}^\intercal$ by $b_1(\cdot,\cdot)$, $\widehat{A}_{22}$ by $a_2(\cdot,\cdot)+\tilde{a}_2(\cdot,\cdot)$, $\widehat{B}_{23}$ by ${b}_2(\cdot,\delta_t(\cdot))$,
$\widehat{B}_{32}$ by $b_2(\cdot,\cdot)$, $\widehat{A}_{33}$ by $a_3(\cdot,\cdot)$, $\widehat{J}_4$ by  the linearisation of $J_{f_h^{n+1}}(\cdot)$ and $J_{g_h^{n+1}}(\cdot)$,
$\widehat{C}_{41}^k$ and $\widehat{C}_{42}^k$ by the linearisation of $c(\cdot,\cdot,\cdot)$ (see its definition in \cite[eq. (2.7)]{verma01}), and $\widehat{A}_{41}$ by $\widetilde{a}_4(\cdot,\cdot)+a_4(\cdot,\cdot)$ and $\widetilde{a}_5(\cdot,\cdot)+a_5(\cdot,\cdot)$.
The right-hand side vectors $\bR_{i,k}$ account for body forces, mass sources, \cblue{terms associated with}  the previous time step, and
residuals from the previous Newton iteration $k$. The system is solved by the GMRES Krylov solver with
incomplete LU factorisation (ILUT) preconditioning. The stopping criterion on the nonlinear iterations is based on a weighted residual norm dropping below the fixed
 tolerance of $1\cdot 10^{-6}$.


	\section{Numerical tests}\label{sec:results}

	\begin{figure}[!t]
		\begin{center}
			\includegraphics[width=0.325\textwidth]{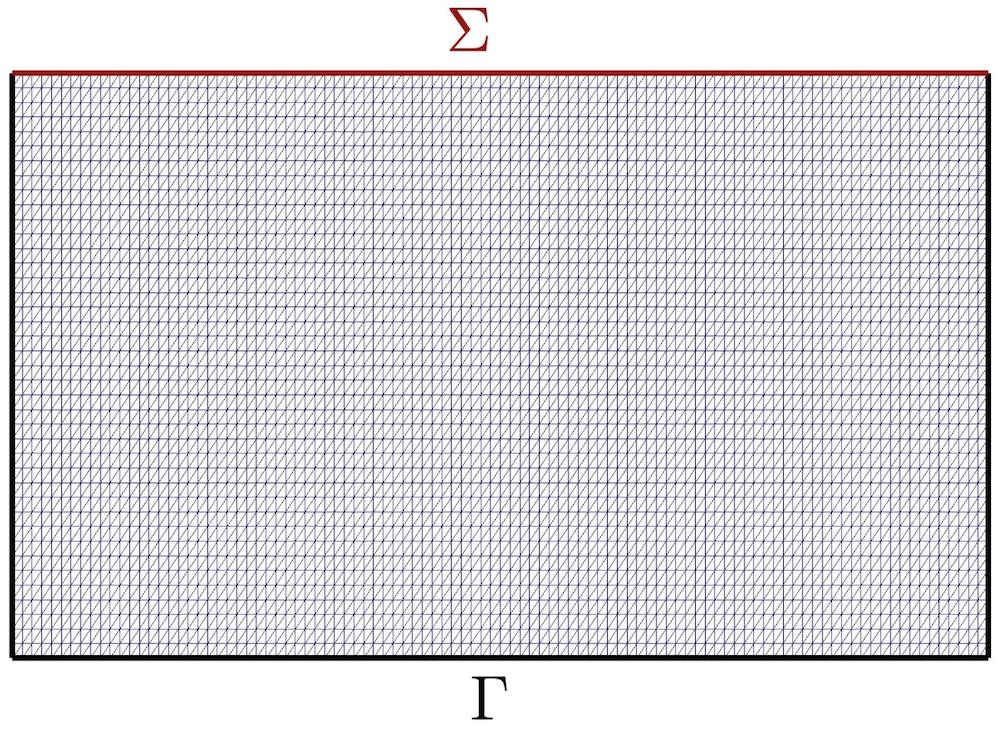}\\
			\includegraphics[width=0.325\textwidth]{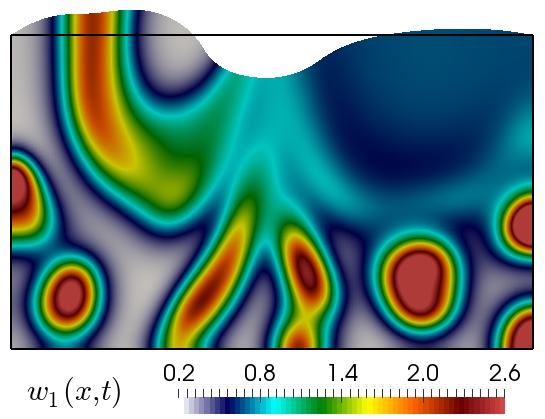}
			\includegraphics[width=0.325\textwidth]{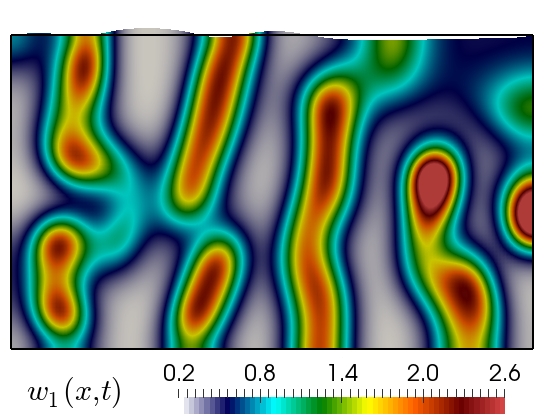}
			\includegraphics[width=0.325\textwidth]{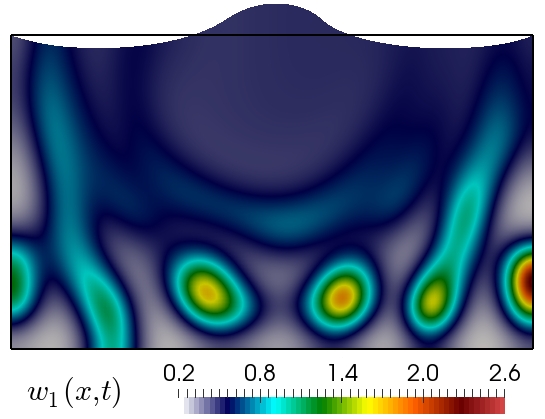}\\
			\includegraphics[width=0.325\textwidth]{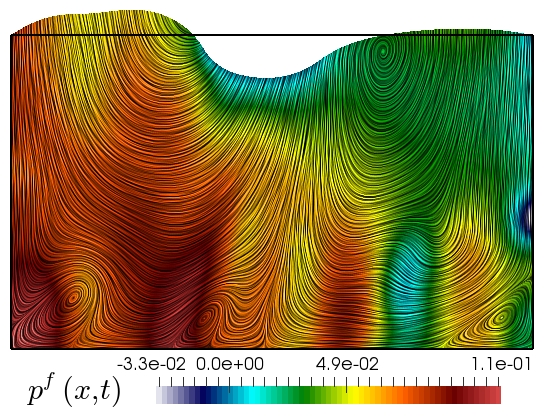}			
			\includegraphics[width=0.325\textwidth]{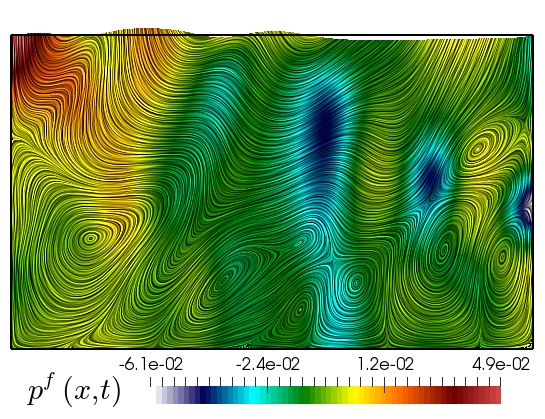}			
			\includegraphics[width=0.325\textwidth]{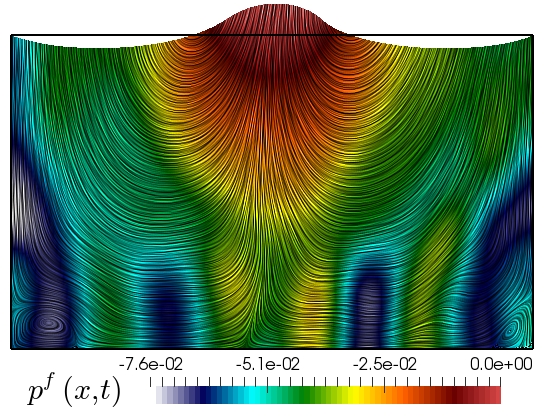}
		\end{center}
		\vspace{-3mm}
		\caption{Test 1. \cblue{Sketched mesh and domain before deformation (top) and sample} of concentrations of $w_1$ \cblue{(middle)} and fluid pressure (bottom)
			at times $t=0.5,1,1.75$ and plotted on the deformed domain according to the solid displacements.}\label{fig:ex03}
	\end{figure}

	\subsection{Test 1: Periodic traction preventing stable patterning}
	In order to investigate the impact that the structural deformation has on the emerging spatial patterns of solutes, we first consider the spatial domain $\Omega=(0,1)\times(0,0.6)$, where
	the clamped boundary is $\Gamma=\{\bx: x_1=0,x_1=1,x_2=0\}$ and the top face constitutes
	$\Sigma$ where we apply a periodic traction defined by
	$$\bt = \begin{cases}
	(0,-s_0 \sin(\pi t))^\intercal & \text{if} \quad 0.4 \leq x_1 \leq 0.6,\\
	\cero & \text{otherwise},\\
	\end{cases}$$
	with $s_0=25000$ (similarly as in the footing problem from \textit{e.g.}  \cite{oyarzua16}). \cblue{A schematic description of the domain can be seen in Figure~\ref{fig:ex03} (top)}. According to \eqref{bc:Gamma}-\eqref{bc:Sigma},
	on $\Gamma$ we also impose zero fluid pressure fluxes, whereas on $\Sigma$
	we set a uniform fluid pressure $p^f = 0$. The parameters that are modified with respect to Test 2 are only
	the coupling constants of active stress modulation $\tau=100$ (using again $r =w_1+w_2$), the
	direction $\bk =(1,0)^\intercal$, the density $\rho = 1$, 
	and the volume-dependent source $\gamma=0.05$.  The resulting patterns (exemplified by
	transients of the activator chemical $w_1$ and final states of poromechanical variables) are
	depicted in Figure~\ref{fig:ex03}. One can readily observe that, apart from altering substantially the distribution of chemical concentrations from
	the beginning of the simulation, the periodic traction applied on part of the top edge (and which only produces less than a
	10\% of vertical stretch) prevents the system from reaching a state with stable spatial patterns.
	For this test we have used a \cblue{uniform} mesh \cblue{(top panel of Figure~\ref{fig:ex03})}.
	\begin{figure}[!t]
\begin{center}
\includegraphics[width=0.24\textwidth]{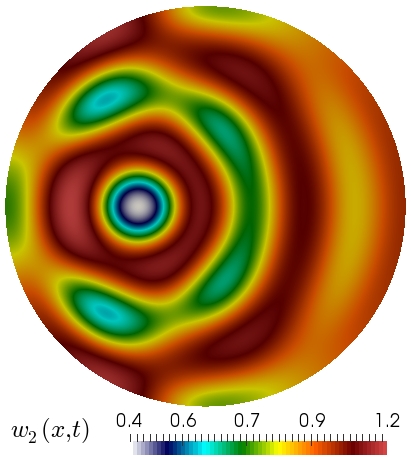}
\includegraphics[width=0.24\textwidth]{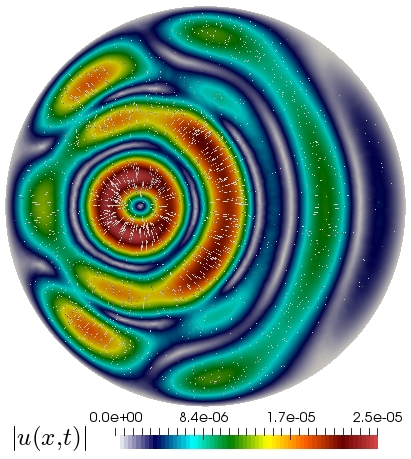}
\includegraphics[width=0.24\textwidth]{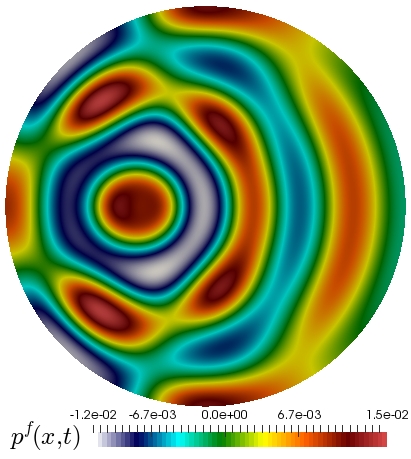}
\includegraphics[width=0.24\textwidth]{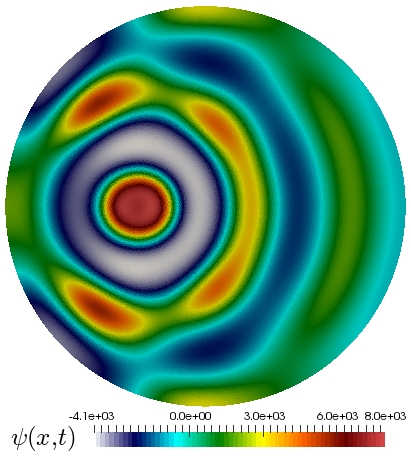}\\
\includegraphics[width=0.24\textwidth]{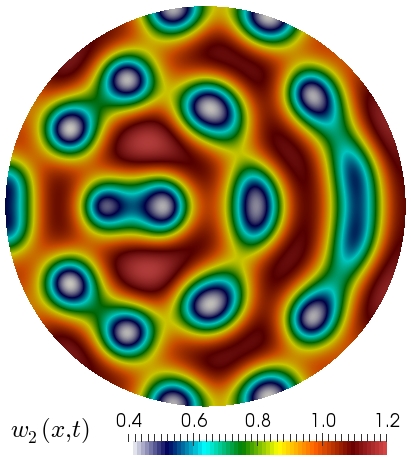}
\includegraphics[width=0.24\textwidth]{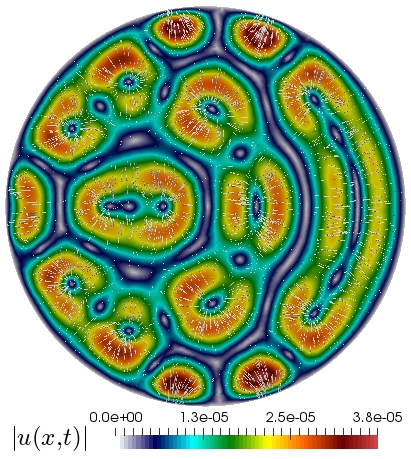}
\includegraphics[width=0.24\textwidth]{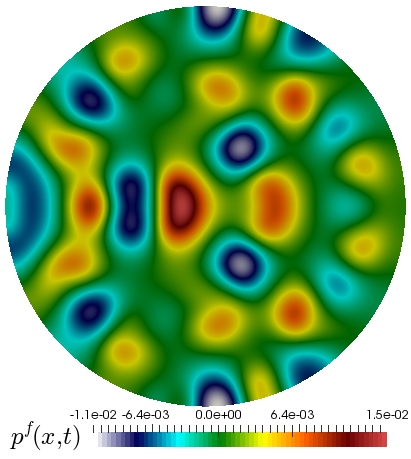}
\includegraphics[width=0.24\textwidth]{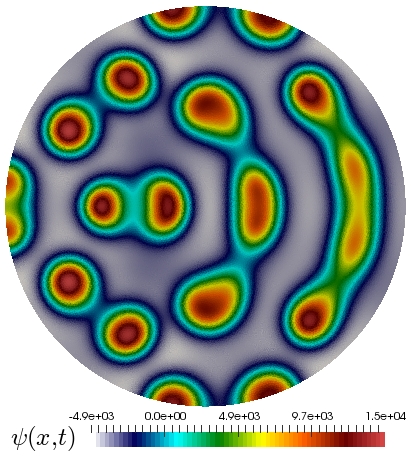}\\
\includegraphics[width=0.24\textwidth]{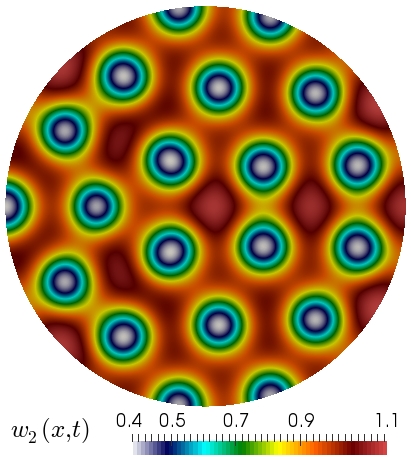}
\includegraphics[width=0.24\textwidth]{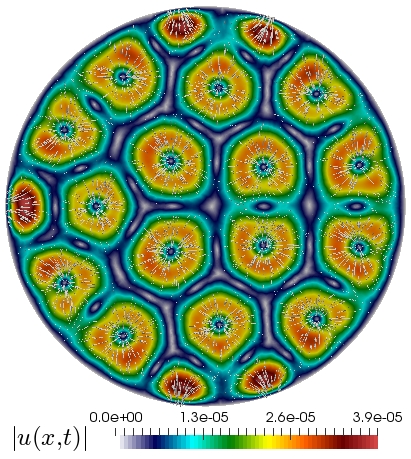}
\includegraphics[width=0.24\textwidth]{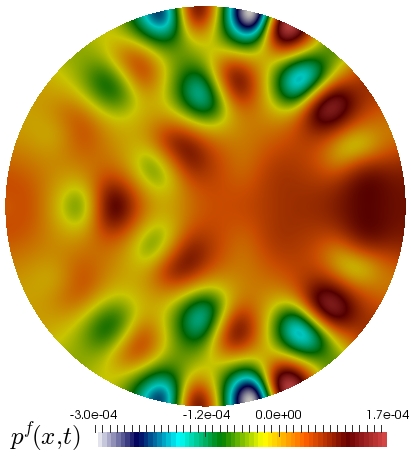}
\includegraphics[width=0.24\textwidth]{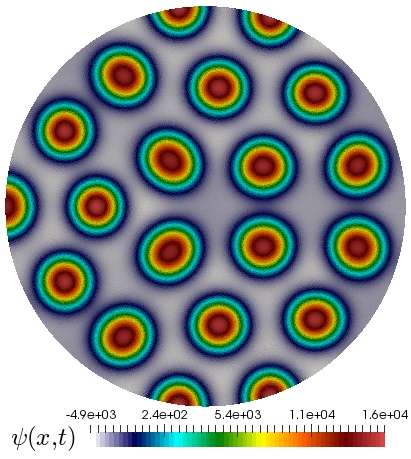}
\end{center}
\vspace{-3mm}
\caption{Test 2. Patterns generated with relatively small chemo-mechanical feedback $\gamma=0.0001$ and
clamped boundary conditions, plotted on the
undeformed domain at times $t=0.25$ (top), $t=0.375$ (middle row), 
and $t=1.5$ (bottom panels).
From left to right: $w_2$ concentration, solid displacement magnitude, fluid pressure, and total pressure.}\label{fig:ex02}
\end{figure}

\subsection{Test 2: Small poromechanical effects}
Next we take the domain as the disk centred at $(0.5,0.5)$ with radius 0.5, and assume that
the boundary coincides with $\Gamma$. Then the displacements are set to zero on the whole
boundary and we take relatively small values for the coupling constants on the chemical
source and on the active stress $\gamma=0.0001$, $\tau = 100$, implying in particular
that the patterns produced by the coupled poroelastic-convection-diffusion-reaction system are
expected to be qualitatively similar to those observed on a fixed domain. In \eqref{eq:active-stress}
we assume a dependence of the form $r=w_1+w_2$, and the remaining model constants
are taken as $D_1=0.05$, $D_2=1$, $\beta_1=170$, $\beta_2=0.1305$, $\beta_3=0.7695$, $E=3\cdot10^4$,
$\nu = 0.495$, $\rho=1$, $c_0=1\cdot 10^{-3}$, $\kappa = 1\cdot 10^{-4}$, $\alpha=0.1$, $\eta=1$,
$\bk =(x_1-0.5,x_2-0.5)^\intercal$. The initial condition for the chemicals is a perturbation of the homogeneous steady state
$w_1^0=\beta_2+\beta_3$, $w_2^0=\beta_3(\beta_2+\beta_3)^{-2}$ and for the displacements and fluid
pressure we use zero initial conditions. The domain is
discretised into an unstructured mesh of 64926 triangles and we employ a fixed time-step
$\Delta t = 0.0025$.
The system is advanced until $t_{\text{final}}=1.5$ and plots with patterns of
$w_2$, small deformations, as well as fluid and total pressures are shown in Figure~\ref{fig:ex02}.
In the bottom row we can see how the initial perturbation of the steady state evolves into
organised dot-shaped spatial structures, seen clearly for the inhibitor chemical $w_2$ and also captured by the
total pressure. No deformation occurs along the domain boundary, but the local deformation patterns
show also tissue contraction near the zones of high concentration of the activator species $w_1$.
\begin{figure}[!t]
\begin{center}
\includegraphics[width=0.24\textwidth]{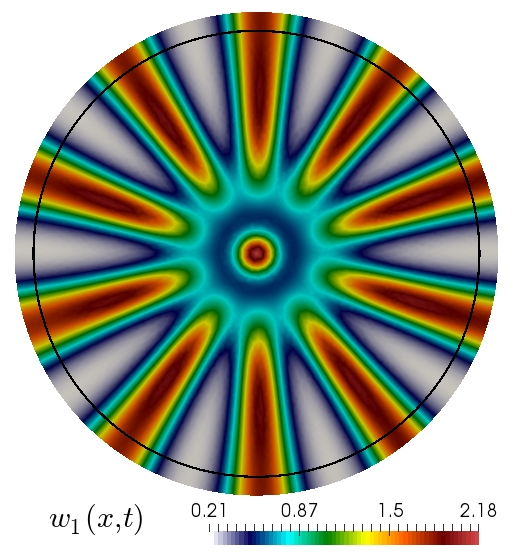}
\includegraphics[width=0.24\textwidth]{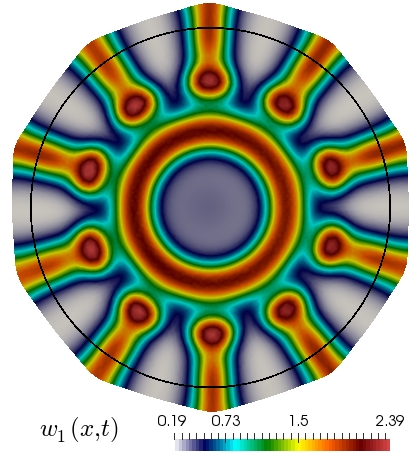}
\includegraphics[width=0.24\textwidth]{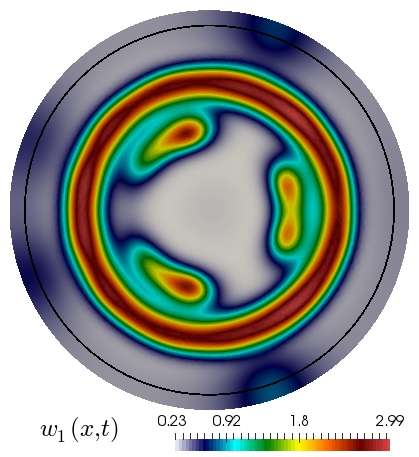}
\includegraphics[width=0.24\textwidth]{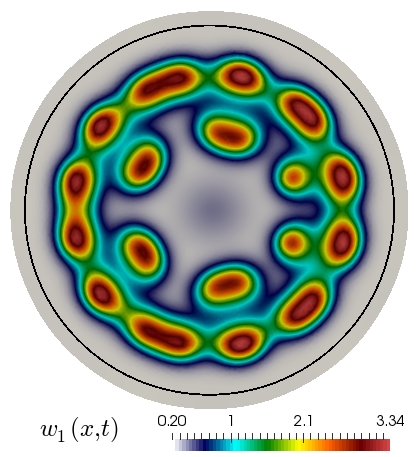}\\
\includegraphics[width=0.24\textwidth]{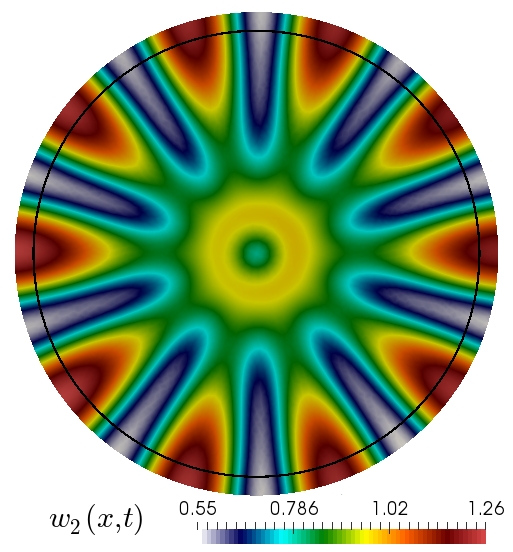}
\includegraphics[width=0.24\textwidth]{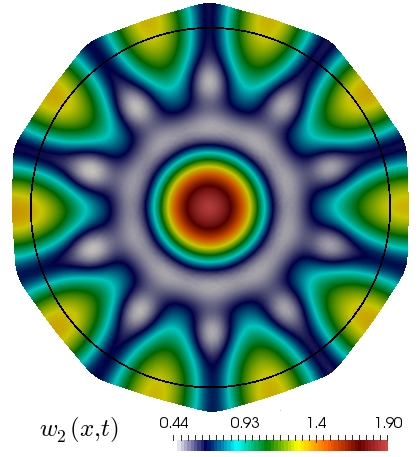}
\includegraphics[width=0.24\textwidth]{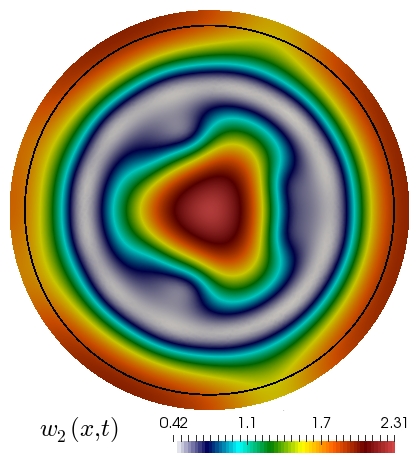}
\includegraphics[width=0.24\textwidth]{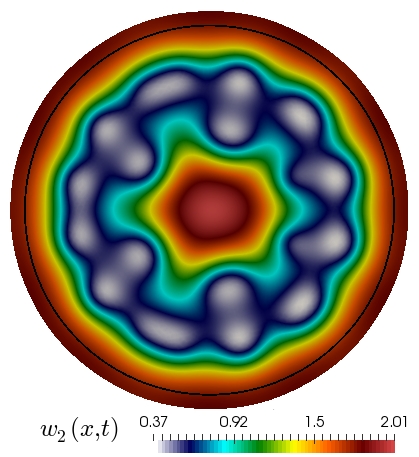}
\end{center}
\vspace{-3mm}
\caption{Test 3. Interplay between linear growth and active stress.
Concentrations of $w_1$ and $w_2$ (top and bottom) on the deformed domain, at time $t=1.25$, for
different mild values of the coupling constants $\tau,\tau_2$.
The black circle in all plots represents the boundary of the initial domain. }\label{fig:ex042d}
\end{figure}

\subsection{Test 3: Linear growth and active stress}
Maintaining the same domain and discretisation parameters \cblue{as in Test 2 above},
we now fix the mechano-chemical coupling constant $\gamma=0.01$ and
study the competing effect between linear growth with radial traction, and
the active stress depending on the concentration of the activator species $w_1$. This is
done with an activation of the type $r = \tau_2 t + w_1^2$, and for this 2D case we consider Robin boundary conditions for the solid motion 
\eqref{eq:Robin} with $\zeta = t \tau_2$, on the circular boundary (whereas for the 3D case below 
we impose zero normal displacement on the bottom of the cylinder and a traction $\bt = \tau_2 t \nn$ 
on the remainder of the boundary).
We set $\bk$ to be the radial vector, and vary $\tau,\tau_2$.  The results are
shown in Figure~\ref{fig:ex042d}. From left to right we display snapshots of the
chemical patterns produced with the parameter choices
 $(\tau  =2\cdot10^5, \tau_2 = 0.2)$,   $(\tau = 10^5,\tau_2 = 2)$,    $(\tau = 10^4,\tau_2 = 10)$,
and $(\tau = 100,\tau_2 = 20)$.

	\begin{figure}[!t]
		\begin{center}
			\includegraphics[width=0.26\textwidth]{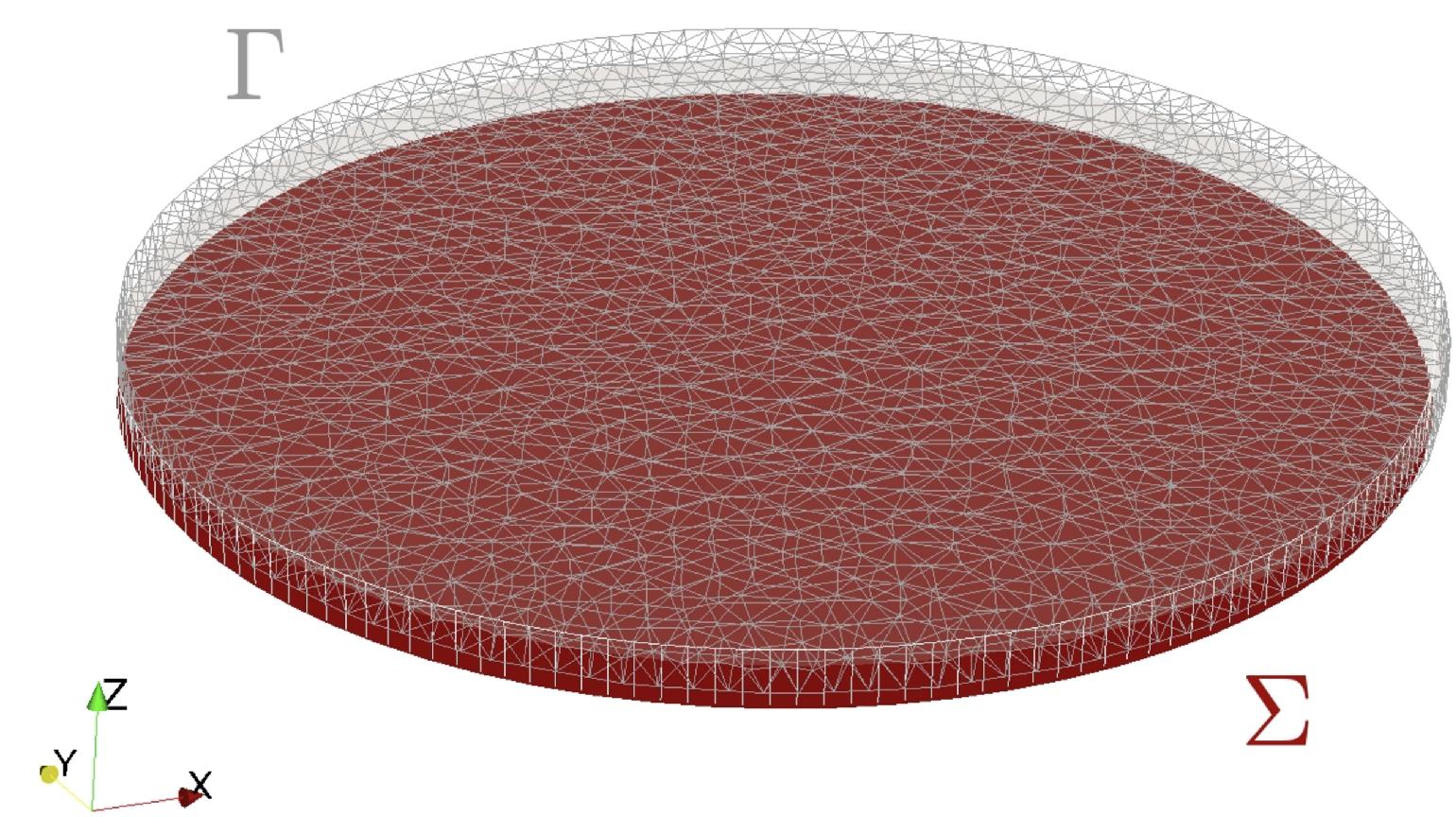}\\
			\includegraphics[width=0.325\textwidth]{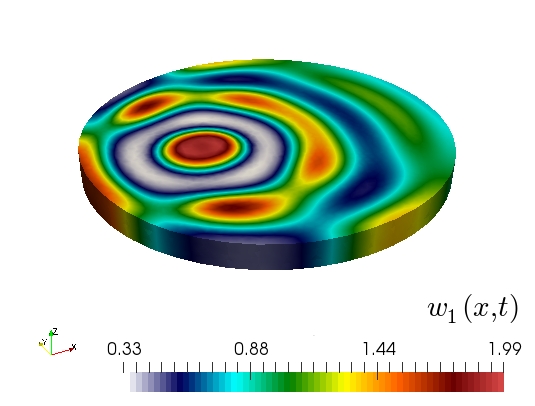}
			\includegraphics[width=0.325\textwidth]{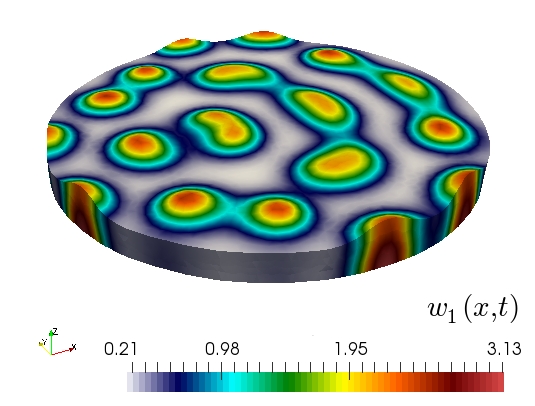}
			\includegraphics[width=0.325\textwidth]{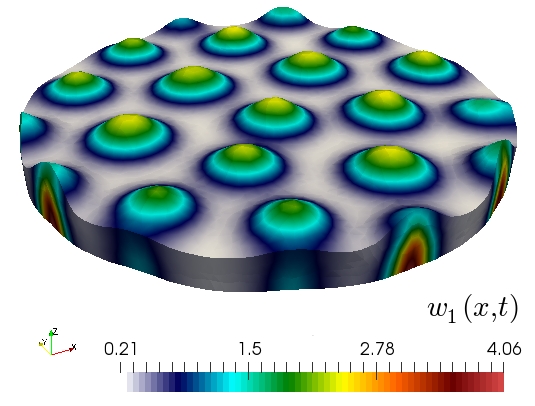}\\
			\includegraphics[width=0.325\textwidth]{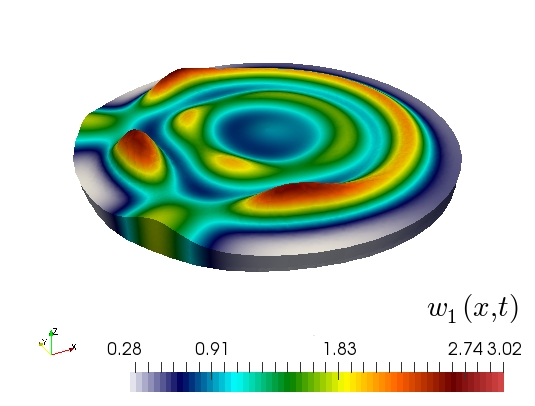}			
			\includegraphics[width=0.325\textwidth]{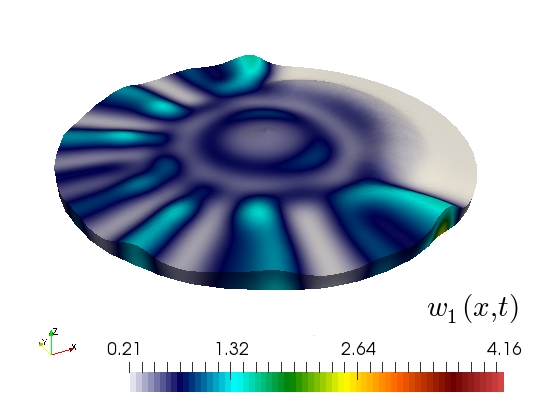}			
			\includegraphics[width=0.325\textwidth]{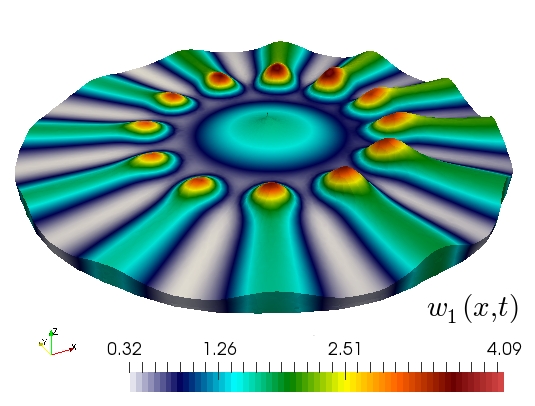}
		\end{center}
		\vspace{-3mm}
		\caption{Test 4. \cblue{Sample coarse mesh and domain/boundary configuration (top), and interplay} between linear growth and active stress for two sets of
			elasticity and coupling parameters \cblue{(middle and bottom)}.
			Concentrations of $w_1$ on the deformed domain, at times $t=0.25, 0.5, 1.25$ (left, centre, and right panels, respectively).}
		\label{fig:ex04}
	\end{figure}

	\subsection{Test 4: Linear growth in 3D}
	We extend the previous test to assess  the behaviour of the model and the finite element scheme
	in a 3D setting. We modify \eqref{eq:active-stress}  to include axial symmetry on another preferential
	direction for active deformation.
	The domain is a cylinder of height 0.05 and radius 0.5 \cblue{(see top panel of Figure~\ref{fig:ex04} showing the different boundaries of the geometry and the generated volumetric mesh)}, and we set
	\begin{equation*}
	\bsigma_{\text{act}} = -\tau \bigl[ \tau_2 t \bk_{12}\otimes\bk_{12} + w_1^2\bk_3\otimes\bk_3 \bigr],
	\end{equation*}
	where $\bk_{12}$ is the radial vector in the plane $x_1x_2$ and $\bk_3=(0,0,1)^\intercal$. This implies that
	the active deformation due to linear growth will occur in the radial direction
	whereas the
	stress due to the chemical concentration will act on the $x_3-$direction. \cblue{On the
	bottom surface (denoted $\Sigma$) we set zero normal displacements $\bu\cdot\nn = 0$, while on the remainder of the
	boundary, which is denoted by $\Gamma$}, we impose zero traction.
	In Figure~\ref{fig:ex04} we show the
	resulting patterns of $w_1$ concentration for two sets of Lam\'e and poromechanical-chemical coupling
	parameters. For the first case we use  $E=1\cdot 10^4$, $\nu=0.499$, and $\tau = 10$, $\tau_2 = 20$, $\gamma=0.05$
	and show the patterns on the deformed domain  in the left panels; while the plots on the right panels were produced with
	$E=1\cdot 10^{3}$, $\nu = 0.3$, and $\tau = 50$, $\tau_2 = 60$, $\gamma=0.1$. We observe stable
	pattern generation with the first set of model parameters, similar to the expected patterning in the case of pure reaction-diffusion
	effects, whereas the patterns on the right exhibit large qualitative differences in $w_1$ (also in the other species)
	as well as in the deformation behaviour.

	\subsection{Test 5: Application to the simulation of brain injuries and calcium propagation}
	We close this section with an example related to the one-way coupling between poroelastic deformations in the brain 
	(induced by a localised high stress) and the subsequent propagation and reaction of two types of calcium concentration, 
	intra-cellular and extra-cellular, throughout the tissue. This illustrative test is based on the kinetic and 1D models 
	recently advanced in  \cite{kant18}. In there, the authors propose that hydrostatic stress build up due to the brain trauma affect (in an exponentially decreasing manner) the reacting fluxes between the calcium concentrations. We do not include acceleration but we propose to incorporate this in \cblue{the coupled} model using a modification 
	of \eqref{eq:ADR1}-\eqref{eq:ADR2} to include a dependence of the reaction terms on the total pressure 
	\begin{align*} 
	f(w_1,w_2,\psi) & =-D_1(w_1-w_2) + \frac{1}{\chi_1}\biggl[-\chi_1+(1+\chi_1)\exp(-k|\psi|)\biggr] \frac{w_2^2}{w_2^2+k_1^2},\\ 
	g(w_1,w_2,\psi) & = -f(w_1,w_2,\psi) + D_2(w_0-w_2) - \frac{1}{\chi_2}\biggl[-\chi_2+(1+\chi_2)\exp(-k|\psi|)\biggr]\frac{w_2}{w_2+k_2},
	\end{align*}
where $w_1,w_2$ represent respectively, the extra-cellular and intra-cellular calcium concentrations (in [mM] units) and the model parameters are as in \cite{kant18} 
\begin{gather*}
D_1=2.94\cdot 10^{-6}\,[\text{1/s}], \quad D_2=3.17\cdot10^{-5}\,[\text{1/s}], \quad k_1=2\cdot10^{-4}\,[\text{mM}], \quad k_2=5\cdot10^{-4}\,[\text{mM}],\\ 
\chi_1 = 2\cdot10^3, \quad \chi_2=4\cdot10^3,\quad k=4.5\cdot10^{-5}\,[\text{Pa}^{-1}], \quad w_0 = 0.1\,[\text{mM}]. 
\end{gather*} 
	On the other hand, 
	the fact that calcium activity effects are negligible in producing deformations of the poroelastic structure (at least, when compared 
	to high stress impacts on the skull or with important kinematic forces building up because of rapid shocks) implies that in \cblue{the proposed}  
	model the total stress \eqref{eq:total-stress} does not contain an active component modulated by $w_1,w_2$. Also, the present model is different than the one in \cite{kant18} in that we do not consider viscoelastic effects but do include poroelasticity of the brain, and we also include diffusion of the calcium concentrations. 
	The remaining constants in the model and  the initial conditions adopt the values 
	\begin{gather*}
	E = 3.15\cdot10^{4}\,[\text{Pa}], \quad \nu = 0.45, \quad \rho = 1130\,[\text{Kg/m}^3], \quad \frac{\kappa}{\eta} = 10^{-5}\,[\text{mm}^2\text{Pa}^{-1}\text{s}^{-1}], \quad \alpha =0.1, \\ c_0 = 3.9\cdot10^{-4}\,[\text{Pa}^{-1}], \quad \bu^s(0)=\cero, \quad p^f=0, \quad \psi(0)=0, \quad w_{1,0} = 1\,[\text{mM}], \quad w_{2,0} = 10^{-4}\,[\text{mM}]. 
	\end{gather*}

\begin{figure}[!t]
		\begin{center}
			\raisebox{2mm}{\includegraphics[width=0.325\textwidth]{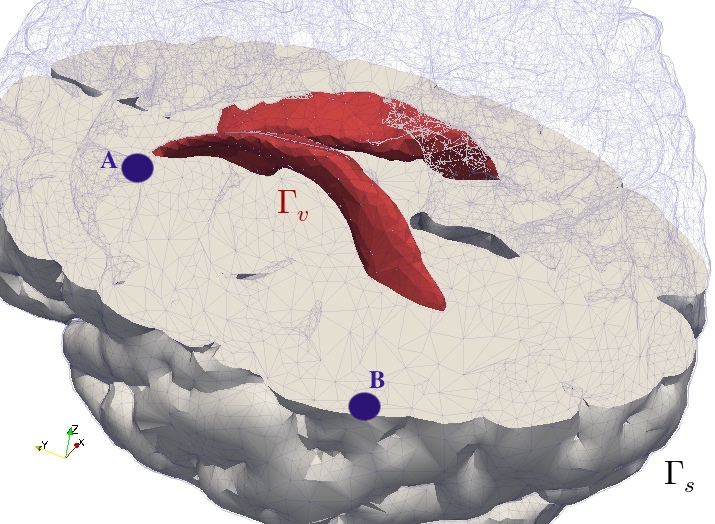}}
			\includegraphics[width=0.325\textwidth]{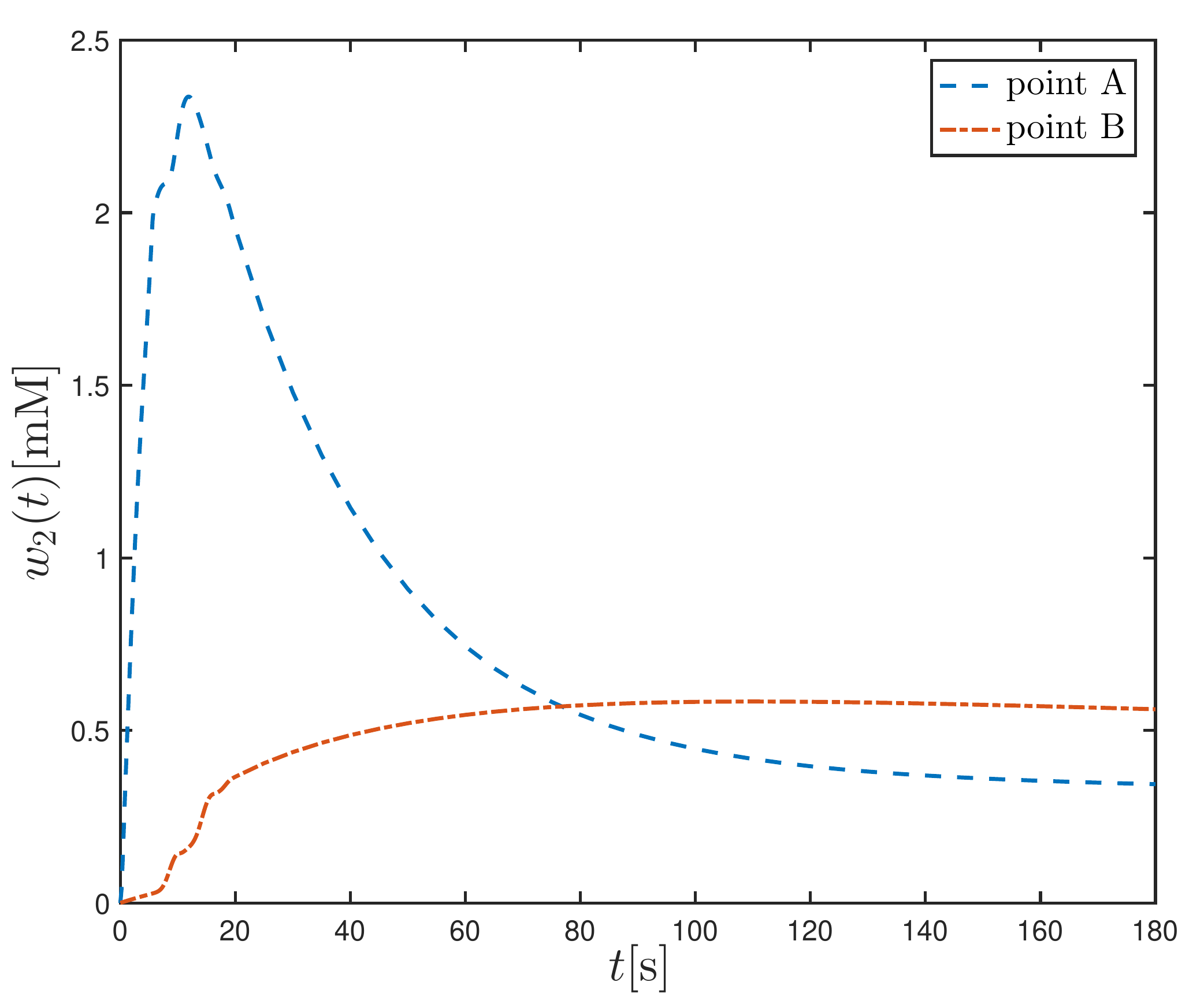}
			\includegraphics[width=0.325\textwidth]{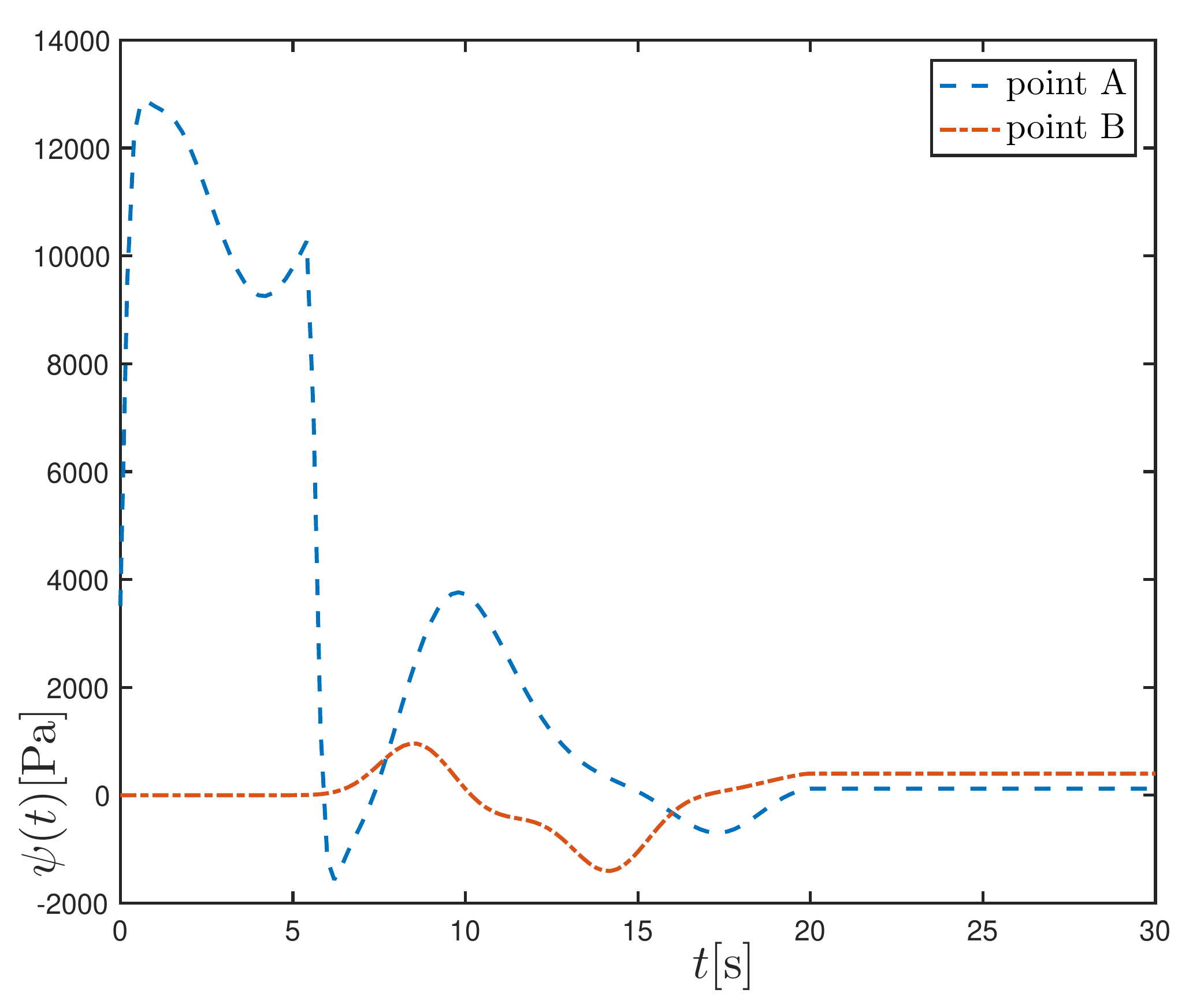}
		\end{center}
		\vspace{-3mm}
		\caption{Test 5. Cut of the human brain geometry and tetrahedral mesh, showing boundaries on the ventricles $\Gamma_v$ and 
		near the skull $\Gamma_s$, also indicating points A and B where we record quantities of interest (left panel). Time evolution of 
			intra-cellular calcium concentration and of total pressure (centre and right plots).}
		\label{fig:ex05}
	\end{figure}

The spatial domain consists of a 3D structure of the human brain and the boundaries are split between ventricles and the outer meningial region of the 
brain, in contact with the  skull. The interstitial 
flow in this case is by cerebrospinal fluid. The domain consists of an adult  human brain atlas \cite{fang09} and we use a tetrahedral mesh 
with 29037 vertices.
	An initial traction of magnitude 1.7$\cdot10^4$\,[Pa] is applied for 5.5\,[ms] on a location near the ventricles and on the skull we impose zero normal displacements and zero 
	fluid pressure, whereas on the ventricles we assume zero fluid pressure flux. We employ a timestep of $\Delta t = 0.1$\,[ms] and 
	run the simulation until $t=180$\,[s]. Transients of the intracellular calcium concentration as well as the total pressure are recorded on two points (one near the ventricles, point A, and another near the meninges, point B), and are displayed in Figure~\ref{fig:ex05}. One can observe an initial peak of several folds the initial homeostatic value of the intracellular calcium followed by a slowly decaying profile (which however does not 
	goes back to the homeostatic value). We also see that the oscillations in total pressure 
	due to the application of high stresses decrease over time. All this 
	is qualitatively consistent with the model predictions from \cite{kant18}.


	\bigskip
	\noindent\textbf{Acknowledgements.} The authors gratefully acknowledge the fruitful discussions with Jenny Dingwall (Oxford), Aayush Kant (Monash), Belle Kim (Oxford), Tim Leach (Oxford), Kent-Andr\'e Mardal (Oslo), 
	and Rodrigo Weber dos Santos (Juiz de Fora),  regarding different aspects of this work.


\section*{References}	
\bibliographystyle{unsrt} 
\bibliography{z-NewBibl}

\end{document}